%% file: manuscript.tex
\documentclass[final,hidelinks,onefignum,onetabnum]{siamart250211}

\input{main_shared}

\ifpdf
\hypersetup{
    pdftitle={Xiao2025RP},
 	pdfauthor={Yao Xiao}
}
\fi

\begin{document}

\maketitle

\begin{abstract}
Lasso regression is a widely employed approach within the $\ell_1$ regularization framework used to promote sparsity and recover piecewise smooth signals $f:[a,b) \rightarrow \mathbb{R}$ when the given observations are obtained from noisy, blurred, and/or incomplete data environments. In choosing the regularizing sparsity-promoting operator, it is assumed that the particular {\em type} of variability of the underlying signal, for example, piecewise constant or piecewise linear behavior across the entire domain, is both {\em known} and {\em fixed}. Such an assumption is problematic in more general cases, e.g.~when a signal exhibits piecewise oscillatory behavior with varying wavelengths and magnitudes. To address the limitations of assuming a fixed (and typically low order) variability when choosing a sparsity-promoting operator, this investigation proposes a novel {\em residual} transform operator that can be used within the Lasso regression formulation.  In a nutshell, the idea is that for a general piecewise smooth signal $f$,  it is possible to design two operators $\mathcal L_1$ and $\mathcal L_2$ such that $\mathcal L_1{\bm f} \approx  \mathcal L_2{\bm f}$, where ${\bm f} \in \R^n$ is a discretized approximation of $f$, but $\mathcal L_1 \not\approx  \mathcal L_2$.  The corresponding {\em residual transform operator}, $\mathcal L  =  \mathcal L_1- \mathcal L_2$, yields a result that (1) effectively reduces the variability dependent error that occurs when  applying either $\mathcal L_1$ or $\mathcal L_2$ to ${\bm f}$, a property that holds even when  $\mathcal L_1{\bm f} \approx \mathcal L_2{\bm f}$ is {\em not} a good approximation to the true sparse domain vector of ${\bm f}$, and (2) does not require $\mathcal L_1$ or $\mathcal L_2$ to have prior information regarding the variability of the underlying signal. Numerical experiments demonstrate the effectiveness of the new residual transform operator when compared to standard sparsity-promoting operators used in Lasso regression for recovering piecewise smooth signals.
\end{abstract}



\section{Introduction} 
\label{sec:introduction}

Recovering piecewise smooth signals from noisy, blurred, and/ or incomplete observations is important across a wide variety of application areas.  Often referred to as compressive sensing (CS) \cite{candes2006robust,candes2006stable,chartrand2008iteratively,donoho2006compressed}, Lasso regression is a widely employed technique within the $\ell_1$ regularization framework that is used to promote sparsity in the recovery of piecewise smooth signals.   

While there are different Lasso regression implementations, each with its own potential benefits and drawbacks, the  sparsity-promoting operator chosen for the corresponding regularization term is almost always based on an assumption regarding the {\em type} of variability in the underlying signal, such as piecewise constant or piecewise linear behavior across the entire domain.  This is problematic in more general cases, for example when the true signal exhibits piecewise oscillatory behavior with varying wavelengths and magnitudes. This assumption is also problematic in applications such as data assimilation for partial differential equations as well as scientific data reduction and compression, where efficient and robust representations of  piecewise smooth highly variable functions are paramount.

To address the limitations of assuming a fixed (and typically low order) variability for the sparsity-promoting operator, this investigation proposes a novel {\em residual} transform operator that can be used within the Lasso regression ($\ell_1$ regularization) framework.  In a nutshell, the idea is that for a general piecewise smooth signal $f:[a,b)\rightarrow \mathbb{R}$ realized at its uniform grid points in its chosen recovery domain, e.g., $\bm f = \{f(s_j)\}_{j = 1}^n$,  it is possible to design two operators $\LL_1$ and $\LL_2$ such that $\LL_1{\bm f} \approx \LL_2{\bm f}$ but $\LL_1 \not\approx \LL_2$.  The corresponding {\em residual transform operator}, $\LL = \LL_1- \LL_2$, yields a result that (1) effectively reduces the variability dependent error that occurs when  applying either $\LL_1$ or $\LL_2$ to ${\bm f}$, a property that holds even when  $\LL_1{\bm f} \approx \LL_2{\bm f}$ is {\em not} a good approximation to the true sparse domain vector of ${\bm f}$, and (2) does not require prior information regarding the variability of the underlying signal when choosing $\LL_1$ or $\LL_2$.  In this investigation we choose $\LL_1$ and $\LL_2$ to be respectively  local and global approximations of the underlying function support in its corresponding edge domain \cite{archibald2005polynomial,gelb2006adaptive}. 
Critically, we observe that while the local and global operators project data differently, they both similarly respond to the spatial variability of the true signal {\em even when that variability is unknown}.  This property is exploited in the development of our new sparsity promoting residual transform operator.
 
It is  important to keep in mind that in and of themselves these local and global edge detection operators can be used as sparsity-promoting transform operators when the variability between discontinuities is {\em known} and {\em of the same type} \cite{archibald2016image,Stefanetal}. For example, a local first order finite differencing operator, or equivalently total variation (TV), is appropriate for recovering piecewise constant signals. However, as already discussed, such fixed assumptions about the signal's variability can degrade the quality of recovery.  On the other hand,  the {\em difference} between their corresponding residuals  to the sparse domain vector of ${\bm f}$  should still be small if  $\LL_1{\bm f} \approx \LL_2{\bm f}$, especially in regions where the true signal is smooth, and provides a compelling argument for choosing $\LL = \LL_1- \LL_2$ as the  sparsity-promoting transform operator in the context of Lasso regression. Our numerical results confirm that employing the residual transform operator within the $\ell_1$ regularization framework yields accurate recovery for piecewise smooth signals {\em even when} neither $\LL_1$ nor $\LL_2$ effectively promote sparsity (for reasons stated above). 

\subsection*{Outline} 
The remainder of this paper is organized as follows: In \Cref{sec:prelim} we introduce the general framework of Lasso regression, highlighting the specific variations used in our approach. The formulations for both the local and global edge detectors are also presented in detail. Our new {\em residual} transform operator $\LL = \LL_1- \LL_2$ is introduced in \Cref{sec:residualprior}.  \Cref{sec:numerics} presents numerical experiments and convergence tests to validate the proposed technique. Finally, some concluding remarks and ideas for future investigations are provided in \Cref{sec:summary}.

\section{Preliminaries} 
\label{sec:prelim} 

Let $f: [-\pi,\pi)\to\R$ be an unknown real-valued $2\pi$-periodic piecewise smooth function.   We seek to recover  $\{f(s_j)\}_{j = 1}^n$  on a fixed set of $n$ uniform grid points $s_j = -\pi + (j-1)\frac{2\pi}{n}$ given observations modeled by   
\begin{equation}
  {\bm y} =  
  A{\bm f} + \boldsymbol\delta.
       \label{eq: data_acquisition}
\end{equation}
Here, ${\bm f} \in \R^n$ is a discretized approximation of $f$ obtained by sampling it at $n$ uniformly spaced points in $[- \pi,\pi)$. In our numerical experiments we will consider the forward operator $A$ to be defined as (1) $A = I_n$ (the $n\times n$ identity matrix), so that \cref{eq: data_acquisition} corresponds to the standard denoising problem,  or (2)  $A = \tilde A \in \R^{n \times n}$, which applies a Gaussian blur defined by its point spread function (psf)
\begin{equation}
 \label{eq:bluroperator}   
\tilde{A}_{ij}= \frac{1}{2\pi\gamma^2}\exp{\left(-\frac{i^2+j^2}{2\gamma^2}\right)}. \end{equation}
Observe that application of \cref{eq:bluroperator} causes the high-frequency details to be ``smoothed out'', making interior structures less easily identifiable while also introducing ill-posedness into the corresponding inverse problem. Finally, we also consider an undersampling problem, where rows of $A = I_n$ are randomly removed at a given ratio, causing undersampling that leads to additional loss of information. In this case, however, we do not assume that we know in advance which rows are removed, so that the forward model used for recovery is still provided by \cref{eq: data_acquisition} for the denoising case.    
The observation $\bm{y}$ is further contaminated by additive, independent, and identically distributed (i.i.d.) Gaussian noise with zero mean and variance $\sigma^2$ given by 
\[\boldsymbol{\delta}\sim\mathcal N\left(0, \sigma^2 I_n\right).\]

\begin{remark}[Assumptions on forward model, $f(s)$, and grid point uniformity]
\label{rem:assumptions}  
    Several assumptions made above are for convenience and not necessary for the development and efficacy of our new residual transform operator.  For example, our method is not limited with regard to the choice of the linear forward operator, with the blurring operator in \cref{eq:bluroperator} chosen as representative of how data may be acquired. A Gaussian matrix, or discrete Fourier or Radon transform could also serve this purpose. Further, $f$ need not be $2\pi$-periodic, although this assumption greatly simplifies the composition of both sparsifying transform operators used in constructing the residual transform operator.  Finally, although not inherently required, having uniform grid points also simplifies the derivation and implementation of each transform operator.  Future investigations will consider underlying functions and data that do not meet these assumptions. 
\end{remark}

\subsection{Lasso Regression}\label{sub:lasso}

Leveraging the sparsity in the sparse transform domain of a piecewise smooth signal naturally leads to the application of CS algorithms \cite{candes2006robust,candes2006stable,candes2006near,donoho2006compressed} for signal recovery. For the forward model  in \eqref{eq: data_acquisition},  a typical CS methodology involves solving the following unconstrained minimization problem, also referred to as Lasso regression:
\begin{equation*}
{\bm x}  =  \argmin_{\tilde{\bm{x}}} \left( \frac{1}{\sigma^{2}} \norm{A\tilde{\bm x}-\bm y}_2 + \tilde\alpha\norm{\mathcal{L}\tilde{\bm x}}_1 \right).
\end{equation*}
The first term, generally known as the fidelity term, measures the discrepancy between the collected data and the forward model. The second component, known as the penalty or regularization term,  integrates prior knowledge about the sparse characteristics of the underlying signal, often in the signal's edge or gradient representation. The operator $\LL$ in the penalty term formalizes such an assumption by presuming that the transformed signal $\LL\bm x$ is sparse.

The Lasso parameter $\tilde\alpha$ (independent of $\sigma^2$) is chosen to weigh the desired influence of the regularization term. A common  equivalent formulation is given by
\begin{equation}
\label{eq:lasso_regression}
{\bm x} =     \argmin_{\tilde{\bm{x}}} \left( \norm{A\tilde{\bm x}-\bm y}_2 + \alpha\norm{\mathcal L \tilde{\bm x}}_1 \right),
\end{equation}
where typically the regularization parameter $\alpha$ is chosen to be directly proportional to the noise variance $\sigma^2$. Here we choose $\alpha$ as proposed in \cite{sanders2020effective},
\begin{equation}\label{eq:lasso_param}
    \alpha=2^{\tfrac32}\frac{\sigma^2}{\beta},
\end{equation}
where the sparse domain average power $\beta^2$ is estimated as
\begin{equation*}
    \beta^2 = \frac{\norm{\mathcal L{\bm x}_{est}}_2^2}{n},
\end{equation*}
and $\bm x_{est}$ is an initial estimate of the solution to \eqref{eq:lasso_regression}. In our numerical experiments we employ 
least squares to obtain $\bm x_{est}$,  although other estimates are also valid, for example, by solving \cref{eq:lasso_regression} with  $\alpha \sim \sigma^2$.   The idea in \cref{eq:lasso_param} is that the weight of the regularization term should also be inversely proportional to the mean of its magnitude.

Due to the non-differentiability of the $\ell_1$-norm, \cref{eq:lasso_regression} does not have a closed-form solution and must be solved iteratively. Numerous efficient algorithms exist for this purpose, such as the alternating direction method of multipliers \cite{wu2010augmented,boyd2011distributed}, the Split-Bregman Algorithm \cite{goldstein2009split}, and the Fast Iterative Shrinkage-Thresholding Algorithm \cite{beck2009fast}.  Since a comparative analysis is outside the scope of this work,  we obtain solutions to \eqref{eq:lasso_regression} using the robust and convenient CVX package for MATLAB \cite{cvx,gb08}.

\begin{figure}[h!]
    \centering
    \begin{subfigure}[b]{.32\textwidth}
\includegraphics[width=\textwidth]{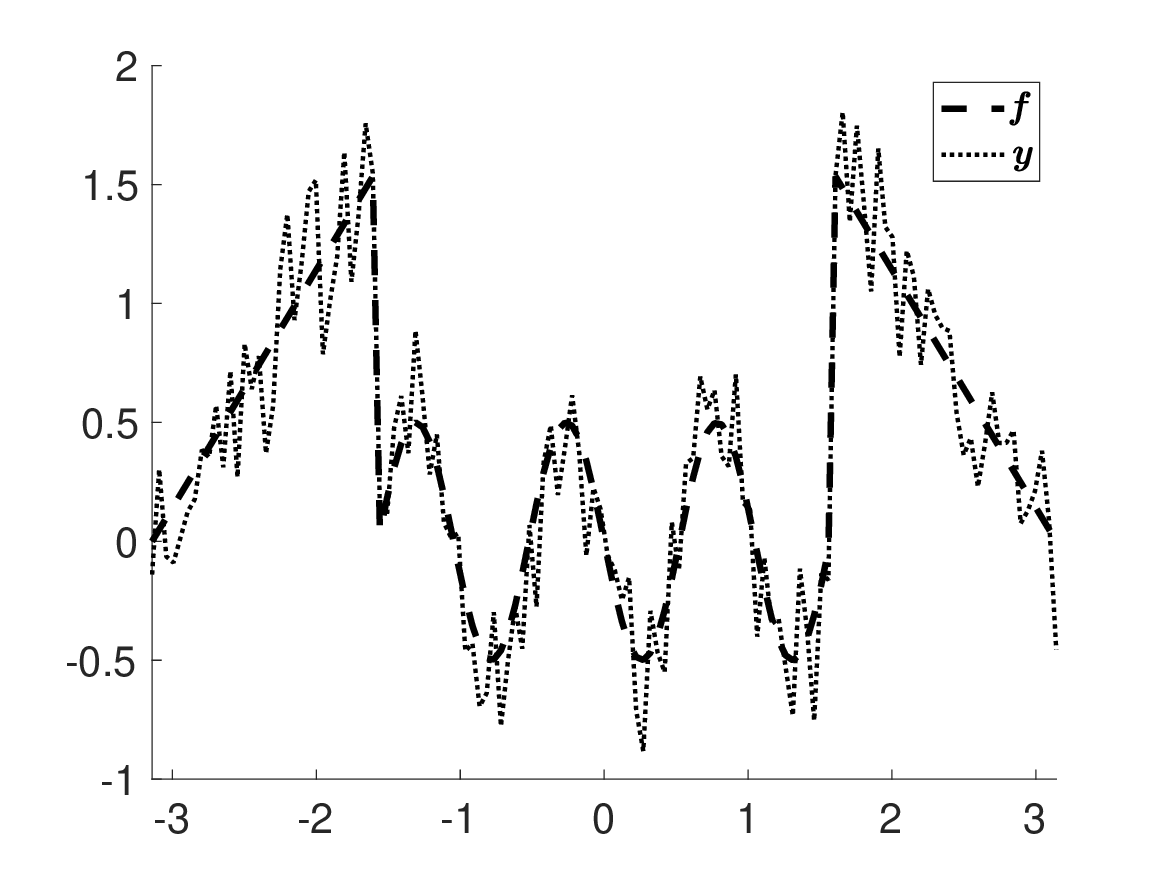}
    \end{subfigure}
    \begin{subfigure}[b]{.32\textwidth}
\includegraphics[width=\textwidth]{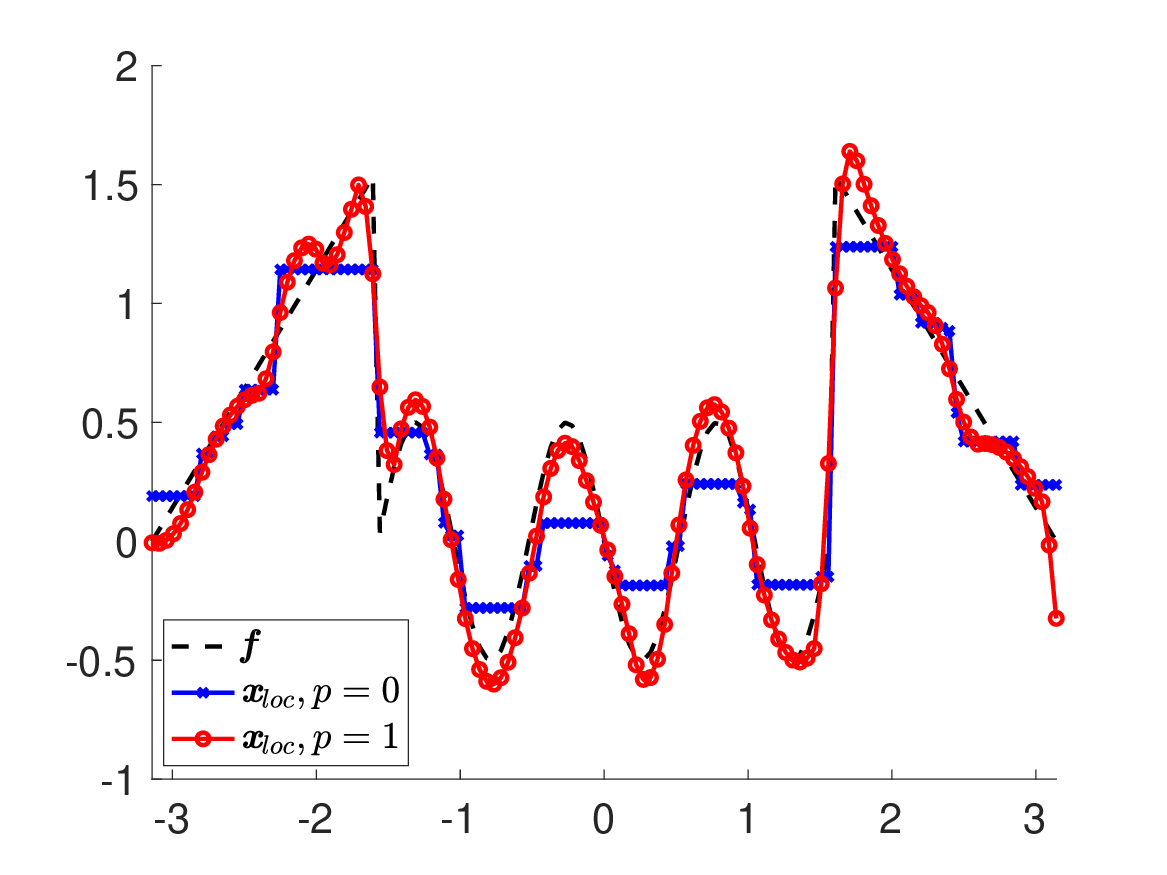}
    \end{subfigure}
    \begin{subfigure}[b]{.32\textwidth}
\includegraphics[width=\textwidth]{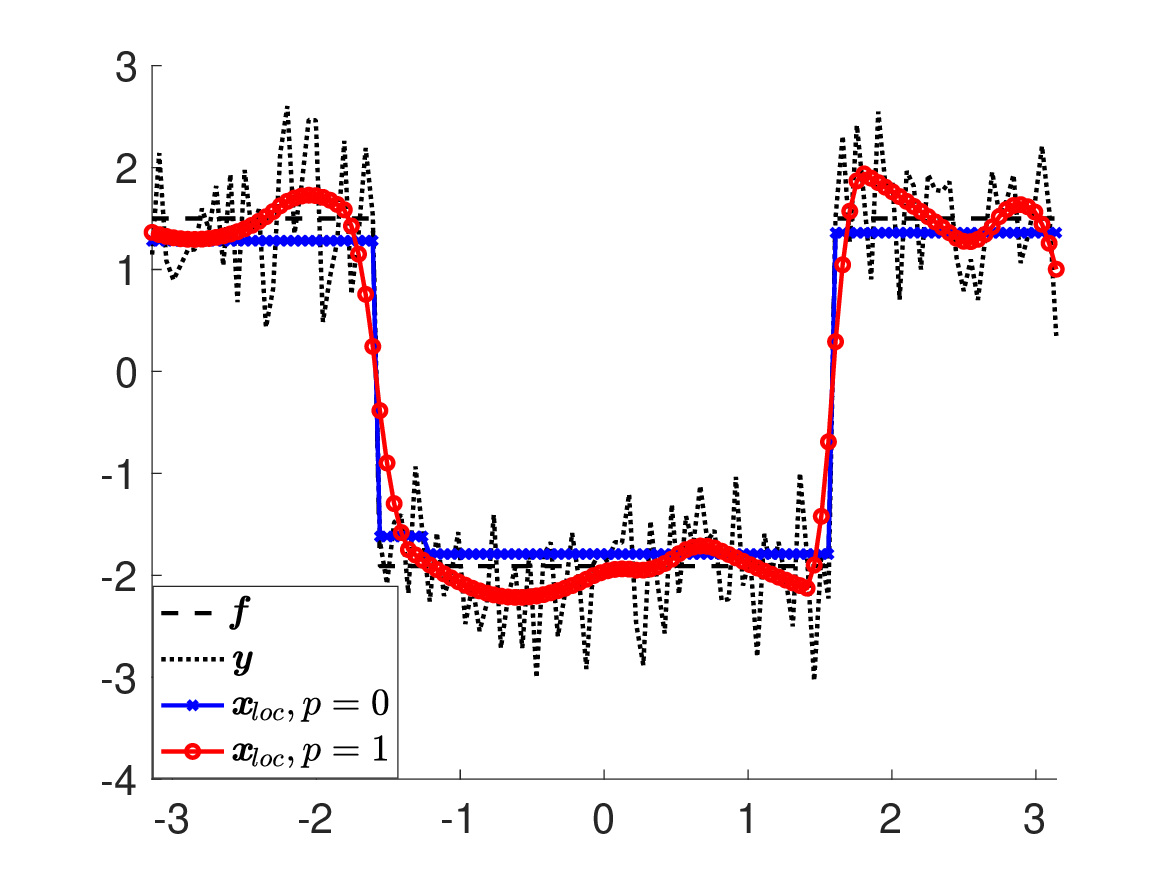}
    \end{subfigure}
    \caption{(left) Discretized underlying signal $\bm f$ (dash) and corresponding measurement $\bm y$ (dot) for the denoising problem ($A = I_n$ in \cref{eq: data_acquisition}). 
    (middle) Comparison of solutions to \cref{eq:lasso_regression} for \cref{ex:example_denoising} with  regularization operators $\mathcal L = T_n^0$ (blue) and $\LL = T_n^1$ (red) with the entries of $T_n^p$ given by \cref{eq:localedgeMV}. 
    (right) Comparison of solutions to \cref{eq:lasso_regression} for \cref{ex:example2}.
    In each case  $n=128$ and $\text{SNR}=10$, giving $\sigma^2=0.038$ (middle) and $\sigma^2=0.295$ (right) in \cref{eq:SNR}, while the regularization parameter is computed from \eqref{eq:lasso_param} with a least squares approximation of ${\bm x}_{est}$.}
    \label{fig:falseOrder}
\end{figure}

Importantly, and a primary motivation for this investigation, is that by construction Lasso regression  assumes that the sparsity promoting operator $\LL$ can be readily determined.  For instance, if the underlying signal ${\bm f}$ is sparse, then choosing $\LL$ to be the identity matrix is appropriate.  Similarly, a first-order differencing matrix is ideal for recovering piecewise constant signals.  However, in many scientific applications ${\bm f}$ represents a discretized piecewise smooth function that may have considerable variability throughout the domain (e.g.~discrete measurements of height or velocity fields that might later be incorporated into simulations of partial differential equations).   A traditional choice for $\LL$ would clearly not promote sparsity in this case.  
\Cref{fig:falseOrder}  illustrates what may happen when the incorrect smoothness assumption is used to choose the regularization operator for the straightforward denoising problem ($A = I_n$ in \cref{eq: data_acquisition}).  Here we compare the solution to \cref{eq:lasso_regression} when $\LL$ is chosen as a first order difference (TV) operator (blue curve, $T_{n}^0$ in \cref{eq:localedgeMV}) and when it is chosen to be a third order local differencing operator (red curve, $T_{n}^1$ in \cref{eq:localedgeMV}).\footnote{As will be described in \cref{sec: local_approach}, the local differencing order is given by $2p+1$, that is $p=0$ yields  first order differencing while $p = 1$ yields third order differencing.} The underlying signal shown on the left (see \cref{ex:example_denoising}) clearly has more variability than can be captured using TV regularization, so using third order differencing greatly improves the results.  Alternatively, when the underlying signal  is piecewise constant, as shown on the right (see \cref{ex:example2}), the high order differencing operator causes artificial oscillations to appear throughout the domain, as it tries to promote piecewise quadratic functions \cite{archibald2005polynomial}. We are therefore motivated to develop an alternative transform operator to handle such variability.  For reasons described in \Cref{sec:introduction}, we refer to our new operator as a {\em residual transform operator}.

\subsection{Edge detection}
\label{sec:edgeDetection}

Before introducing our new residual transform operator, we first lay out some key ingredients that will be useful in its construction. To this end we note that  sparsity promoting operators for piecewise smooth functions can collectively be viewed as {\em edge detectors}.  Moreover, an effective edge detector must be able to discriminate true edges from smooth scale variations in the underlying signal.  This is accomplished by constructing the edge detectors to ``match'' the assumed variability in the smooth regions. An important observation for our new residual transform operator is that such constructions are non-unique.  In particular, two edge detectors designed to match the same type of variability in the smooth regions should produce comparable {\em residuals} in the edge detector approximation.  Hence even when the variability is not known in advance, the {\em difference between two edge detectors designed for the same type of variability will be small.}

It is important to note that there are myriad operators (beyond edge detectors) that can potentially be used to construct a residual transform operator.  Due to their well-known convergence properties that specifically characterize their residuals, this paper employs (i) a local differencing operator and (ii) the  (discrete) Fourier concentration factor edge detection method, which is inherently a global approach. The explicit relationship between these two edge detectors is described in detail in \cite{gelb2006adaptive}.\footnote{A general local differencing formulation (also called the polynomial annihilation method) can be found in  \cite{archibald2005polynomial}, while the concentration factor edge detection method is originally described in \cite{gelb1999detection}.}  Employing these particular methods exploits the contrast of local and global edge detectors, and helps to ensure that the residual transform operator is well-conditioned. Each method is briefly reviewed below for self containment purposes. 

\subsubsection{Jump function approximation}
\label{sec:Edgesgeneral}
Let ${\bm f}\in\mathbb{R}^n$ represent a discretization  of a real-valued $2\pi$-periodic piecewise smooth function $f(s)$ over $n$ uniformly spaced grid points, $s_j = -\pi + (j-1)\Delta s = -\pi + (j-1)\frac{2\pi}{n}$, $j = 1,\dots,n$.  That is,
\begin{equation}
\label{eq:discretization_f}
{\bm f} = \{ f(s_j) \}_{j=1}^n = \{ f_j \}_{j=1}^n. 
\end{equation} 
The corresponding {\em jump function} is defined as 
\begin{equation} \label{eq:jumpfunction}
[f](s):=f(s^+)-f(s^-),
\end{equation}
where $f(s^+)$ and $f(s^-)$ represent the right- and left-hand limits of $f$ at the point $s$.  If $f$ has $\mathcal{K}$ simple discontinuities located at $\{\xi_k\}_{k=1}^{\mathcal{K}} \subset [-\pi,\pi)$, then an equivalent representation  to \cref{eq:jumpfunction} is given by 
\begin{equation}
\label{eq:jumpfunction2}
[f](s) =\sum_{k=1}^{\mathcal K} f(\xi_k)I_{\xi_k}(s), 
\end{equation}
where $I_{\xi_k}(s)$ is the indicator function defined as
\[I_{\xi_k}(s) = \begin{cases}
1, & \text{if }  s=\xi_k, \\
0, & \text{otherwise}.
\end{cases}\]
We seek to approximate $[f](s_j)$ in each cell $[s_j,s_{j+1})$, $j = 1,\dots,n$, with $f(s_{n+1}) = f(s_1)$ due to periodicity.   Each cell entry approximation of $[f](s)$ populates the corresponding {\em edge vector}. The ground truth edge vector will be defined as ${\bm g} \in \mathbb{R}^n$, with a discretized representation ${\bm{g}}_j = [f](s_j)$, $j = 1,\dots,n$ defined on $n$ uniform grid points.

\subsubsection{A local approach to edge detection}
\label{sec: local_approach}
We will employ a {\em local differencing formulation} to detect edges in a cell location between two adjacent grid points $s_j$ and $s_{j+1}$. Consider a pair of adjacent grid points $s_j \leq \xi \leq s_{j+1}$, where $\xi$ is an edge location described in \cref{eq:jumpfunction2}. Under suitable smoothness assumptions on the piecewise regions of $f$ we can use the asymptotic statement

\begin{equation}
\label{eq:asym_stmt}
\Delta f_{j+\frac{1}{2}} := f_{j+1}-f_j=
\begin{cases}
[f](\xi) + \mathcal{O}(\Delta s), & \text{if } j=j_\xi : \xi\in[s_j, s_{j+1}), \\
\mathcal{O}(\Delta s), & \text{otherwise},
\end{cases}
\end{equation}
where $f_j$ is defined in \cref{eq:discretization_f} for $j = 1,\dots,n$ with $\Delta s = \frac{2\pi}{n}$. Note that \cref{eq:asym_stmt} motivates the typical use of TV or equivalently first order difference operators for promoting sparsity in the gradient  domain.  Clearly, \eqref{eq:asym_stmt} also assumes that  $\abs{[f](\xi_k)}\gg\mathcal{O}(\Delta s)$, $k = 1,\dots,\mathcal{K}$, so that variations of $f$ in smooth regions are not falsely identified as edges.

While first order local differencing is appropriate for piecewise constant functions,  it is less effective when $f$ has more variation in smooth regions.  In this case, higher-order differencing can be beneficial since they are explicitly designed to find edges in piecewise polynomial functions.  Indeed, the so-called  {\em polynomial annihilation} edge detectors in \cite{archibald2005polynomial} are constructed using appropriate orders of finite differencing (also formulated for non-uniform grids and multi-dimensional domains).

By expanding \cref{eq:asym_stmt} we obtain the  $2p+1$-order local differencing operator for $j = 1,\dots, n$ as 
\begin{equation}
\label{eq:general_differencing}
\Delta^{2p+1} f_{j+\frac{1}{2}} 
= \sum_{l=0}^p (-1)^{l}\binom{2p+1}{p-l
} (f_{j+1+l} - f_{j-l}),
\end{equation}
which is essentially a scaled version of the Newton divided difference formula.  It follows that in smooth regions, \eqref{eq:general_differencing} provides $\mathcal{O}(\Delta s)^{2p+1}$ accuracy.  Convergence to the jump function can be characterized by the asymptotic statement \cite{gelb2000detection,gelb2002spectral,gelb2006adaptive} 
\begin{equation}
\label{eq:differencing_jump}
\Delta^{2p+1} f_{j+l+\frac{1}{2}}=
\begin{cases}
(-1)^{l} q_{l,p}[f](\xi) + \mathcal{O}(\Delta s), & \abs{l}\leq p, \\
\mathcal{O}(\Delta s)^{2p+1}, & \abs{l} > p,
\end{cases}
\end{equation}
where
\[
q_{l,p} = (-1)^{l}\sum_{k=0}^{p-\abs{l}}\binom{2p+1}{k}(-1)^k = \binom{2p}{p+\abs{l}}.
\]
Note that \cref{eq:differencing_jump} describes the  behavior of \cref{eq:general_differencing} both within and outside of the $2p$-cell neighborhood containing each jump discontinuity.  Finally, we define the $2p+1$-order local edge detection method for $f$ discretized on $n$ equally spaced grid points as \cite{archibald2005polynomial,gelb2006adaptive}
\begin{equation}
\label{eq:local_edge_detector}
T_n^{\Delta_{2p+1}}{f}_j = \frac{1}{q_{0,p}}\Delta^{2p+1} f_{j+\frac{1}{2}},  \quad j = 1,\dots, n,
\end{equation}
where 
\begin{equation}\label{eq:q0P}
 q_{0,p}=\binom{2p}{p}.
\end{equation}
Convergence to edge vector ${\bm g}  = \{[f](s_j)\}_{j = 1}^n$ follows from \cref{eq:differencing_jump}, and we use the $2\pi$-periodicity of $f$ to construct $f_{n+1+p}$ and $f_{n-p}$ as needed in \cref{eq:general_differencing}. Note that $p = 0$ returns standard first order differencing. We can also conveniently express \cref{eq:local_edge_detector} in matrix vector form, $T_{n}^p {\bm f} \approx {\bm g}$,  with the entries  $T_{n}^p(j,l)$, $j,l = 1,\dots,n$, given by \cite{archibald2016image}
\begin{equation}
    \label{eq:localedgeMV}
    {T_{n}^p(j, l)} = \begin{cases}
        \tfrac{(-1)^{l+1}}{q_{0,p}}\binom{2p+1}{p-l}, & j-p\leq l\leq j\\
        \tfrac{(-1)^{l}}{q_{0,p}}\binom{2p+1}{p-l}, & j+1\leq l\leq j+p+1\\
        0, & \text{otherwise,}
    \end{cases}. 
\end{equation} 
providing $(T_n^p{\bm f})_j = T_n^{\Delta_{2p+1}}{f}_j$. Even orders for \cref{eq:local_edge_detector} are similarly generated but are omitted here since they require non-symmetric stencils.  While not instructive in the development of the residual transform operators, even order differencing operators can be just as easily employed in their construction.

\begin{remark}[Higher order differencing]\label{rem:oscillations_local}
While higher order differencing promotes faster convergence outside of the $2p$-cell neighborhood of each edge, spurious oscillations are introduced within these neighborhoods, as described by \cref{eq:local_edge_detector}, which can be further exacerbated by noise in the given measurements. Moreover it may be difficult to resolve multiple edges in sparsely sampled data.  Because of these drawbacks \eqref{eq:lasso_regression} typically employs first order differencing ($p = 0$ in \cref{eq:local_edge_detector})  as the sparsity promoting operator, even when the underlying piecewise smooth signal is known to have variable (rather than constant) behavior between edge locations.  While the new residual transform operator developed here also must consider these issues, as we will see it is more robust to both varying behavior within smooth regions (independent of the chosen $p$) as well as additive noise.
\end{remark}

\subsubsection{A global approach to edge detection}\label{sec:globaledge}
The discrete Fourier transform (DFT) provides a mechanism for obtaining a {\em globally}-based  edge detector from the discretized data \cref{eq:discretization_f}.  Specifically, we begin by calculating  the uniformly sampled discrete Fourier coefficient vector  $\hat{\bm f}$ $\in \mathbb{C}^n$ with components given by 
\begin{equation}
\hat{f}_k = \frac{1}{n} \sum_{j=1}^n f(s_j) e^{-i k s_j}, \quad k = -\frac{n}{2}, \dots, \frac{n}{2}-1,
\label{eq:fourcoefdiscrete}
\end{equation}
which is conveniently written in matrix vector notation as
\begin{equation}
    \label{eq:discretefourMV}
\hat{\bm f} = F{\bm f},
\end{equation}
where $F$ is the discrete Fourier transform matrix with entries given by $F(j,k) = \frac{1}{n}e^{-iks_j}$ and ${\bm f}$ is defined by \cref{eq:discretization_f}. The  {\em concentration factor edge detection method} \cite{gelb1999detection,gelb2006adaptive}  is given by
\begin{equation}
\label{eq:consum}
S^{\mu}_n f(s) 
:=  \pi i \sum_{k=-\frac{n}{2}}^{\frac{n}{2}} \operatorname{sgn}(k) \, \mu\left(\frac{|k|\Delta s}{\pi}\right)sinc\left(\frac{|k|\Delta s}{2}\right)\hat{f}_k \, e^{i k s}.
\end{equation}
Here \(\operatorname{sgn}(k)\) is the sign function and $\mu(\cdot)$ is {\em an admissible} concentration factor function that modulates  $\hat{\bm f}$ so that \cref{eq:consum} converges to $[f](s)$ in \cref{eq:jumpfunction} for any $s \in [-\pi,\pi)$.  In general $\mu(\cdot)$ can be described as an appropriately scaled band pass filter in the Fourier domain with tuning parameters that affect convergence order, localization, and smoothing properties.  A concentration factor $\mu(\cdot)$ is admissible if it satisfies the following requirements: 
\begin{enumerate}
    \item  $\displaystyle K^\mu_n(s) = \sum_{k=1}^{\frac{n}{2}} \mu{\left(\frac{2k}{n}\right)} \sin(ks)$ must be an odd function.

    \item $\frac{\mu(\eta)}{\eta}\in C^2(0,1)$, ensuring the continuity of both the first and second derivatives.

    \item $\int^1_\epsilon \frac{\mu(\eta)}{\eta} \, d\eta \to 1$, $\epsilon = \epsilon(n) > 0$ being small, which guarantees the appropriate normalization in the limit.
\end{enumerate}
Observe that \cref{eq:consum}  is constructed for the {\em discrete} Fourier coefficients in \cref{eq:fourcoefdiscrete}.  When given {\em continuous} Fourier coefficients, i.e. $\hat{f}_k = \frac{1}{2\pi}\int_{-\pi}^{\pi}f(s)e^{-i k s}ds$, integration by parts can be used to derive an admissible concentration factor $\mu(\cdot)$.  Summation by parts is analogously applied in the discrete case, leading to the additional sinc factor \cite{gelb1999detection}.   Various concentration factors $\mu(\cdot)$ satisfying the admissibility requirements have been studied for their localization and convergence properties \cite{gelb1999detection,gelb2002spectral,concentration_paper}. 

We  can  apply \cref{eq:consum} at $s = s_{j+\zeta}$,   $\zeta \in [0,1)$, to approximate  ${\bm g}_j = [f](s_j)$, $j = 1,\dots, n$,  as this is consistent with identifying a jump discontinuity in its enclosed grid cells $[s_j,s_{j+1})$.\footnote{In fact we can prescribe $\zeta_j \in [0,1)$, $j = 1,\dots,n$, but for ease of presentation we choose $\zeta$ to be constant.}   Specifically, direct substitution of \eqref{eq:discretefourMV} into \eqref{eq:consum}  at $s = s_{j+ \zeta}$ yields $S_{n,\zeta}^\mu {\bm f} \approx {\bm g}$ with corresponding matrix coefficients  given by
\begin{equation}\label{eq:matrix_S_mu}
    S_{n,\zeta}^{\mu}(j,l)=\sum_{k=1}^{\tfrac{n}{2}}\frac{\mu\left(\frac{k\Delta s}{\pi}\right)}{k}\left(\cos{(k\Delta s 
    \theta^+_{j,l} 
    )}-\cos{(k\Delta s
    \theta^-_{j,l}}
    )\right), \quad j,l = 1,\dots,n,
\end{equation}
where $\theta^+_{j,l}=j-l+(\tfrac12+\zeta$) and $\theta^-_{i,j}=j-l-(\tfrac12-\zeta$). 

The {\em residual transform operator} developed in \Cref{sec:residualprior} is based on the observation  that  $\mu(\cdot)$ can be chosen so that \cref{eq:consum} evaluated at $s_{j+\frac{1}{2}}$ {\em exactly} matches the differencing behavior of \cref{eq:local_edge_detector}.
In particular, 
\begin{equation}
\label{eq: trig_conc}
    \mu(\eta) := \mu_{2p+1}(\eta) = \frac{2^{2p}\eta \sin^{2p}(\frac{\pi}{2}\eta)}{q_{0,p}}
\end{equation}
 was first identified in \cite{gelb2006adaptive} as the concentration factor yielding the equivalency to the local differencing operator
 \[S^{\mu_{2p+1}}_nf(s_{j+\frac{1}{2}}) = T_n^{\Delta_{2p+1}}f_j, \quad j = 1,\dots,n.\]  
Employing \cref{eq: trig_conc} in 
 \eqref{eq:matrix_S_mu}  provides
\begin{equation}
    \label{eq:matrixconcsum}
S^{\mu_{2p+1}}_{n,{\zeta}}(j,l) = \frac{2^{2p+1}}{n q_{0,p}}\sum_{k = 1}^{\frac{n}{2}} \sin^{2p}{\left(\frac{k\Delta s}{2}\right)}\left(\cos{(k\Delta s\theta^+)}- \cos{(k\Delta s\theta^-)}\right),
 \end{equation}
for $j,l = 1,\dots, n$. It follows that  $(S^{\mu_{2p+1}}_{n,\zeta}{\bm f})_j=S_n^{\mu_{2p+1}}f(s_{j+\zeta})$, and moreover that
\begin{equation}
    T_{n}^p{\bm f} = S^{\mu_{2p+1}}_{n,\frac{1}{2}}{\bm f}
\label{eq:equivalencies}
\end{equation}
for $T_{n}^p$ defined in \cref{eq:localedgeMV} and ${\bm f}$ defined in \cref{eq:discretization_f}. Finally, by construction  $T_{n}^p{\bm f} \approx {\bm g}$ and $S^{\mu_{2p+1}}_{n,\frac{1}{2}}{\bm f}  \approx {\bm g}$, where ${\bm g}_j = [f](s_j)$ are the components of the edge vector corresponding to the piecewise smooth function $f(s)$ on $[-\pi,\pi)$. The accuracy of either 
approximation depends on the underlying variability of $f(s)$, the grid point resolution $\Delta s$, and the suitability of the chosen order $2p+1$.  Since it is instructive in the context of our new residual transform operator, we rederive \cref{eq:equivalencies} below.

\subsubsection{Relationship between local and global edge detectors}
\label{sec:relationlocalglobal}
The main criteria used in developing our new residual transform operator is that we are able to construct two distinct operators of the same order that produce comparable behavior when applied to the numerical solution of \cref{eq: data_acquisition}.  We now demonstrate that the local edge detector in \cref{eq:local_edge_detector} and the global edge detector in \cref{eq:consum} are two such operators.

We begin by rewriting \cref{eq:consum} as
\[S^{\mu_{2p+1}}_n f(s)=\sum_{k=-\frac{n}{2}}^{\frac{n}{2}}
G^{\mu_{2p+1}}_n(k)\hat{f}_k e^{i k s},\]
where $G^{\mu_{2p+1}}$ combines the multiplicative factors in \cref{eq:consum} with the specific concentration factor $\mu_{2p+1}$  in \cref{eq: trig_conc} so that
\begin{equation}
\label{eq:G_mu}
G^{\mu_{2p+1}}_n(k)=i\operatorname{sgn(k)}\tfrac{2^{2p+1}}{q_{0,p}}\sin^{2p+1}\left(\tfrac{\abs{k}\Delta s}{2}\right),
\end{equation}
with $q_{0,p}$  defined in \cref{eq:q0P}.  Using the equivalent representation
\[\operatorname{sin}^{2p+1}(\xi) = \tfrac{1}{2^{2p}}\sum_{m=0}^p(-1)^m\binom{2p+1}{p-m}\operatorname{sin}\left((2m+1)\xi\right),\]
we can rewrite \eqref{eq:G_mu} as
\begin{eqnarray*}
    G^{\mu_{2p+1}}_n(k) &=& i\operatorname{sgn(k)}\tfrac{2^{2p+1}}{q_{0,p}}\left(\tfrac{1}{2^{2p}}\sum_{l=0}^p(-1)^l\binom{2p+1}{p-l}\operatorname{sin}\left(\abs{k}\Delta s\left(l+\tfrac{1}{2}\right)\right)\right)\\
    &=& \operatorname{sgn(k)}\tfrac{2 i}{q_{0,p}}\sum_{l=0}^p(-1)^l\binom{2p+1}{p-l}\operatorname{sin}\left(\abs{k}\Delta s\left(l+\tfrac{1}{2}\right)\right).
\end{eqnarray*} 

Finally, since the jump value must occur between two adjacent grid points, by using the midpoint $s_{j+\frac{1}{2}}$ as the point of approximation in \cref{eq:consum} we obtain the approximation to each component  of sparse domain vector ${\bm g} = \{[f](s_j)\}_{j= 1}^n$ as
\begin{align*}
    S^{\mu_{2p+1}}_n {f}(s_{j+\frac{1}{2}})&=\sum_{k=-\frac{n}{2}}^{\frac{n}{2}}G^{\mu_{2p+1}}_n(k)\hat{f}_k e^{iks_{j + \frac{1}{2}}} \\
    &= \sum_{k=-\frac{n}{2}}^{\frac{n}{2}}
e^{i\frac{k\Delta s}{2}}G^{\mu_{2p+1}}_n(k)\hat{f}_k e^{iks_{j}} =: (S_{n,\frac{1}{2}}^{\mu_{2p+1}}{\bm f)_j}.
\end{align*}

From the local direction, we observe that  since $s_j = -\pi + (j-1)\Delta s$ are equally spaced, by directly applying the DFT \cref{eq:discretefourMV} and noting its interpolatory properties (see e.g.~\cite{HGG}), we can express \cref{eq:local_edge_detector} as
\[
T^{\Delta_{2p+1}}_n {f}_j=\sum_{k=-\tfrac{n}{2}}^{\tfrac{n}{2}}
G^{\Delta_{2p+1}}_n(k)\hat{f}_ke^{iks_j}.
\]
It then follows from the differencing formula in \cref{eq:general_differencing} that

\begin{align}
 G^{\Delta_{2p+1}}_{n}(k) 
&= i \frac{2 k \Delta s}{q_{0,p}} e^{i\frac{k\Delta s}{2}}\sum_{l = 0}^p (-1)^l \left(l + \frac{1}{2}\right) \binom{2p+1}{p - l} sinc{\left(k\Delta s\left(l+\frac{1}{2}\right)\right)}\nonumber\\
   &= e^{i\frac{k\Delta s}{2}}\frac{2 i }{q_{0,p}} \sum_{l = 0}^p (-1)^l \binom{2p+1}{p - l} \operatorname{sin}{\left(k\Delta s\left(l+\frac{1}{2}\right)\right)} =e^{i\frac{k\Delta s}{2}}G^{\mu_{2p+1}}_n(k),\label{eq:G_Delta}
\end{align}
establishing \cref{eq:equivalencies}. Observe that this does not mean that $T_{n}^p = S^{\mu_{2p+1}}_{n,\frac{1}{2}}$, but rather {\em only} that $T_{n}^p{\bm f} = S^{\mu_{2p+1}}_{n,\frac{1}{2}}\bm f$ for ${\bm f}$ in \cref{eq:discretization_f}.  Nor does it imply that either construction is a good approximation to the true sparse domain vector ${\bm g}$.  Regardless, unlike the case for a standard sparsity promoting transform operator (such as TV), the residual transform operator that we develop can  still be effective, even if the assumptions regarding the true sparse domain vector are  incorrect.

We can also show this equivalency from the perspective of Fourier collocation methods \cite{HGG}.   In this regard,
starting from \eqref{eq:consum} at grid point $s_{j+\frac{1}{2}}$, we have
\begin{align*}
    &S_n^{\mu_{2p+1}}f\left(s_{j+\tfrac{1}{2}}\right) \\ 
    &=  \pi i \sum_{k=-\frac{n}{2}}^{\frac{n}{2}} \operatorname{sgn}(k) \, \mu\left(\frac{|k|\Delta s}{\pi}\right)sinc\left(\frac{|k|\Delta s}{2}\right)\hat{f}_k \, e^{i k s_j}e^{i\tfrac{k\Delta s}{2}}\\
    &=i\tfrac{2^{2p+1}}{q_{0,p}}\sum_{k=-\frac{n}{2}}^{\frac{n}{2}}\operatorname{sin}^{2p+1}\left(\frac{|k|\Delta s}{2}\right)\hat{f}_k \, e^{i k s_j}e^{i\tfrac{k\Delta s}{2}}\\
    &=i\tfrac{2^{2p+1}}{q_{0,p}}\sum_{k=-\frac{n}{2}}^{\frac{n}{2}}\left(\tfrac{1}{2^{2p}}\sum_{l=0}^p(-1)^l \binom{2p+1}{p-l}\operatorname{sin}\left(\frac{|k|\Delta s_j}{2}\right)\right)\hat{f}_k \, e^{i k s_j}e^{i\tfrac{k\Delta s}{2}}.
\end{align*}
Use of the identity
\[e^{ik\Delta s}-1=2\,i\,e^{i\frac{k\Delta s}{2}}\sin\left(\frac{k\Delta s}{2}\right)\] 
further simplifies this expression to
\begin{align}\label{eq:Sf_Tf}
\begin{split}
    S_n^{\mu_{2p+1}}f\left(s_{j+\tfrac{1}{2}}\right) 
    &= \tfrac{2\; i}{q_{0,p}}\sum_{l=0}^p(-1)^l \binom{2p+1}{p-l}\sum_{k=-\frac{n}{2}}^{\frac{n}{2}}\tfrac{e^{ik\Delta s(2l+1)-1}}{2 i e^{i(2l+1)\tfrac{k\Delta s}{2}}} \,\hat{f}_k \, e^{i k s_j}e^{i\tfrac{k\Delta s}{2}}\\
    &=\tfrac{1}{q_{0,p}}\sum_{l=0}^p(-1)^l \binom{2p+1}{p-l}\sum_{k=-\frac{n}{2}}^{\frac{n}{2}}\hat{f}_k\left( e^{ik\Delta s(j+l+1)}-e^{ik\Delta s(j-l)} \right)\\
    &= \tfrac{1}{q_{0,p}}\sum_{l=0}^p(-1)^l \binom{2p+1}{p-l}\left( f_{j+l+1}-f_{ j-l} \right) =: T_n^{\Delta_{2p+1}}{f}_j, 
\end{split}
\end{align}
which follows from \eqref{eq:local_edge_detector}. Hence once again we see that  for {\em the exact realization} ${\bf f} = \{f(s_j)\}_{j = 1}^n$ in \cref{eq:discretization_f}, $T_{n}^p{\bm f} = S^{\mu_{2p+1}}_{n,\frac{1}{2}}{\bm f}$, even though $T_{n}^p \ne S^{\mu_{2p+1}}_{n,\frac{1}{2}}$. 

\subsubsection*{Near equivalent formulation}
While the equality of \cref{eq:equivalencies} is what motivates our new residual transform operator, for conditioning purposes, which will be further explained in \Cref{sec:residualprior}, we instead consider $S_{n,\zeta}^{\mu_{2p+1}} {\bm f}$, $\zeta \ne \frac{1}{2}$.  In this case,  although $T_n^p{\bm f} \ne S_{n,\zeta}^{\mu_{2p+1}} {\bm f}$, it is still true that  
\begin{equation}\label{eq:nearequivalent}
T_n^p{\bm f} \approx  S_{n,\zeta}^{\mu_{2p+1}} {\bm f}.\end{equation}
We show this by first observing that the derivation of \cref{eq:Sf_Tf} can also be used for general $\zeta \in [0,1)$  to obtain
\begin{align*}
(S_{n,\zeta}^{\mu_{2p+1}}\bm f)_j &= S_n^{\mu_{2p+1}}f\left(s_{j+\zeta}\right) \\
    &= \tfrac{2\; i}{q_{0,p}}\sum_{l=0}^p(-1)^l \binom{2p+1}{p-l}\sum_{k=-\frac{n}{2}}^{\frac{n}{2}}\tfrac{e^{ik\Delta s(2l+1)-1}}{2 i e^{i(2l+1)\tfrac{k\Delta s}{2}}} \,\hat{f}_k \, e^{i k s_j}e^{ik\zeta\Delta s}\\
    &=\tfrac{1}{q_{0,p}}\sum_{l=0}^p(-1)^l \binom{2p+1}{p-l}\sum_{k=-\frac{n}{2}}^{\frac{n}{2}}\hat{f}_k\left( e^{ik\Delta s(j+l+1)}-e^{ik\Delta s(j-l)} \right)e^{ik\Delta s(\zeta-\frac12)}.
\end{align*}
It then follows from \cref{eq:G_Delta} that 
\begin{align*}
    &\abs{(S_{n,\zeta}^{\mu_{2p+1}}\bm f)_j-(T_n^{p}\bm f)_j}\\
    \leq& \sum_{l=0}^p\tfrac{(-1)^l}{q_{0,p}} \binom{2p+1}{p-l}\sum_{k=-\frac{n}{2}}^{\frac{n}{2}}\hat{f}_k\left( e^{ik\Delta s(j+l+1)}-e^{ik\Delta s(j-l)} \right)\abs{e^{ik\Delta s(\zeta-\frac12)}-1}\\
    =& \sum_{l=0}^p\tfrac{(-1)^l}{q_{0,p}} \binom{2p+1}{p-l}\sum_{k=-\frac{n}{2}}^{\frac{n}{2}}\hat{f}_k\left( e^{ik\Delta s(j+l+1)}-e^{ik\Delta s(j-l)} \right)\abs{2\sin\left(\frac{k\Delta s(\zeta-\frac12)}{2}\right)}\\
    \leq& \frac{\abs{\Delta s}}{2}\sum_{l=0}^p\tfrac{(-1)^l}{q_{0,p}} \binom{2p+1}{p-l}\sum_{k=-\frac{n}{2}}^{\frac{n}{2}}\abs{k}\hat{f}_k\left( e^{ik\Delta s(j+l+1)}-e^{ik\Delta s(j-l)} \right)=\mathcal{O}(\Delta s),
\end{align*}
resulting in \cref{eq:nearequivalent}, which will be directly exploited in the development of the new residual transform operator introduced in \Cref{sec:residualprior}.
\begin{remark} [Shifted concentration factor]\label{rem:shifted}
    Setting $\zeta\ne \tfrac12$ leads to a modified partial sum approximation in \cref{eq:consum} with $\mu(\eta) = \mu^{shift}(\eta)$.  It then follows from \cref{eq: trig_conc} that when  $\mu(\eta) = \mu_{2p+1}(\eta)$, we obtain  $\mu_{2p+1}^{shift}{(\eta)} = e^{ik\zeta\Delta s}\mu_{2p+1}(\eta) = e^{i\pi\zeta \eta}\mu_{2p+1}(\eta)$, where $\eta=\frac{k\Delta s}{\pi}$ as in \eqref{eq:matrix_S_mu}. 
    While the first two admissibility conditions are satisfied, $\mu_{2p+1}^{shift}(\eta)$ is {\em not}  admissible since by construction the peaks of \cref{eq:consum} occur at the cell average ($\zeta = \tfrac12$). That is, the shift causes the peak value to be replaced with something smaller in magnitude.    Renormalization is possible but there is no closed form for the corresponding constant.  Moreover, related to the discussion in \Cref{rem:deficientmatrix}, the resulting residual transform operator introduced in \Cref{sec:residualprior} for an admissible shifted concentration factor may again be singular or ill-conditioned, which we are trying to avoid. 
\end{remark}

\section{A new residual transform operator}
\label{sec:residualprior}

We now have all of the ingredients needed for our new residual transform operator. We begin by defining the estimation error $\boldsymbol\epsilon$ as \[\boldsymbol\epsilon = \bm x-\bm f,\] with $\bm x$ being the solution obtained in \cref{eq:lasso_regression} and $\bm f$ being the true signal discretized by \eqref{eq:discretization_f}.  From \cref{eq:nearequivalent} it follows that
\[|T_{n}^p{\bm x} - S_{n,\zeta}^{\mu_{2p+1}}{\bm x}| \approx |T_{n}^p{\boldsymbol{\epsilon}} - S_{n,\zeta}^{\mu_{2p+1}}{\boldsymbol{\epsilon}}| = |(T_{n}^p - S_{n,\zeta}^{\mu_{2p+1}}){\boldsymbol{\epsilon}}|.\]
Therefore, by defining a new {\em residual transform  operator}  as
\begin{equation}
    \label{eq:residualprior}
    R_{n,\zeta}^p := T_{n}^p-S_{n,\zeta}^{\mu_{2p+1}},
\end{equation}
we have by design that
\[|R_{n,\zeta}^p{\bm x}| \approx |(T_{n}^p - S_{n,\zeta}^{\mu_{2p+1}}){\boldsymbol{\epsilon}}|.\]
Moreover, an immediate consequence of the definition \eqref{eq:residualprior} is that 
\begin{equation}\label{eq:residualerror} |R_{n,\zeta}^p{\bm x}| \le |T_{n}^p{\bm x}| \mbox{ and } |R_{n,\zeta}^p{\bm x}| \le |S_{n,\zeta}^{\mu_{2p+1}}{\bm x}|,
\end{equation}
which has direct implications on their respective utility in \cref{eq:lasso_regression} as  sparsifying transform operator ${\mathcal L}$.
As we proceed to describe, for a general $2\pi$-periodic piecewise smooth function $f(s)$ on $[-\pi,\pi)$, the residual transform operator yields a more accurate description of the underlying discretized signal ${\bm f}$ in \cref{eq:discretization_f} than standard sparsifying transform operators do. Specifically, both the more traditional local sparsity-promoting  operators (represented here by $T_{n}^p$) as well as the non-standard global sparsity-promoting operator ($S_{n,\zeta}^{\mu_{2p+1}}$) rely on the assumption that the variability is of {\em fixed  type} in {\em all} smooth regions, for example the underlying function is {\em only} piecewise constant or {\em only} piecewise linear throughout the domain.  Therefore {\em by construction} they fail to capture the inherent behavioral irregularity  in real-world signals,  which may lead to a complete mischaracterization of the underlying signal. Moreover, since the true behavior in smooth regions is often unknown in practice, such transform operators  must rely on an assumed order $2p+1$ (often $p = 0$ corresponding to a sparse gradient domain), further limiting their effectiveness even in fixed variability cases.  By contrast, the residual transform operator is designed to better account for potential behavioral variability in smooth regions.

\begin{remark}[Residual Transform Operator]\label{rem:residualterminology}
While the discussion above provides a compelling reason to use \cref{eq:residualprior}, the name {\em residual transform operator} is motivated as follows.  Suppose ${\bm g}$ is the edge vector corresponding to the signal vector ${\bm f}$. The corresponding residual vectors for each approximation $T_{n}^p{\bm f}$ and $S_{n,\frac{1}{2}}^{\mu_{2p+1}}{\bm f}$ are then defined as
\begin{equation*}
    {\bm r}_1 = {\bm g} - T_{n}^p{\bm f}, \quad {\bm r}_2 = {\bm g} - S_{n,\frac{1}{2}}^{\mu_{2p+1}}{\bm f}.
\end{equation*}
Therefore,
\[R_{n,\frac{1}{2}}^p{\bm f} = {\bm r}_1 - {\bm r}_2 = T_{n}^p{\bm f} -S_{n,\frac{1}{2}}^{\mu_{2p+1}}{\bm f}= {\bm 0},\]
thus inspiring \cref{eq:residualprior} and its name. Importantly, we observe that the approximation to the sparse domain vector ${\bm g}$ by either $T_{n}^p{\bm f}$ or $S_{n,\frac{1}{2}}^{\mu_{2p+1}}{\bm f}$ does not have to be accurate, as this never enters into the calculation of $R_{n,\frac{1}{2}}^p{\bm f}$.  This is in sharp contrast to more traditional sparsifying transform operators whose heavy reliance on the accuracy of this approximation is evident from \cref{eq:residualerror} and observed in \Cref{fig:falseOrder}.
\end{remark}

\begin{remark}[singularity  of $R_{n,\frac{1}{2}}^p$]
\label{rem:deficientmatrix}
The residual transform operator $R_{n,\zeta}^p$ as defined in \cref{eq:residualprior} is intuitively appealing due to the equality of \cref{eq:equivalencies} for $\zeta = \frac{1}{2}$. Unfortunately  $R_{n,\frac{1}{2}}^p$ is a  severely rank-deficient singular matrix with a large null space. In particular, $R_{n,\frac{1}{2}}^p$ is rank 2 for $p = 0,1$.  Therefore the regularizer only attenuates noise  within a specific plane of the n-dimensional discretized signal, rendering it ineffective for most practical purposes. While a small perturbation   around $\zeta = \tfrac12$ in $S_{n,\zeta}^{\mu_{2p+1}}$ would technically increase the rank, most of resulting singular values would remain  close to zero, severely limiting its utility.
\end{remark}

To address the rank deficiency of $R^p_{n,\tfrac12}$  and create a robust operator $R^p_{n,\zeta}$, we instead choose $\zeta=\tfrac14$. While this choice means the convenient equality between $T_{n}^p\bm f$ and $S_{n,\zeta}^{\mu_{2p+1}}\bm f$ no longer holds, it does provide the crucial benefit of ensuring that the operator  $R_{n,\zeta}^p$ is well-conditioned. Moreover,  having the  {\em approximate} equality $T_{n}^p\bm f \approx S_{n,\zeta}^{\mu_{2p+1}}\bm f$ in \cref{eq:nearequivalent} still means that the smooth regions in the underlying solution do not have to be characterized by a {\em known} particular type of variability. As  will be evident in our numerical experiments, the shift in the approximation cell value to $\zeta = \tfrac14$ has the most impact near the discontinuities, where some smoothing occurs.

\section{Numerical results} 
\label{sec:numerics} 

Our numerical experiments are designed to consider two representative problems. The first describes a common situation where the underlying piecewise smooth function exhibits differences in the types of variability in the smooth regions, and can be illustrated by
\begin{example}
    \label{ex:example_denoising}
$f_1(s)$ is a $2\pi$-periodic piecewise smooth on $[-\pi,\pi)$ such that
    \begin{equation*}
        f_1(s) =
\begin{cases} 
s+\pi, &  s < -\frac{\pi}{2}, \\
-\frac{1}{2}\sin\left(6s\right), &  -\frac{\pi}{2}  \leq s < \frac{\pi}{2}, \\
-s+\pi, &  s \geq \frac{\pi}{2}.
\end{cases}
    \end{equation*}
\end{example}   
We demonstrate the advantage of our proposed residual transform operator $R^p_{n,\zeta}$, with $\zeta = \frac{1}{4}$ in \cref{eq:lasso_regression} over both first order differencing ($p = 0$ in \cref{eq:local_edge_detector}, equivalent to using  TV) and third order differencing ($p = 1$ in \cref{eq:local_edge_detector}).We compare performances for three  basic tasks in the recovery of $f_1(s)$: (1) denoising, (2) deblurring, and (3) undersampling.   When the input data are noisy and/or undersampled we compute the regularization parameter $\alpha$ using \eqref{eq:lasso_param}, with $\bm x_{est}$  calculated as the least squares estimate employing the corresponding forward model in \cref{eq: data_acquisition}. As there is no additive noise in the deblurring problem, in this case we instead heuristically choose  $\alpha$.

The second scenario applies this same suite of tests to a piecewise constant function, thereby establishing a baseline performance for our residual transform operator on simpler structures.   The respective illustrative example is given by
\begin{example}
    \label{ex:example2}
    $f_2(s)$ is a $2\pi$-periodic piecewise smooth on $[-\pi,\pi)$ such that  
    \begin{align*}
        f_2(s) =
\begin{cases} 
1.5, &  s < -\frac{\pi}{2}, \\
-\frac6\pi, & -\frac{\pi}{2}  \leq s < \frac{\pi}{2}, \\
1.5, &  s \geq \frac{\pi}{2},
\end{cases}
    \end{align*}
\end{example}
Since TV is designed specifically for recovering piecewise constant signals, we do not expect our new method to outperform TV in approximating \cref{ex:example2}.

\subsection*{Analyzing our numerical results}
We use the signal to noise ratio $\text{SNR}$ defined as 
\begin{equation}
    \label{eq:SNR}
    \text{SNR}_{\text{dB}} = 10\log_{10}\left(\frac{P_{signal}}{P_{noise}}\right) = 10\log_{10}\left(\frac{\norm{\bm f}^2_2}{n\sigma^2}\right),
\end{equation}
to determine the values of $\sigma^2$ in our experiments.  In each suite of tests we consider various SNR levels  to evaluate our method's robustness to noise.

We also utilize metrics for the relative error $E^{rel}$
\begin{equation}
    \label{eq:relative_err}
    E^{rel}= \sum_{j = j_{min}}^{ j_{max}}\frac{\abs{\bm x_j-\bm f_j}}{\abs{\bm f_j}}, 
\end{equation}
where $1 \le j_{min} < j_{max} \le n$ reflects the interval $[s_{min},s_{max}]$ in which we compute the relative error,
along with the pointwise absolute error $E^{abs}$ with components
\begin{equation}
    \label{eq:absolute_err}
    E^{abs}_j=\abs{{\bm x}_j-\bm f_j},\quad j = 1,\dots, n.
\end{equation}
In each case $\bm x$ is calculated using \eqref{eq:lasso_regression} for the given regularization operator $\mathcal L$ and ${\bm f}$ is the discretization of the underlying piecewise smooth function $f(s)$ given in \cref{eq:discretization_f}.

In all figures, ${\bm x}_{loc}$ refers to the solution obtained using \cref{eq:lasso_regression} with $\mathcal L = T_n^p$, for either $p = 0$ (first order) or $p = 1$ (third order) local differencing, as will be specified in each figure caption. We  will similarly use ${\bm x}_R$ to represent the solution using \cref{eq:lasso_regression} with $\mathcal L = R_{n,\zeta}^p$ as defined by \cref{eq:residualprior}. For ease of presentation, our numerical results are only displayed for $n = 128$ and $\zeta = \tfrac14$ for $R^p_{n,\zeta}$ in \cref{eq:residualprior}.  Other choices for resolutions and shifts yield consistent outcomes for all of our tests and are not shown here.

\subsection{Denoising}
\label{sec:denoising}
We first consider the denoising problem, i.e. $A=I_n$ in \eqref{eq: data_acquisition}. In this case each  sampled measurement is corrupted by additive independent and identically distributed (i.i.d.) zero-mean Gaussian noise with decreasing SNR levels.

\begin{figure}[h!]
    \centering
\begin{subfigure}[b]{.24\textwidth}
\includegraphics[width=\textwidth]{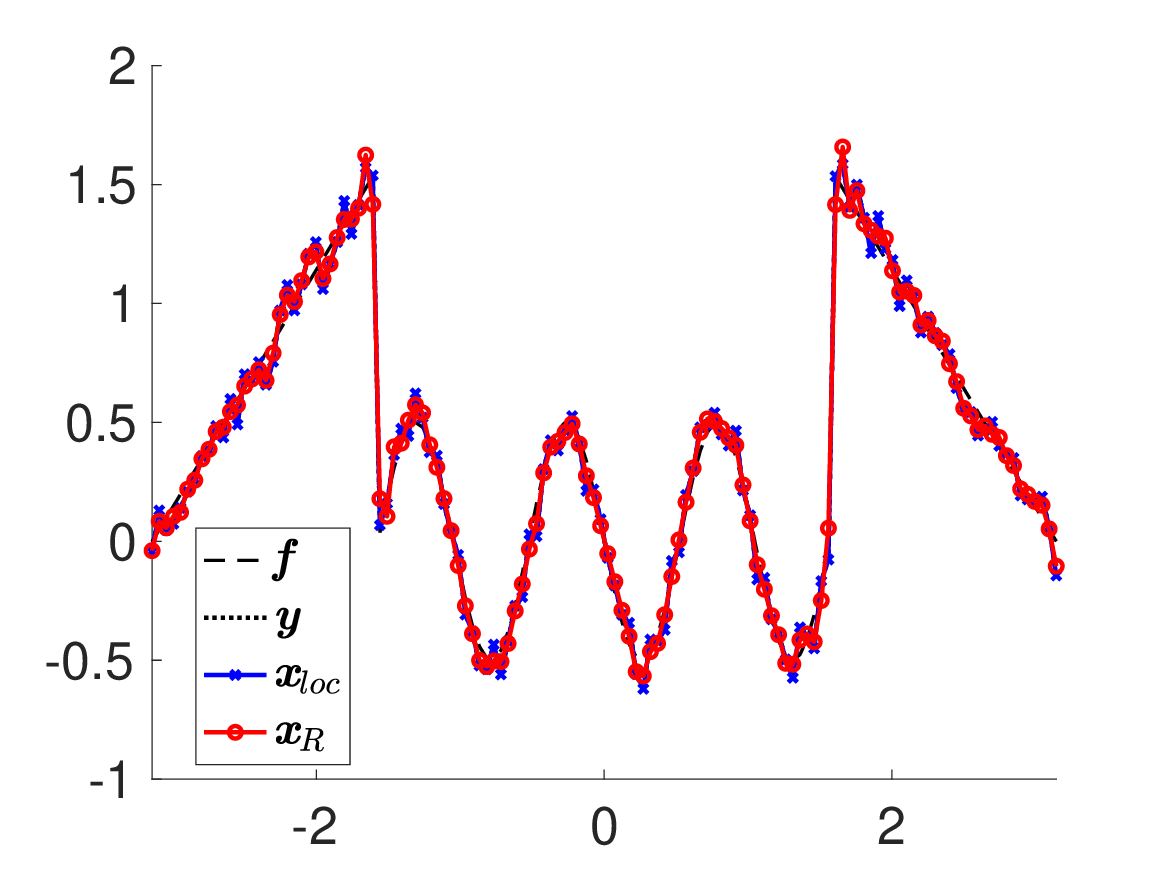}
\end{subfigure}
\begin{subfigure}[b]{.24\textwidth}
\includegraphics[width=\textwidth]{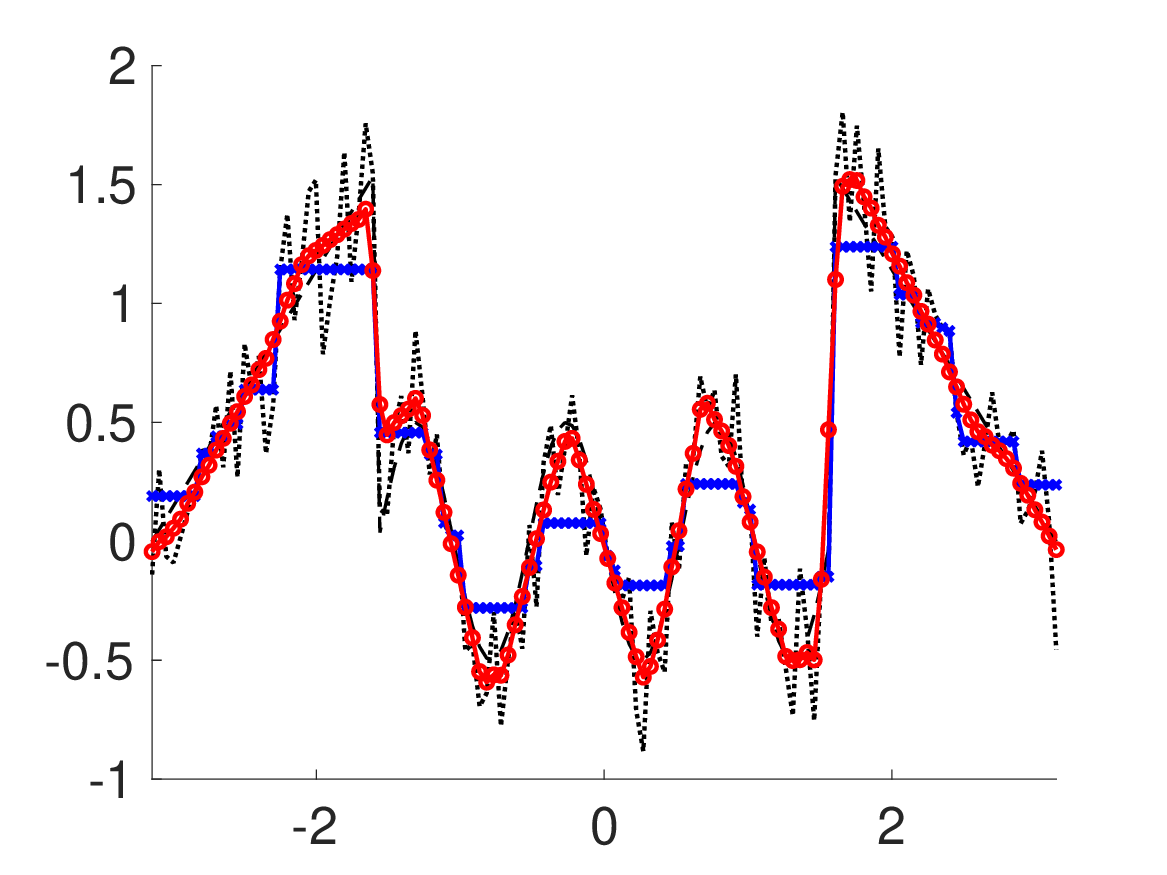}
\end{subfigure}
\begin{subfigure}[b]{.24\textwidth}
\includegraphics[width=\textwidth]{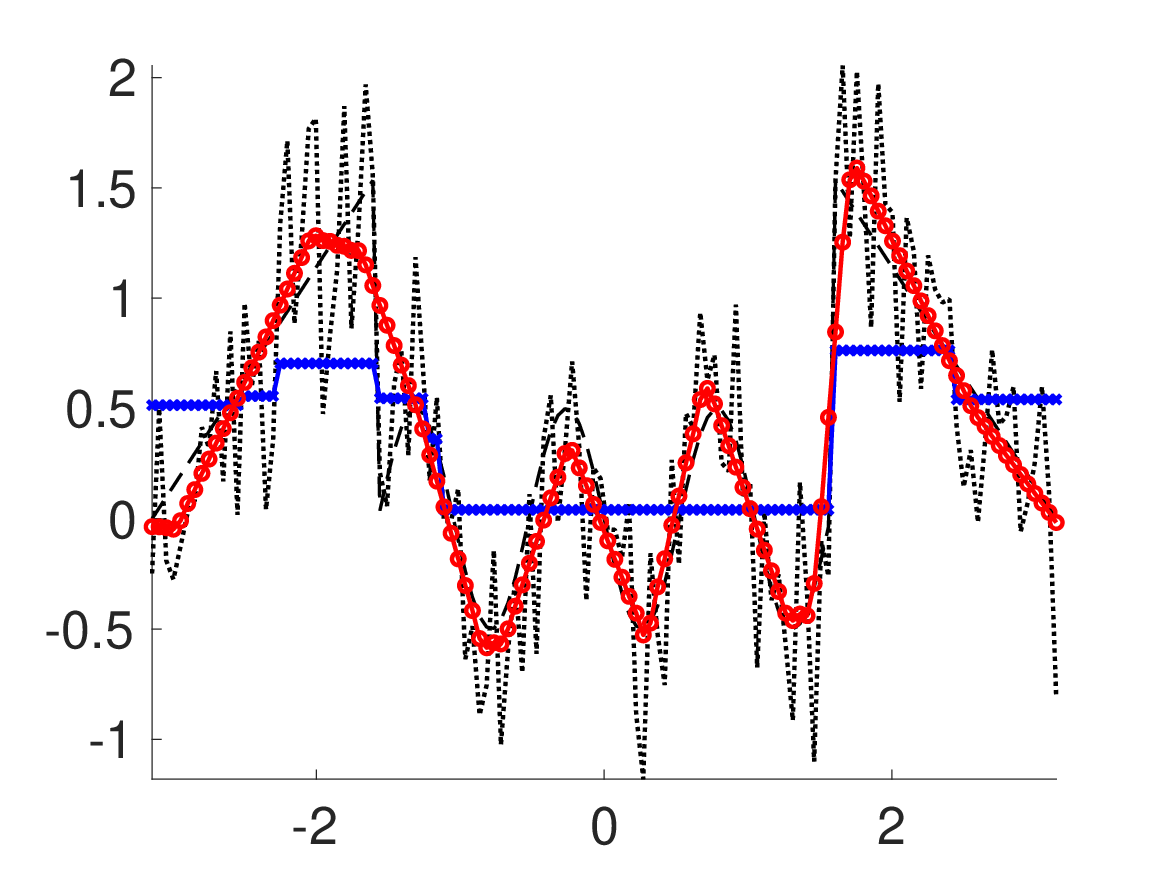}
\end{subfigure}
\begin{subfigure}[b]{.24\textwidth}
\includegraphics[width=\textwidth]{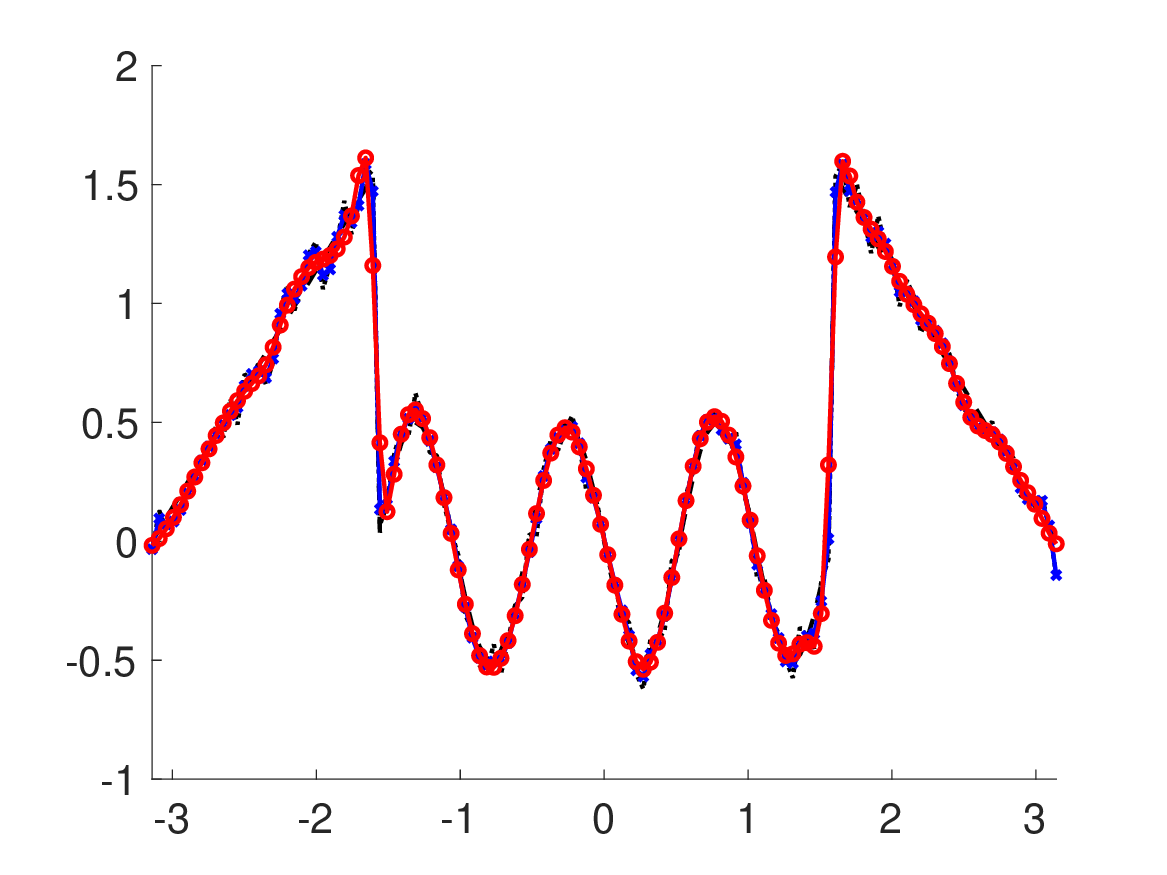}
\end{subfigure}
\\
\begin{subfigure}[b]{.24\textwidth}
\includegraphics[width=\textwidth]{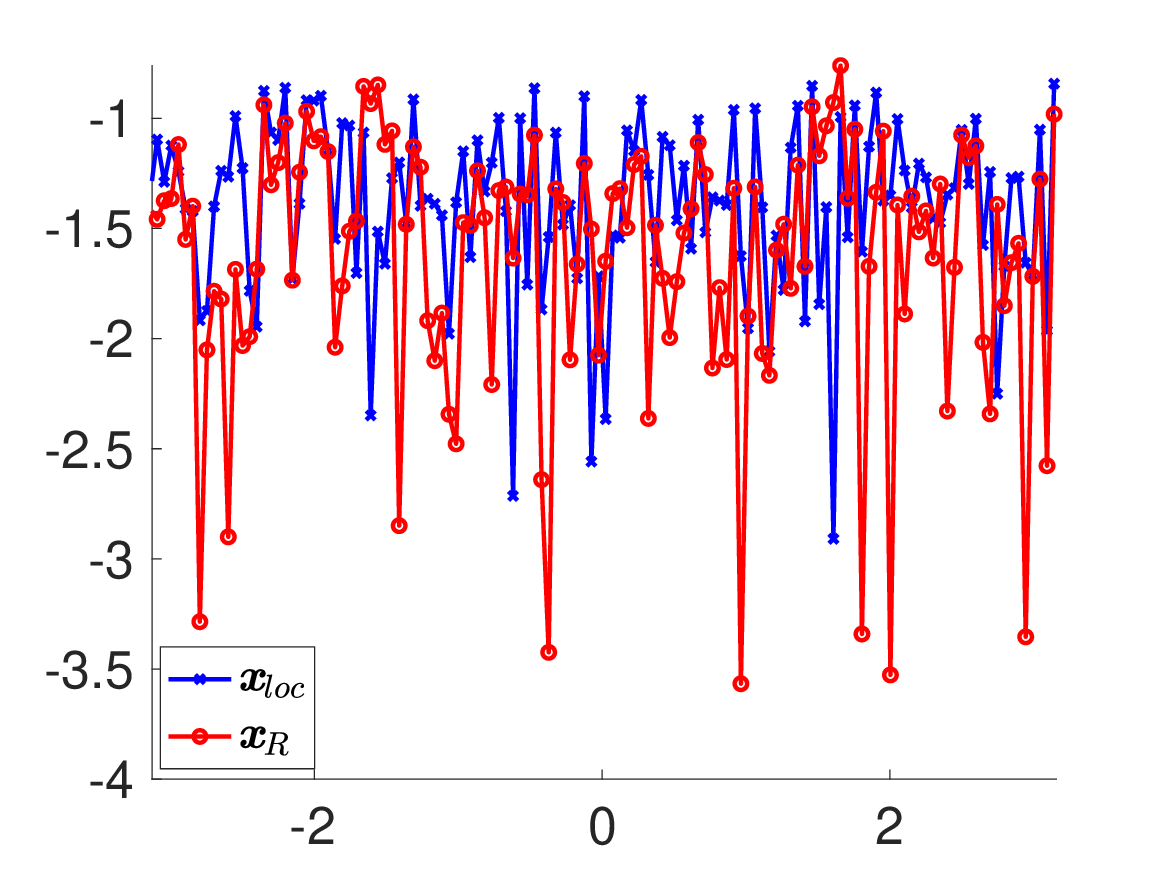}
\end{subfigure}
\begin{subfigure}[b]{.24\textwidth}
\includegraphics[width=\textwidth]{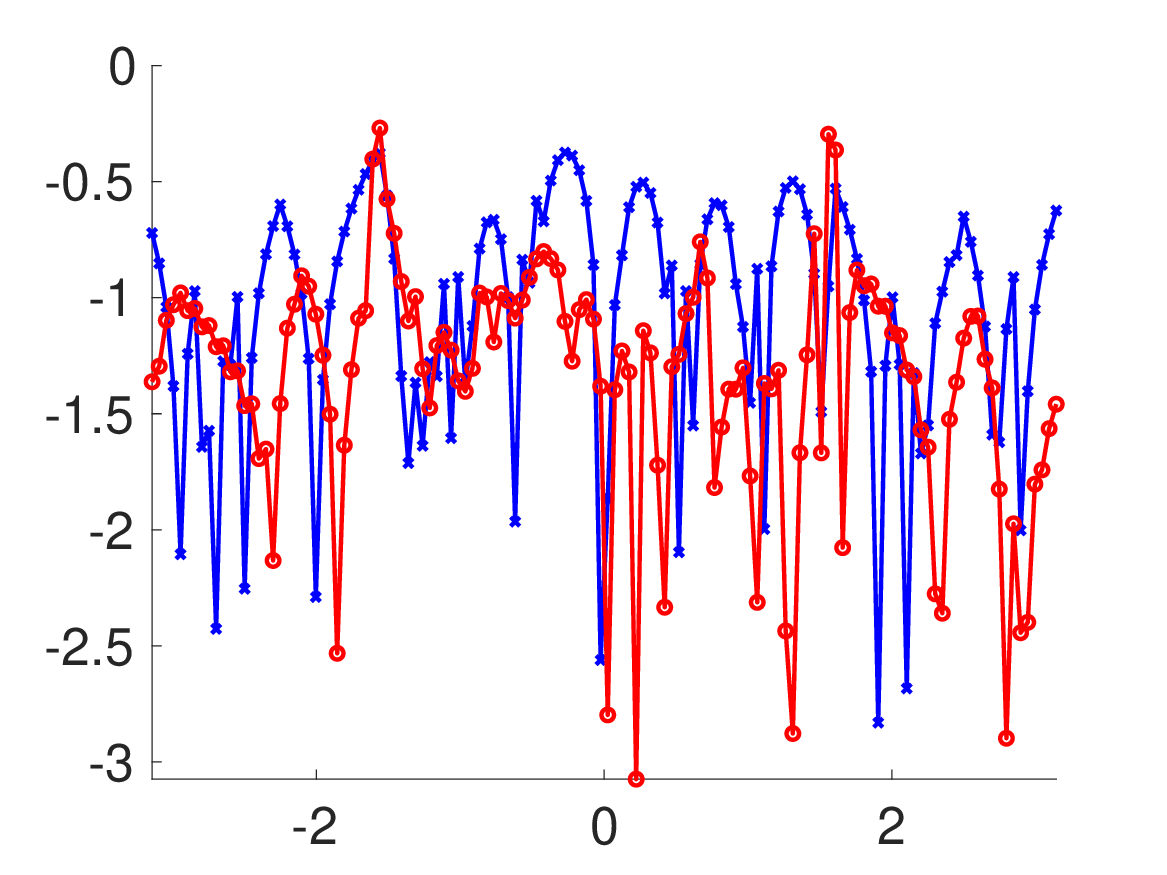}
\end{subfigure}
\begin{subfigure}[b]{.24\textwidth}
\includegraphics[width=\textwidth]{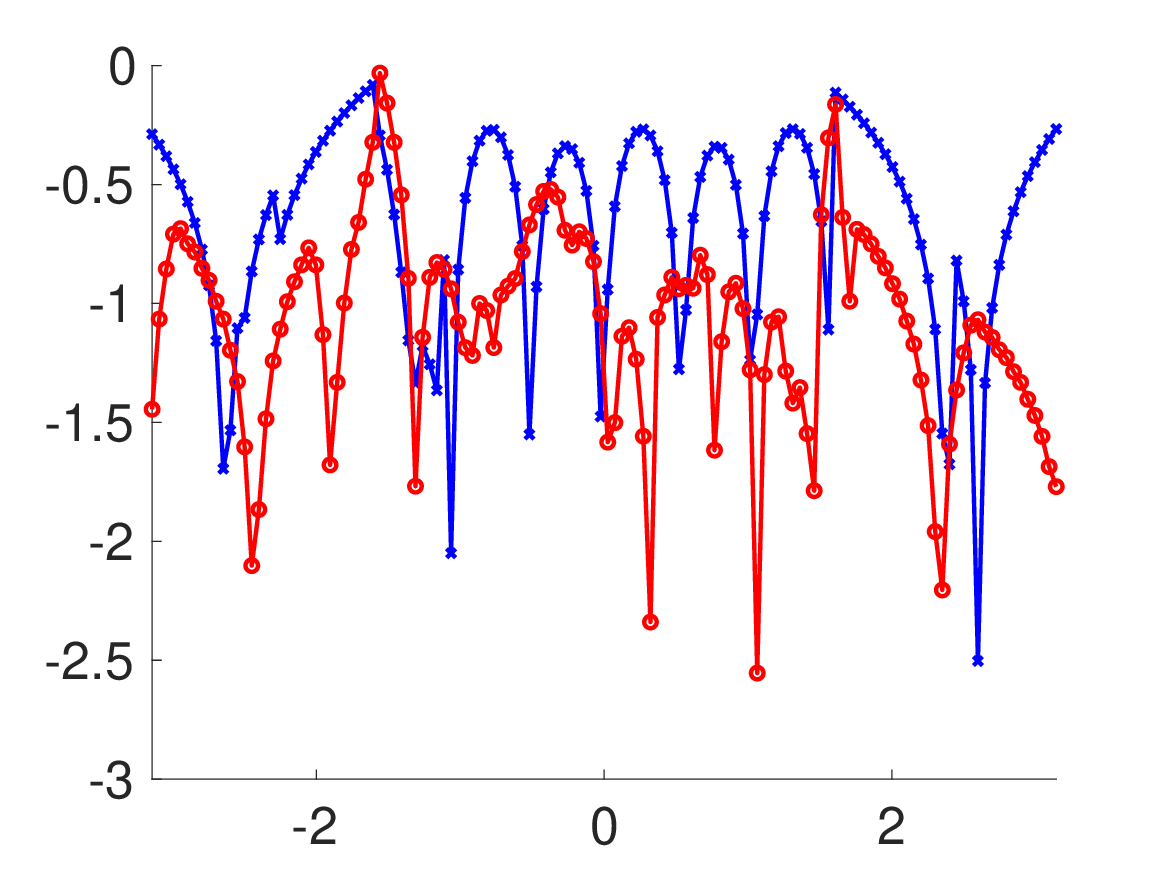}
\end{subfigure}
\begin{subfigure}[b]{.24\textwidth}
\includegraphics[width=\textwidth]{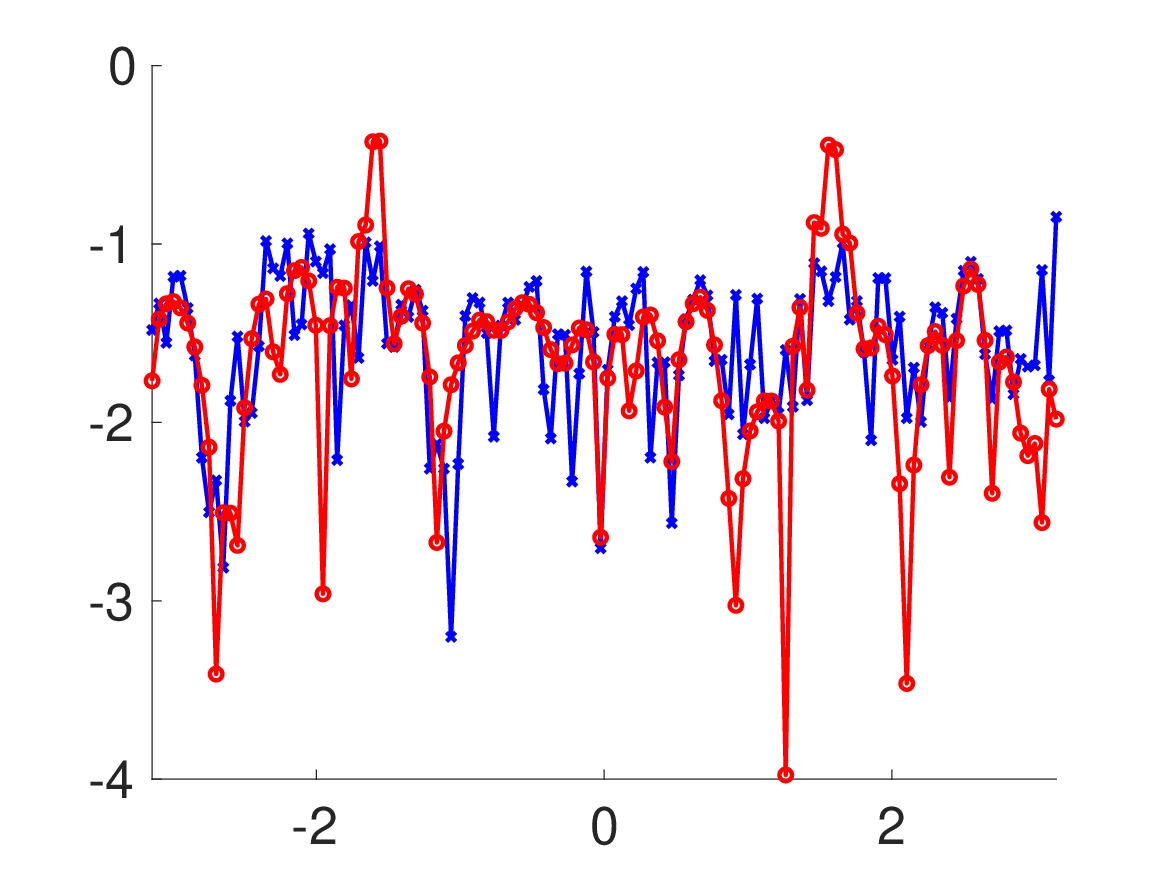}
\end{subfigure}
\\
\begin{subfigure}[b]{.24\textwidth}
\includegraphics[width=\textwidth]{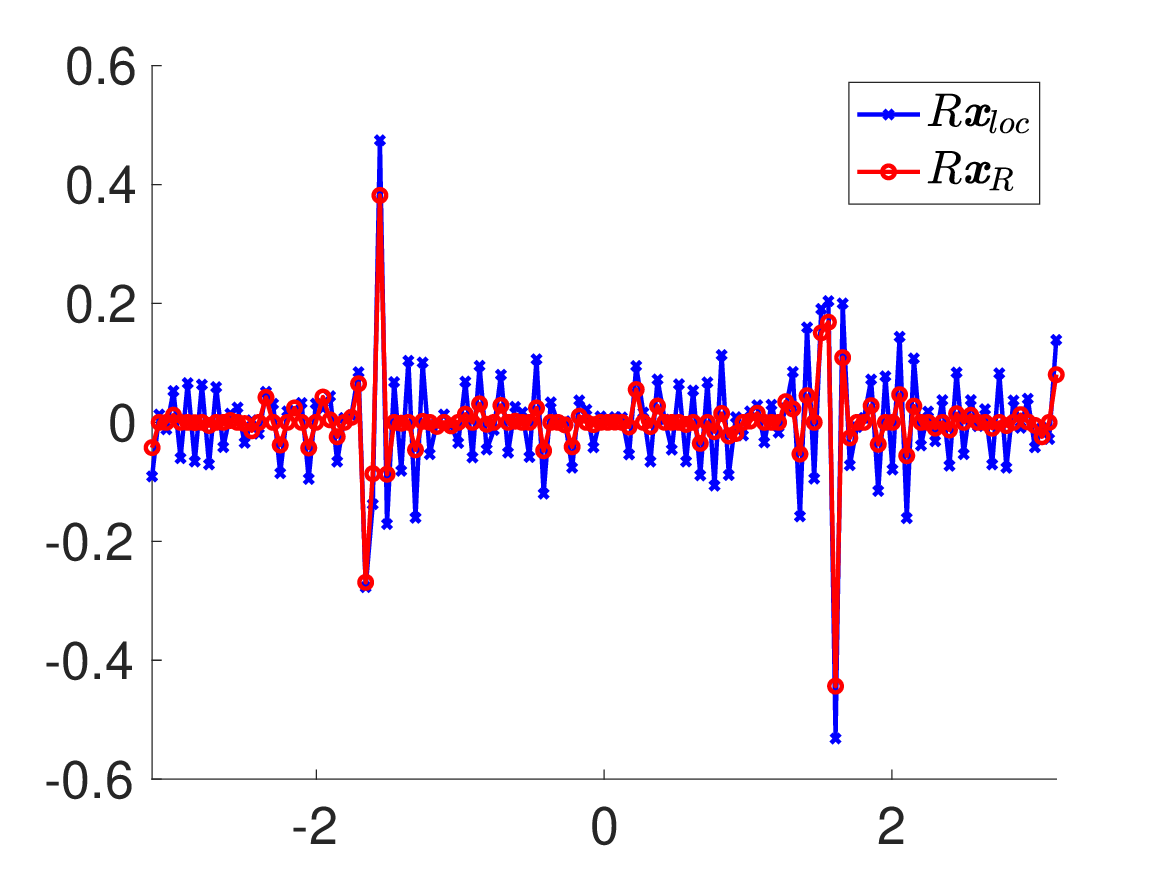}
\end{subfigure}
\begin{subfigure}[b]{.24\textwidth}
\includegraphics[width=\textwidth]{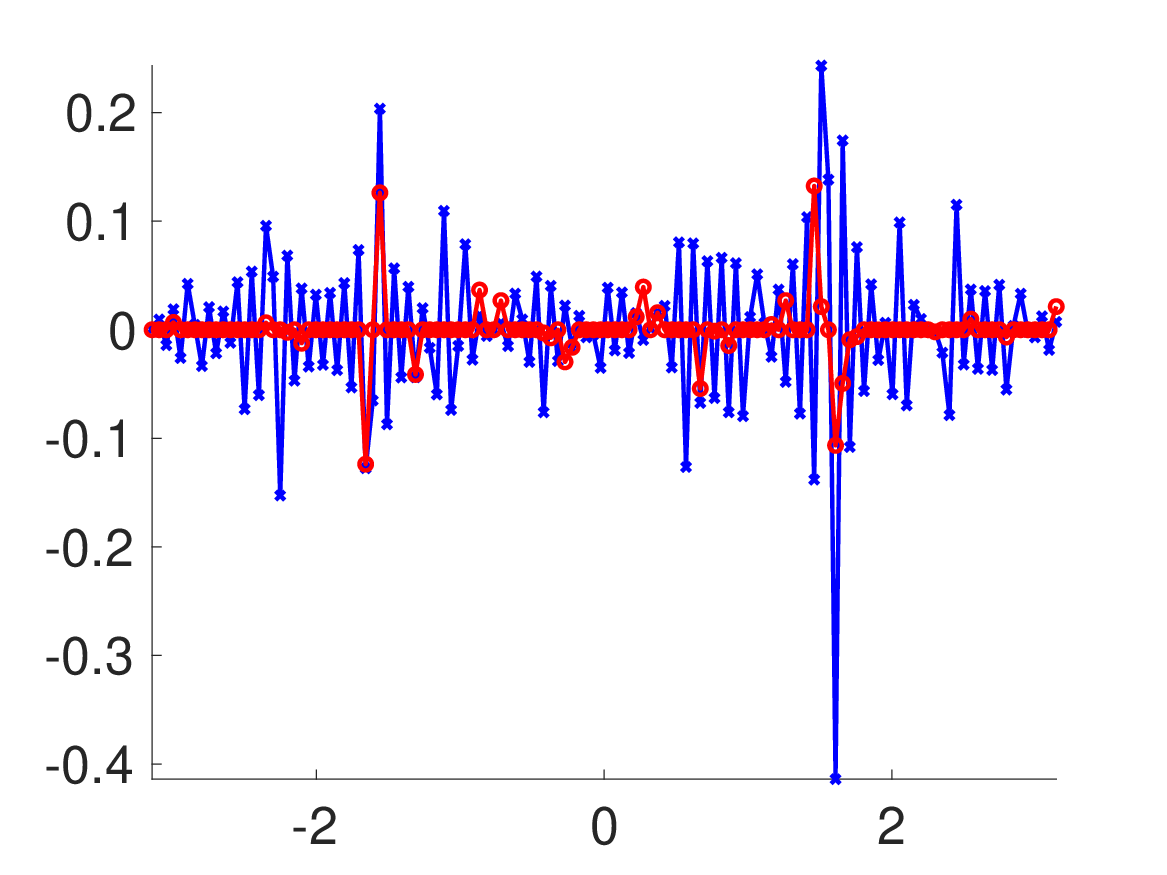}
\end{subfigure}
\begin{subfigure}[b]{.24\textwidth}
\includegraphics[width=\textwidth]{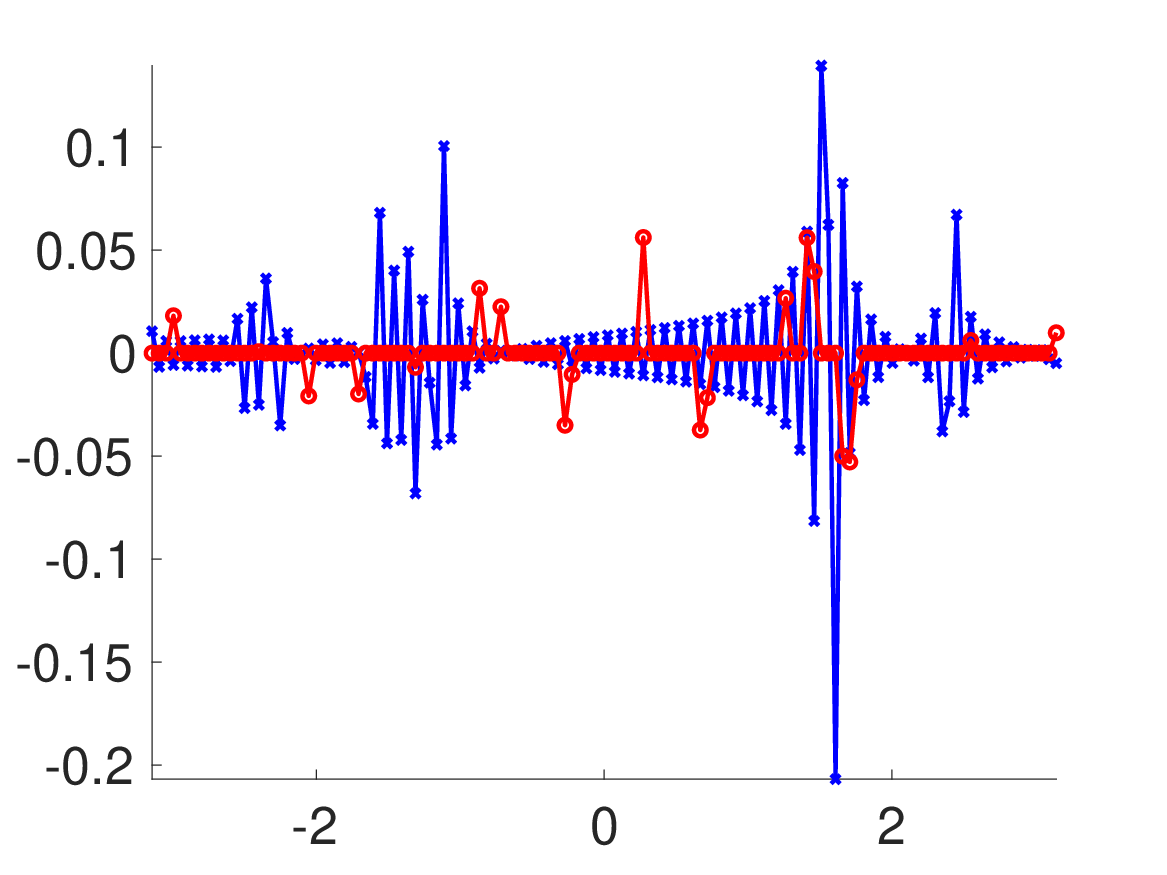}
\end{subfigure}
\begin{subfigure}[b]{.24\textwidth}
\includegraphics[width=\textwidth]{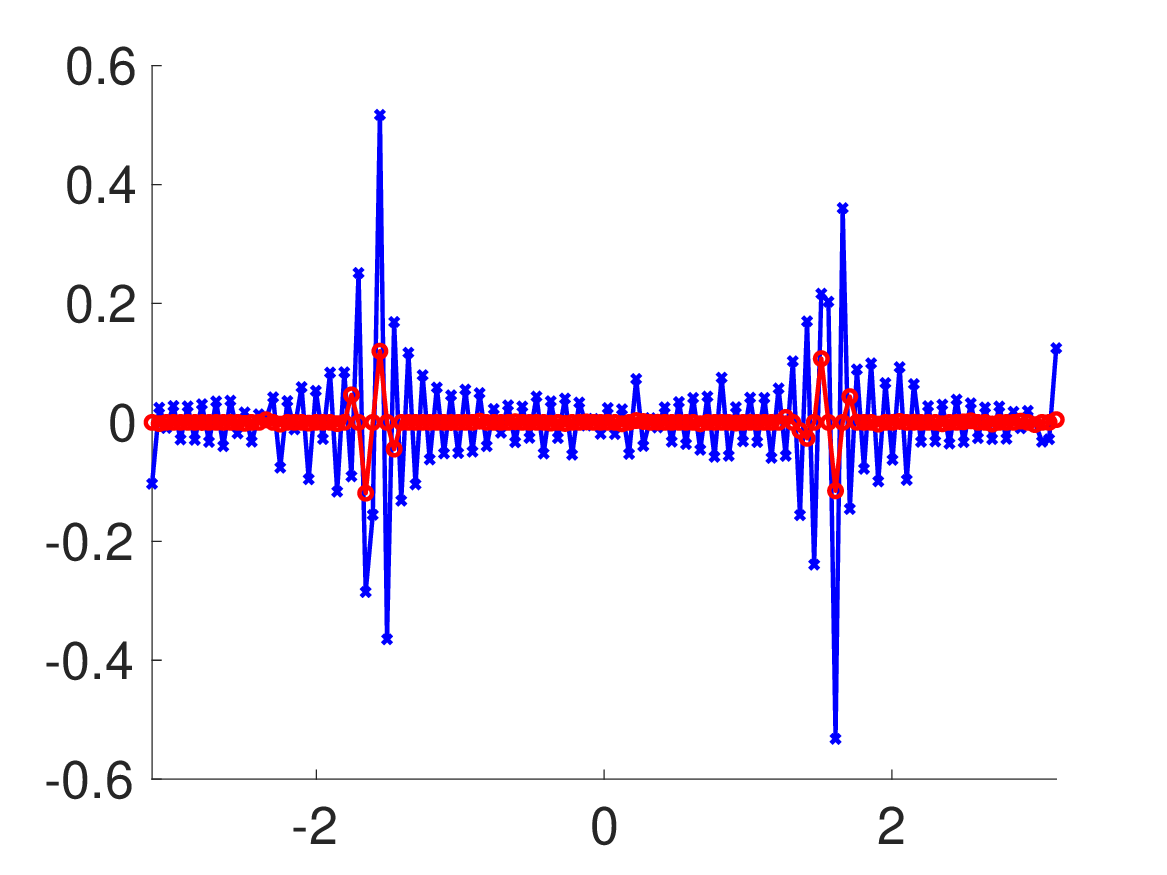}
\end{subfigure}
\caption{(top row) Recovery of \cref{ex:example_denoising} for (left) SNR = 20 dB; (middle left) SNR = 10 dB; (middle right) SNR = 5 dB; (right) SNR = 20 dB in \cref{eq: data_acquisition} with $A = I_n$.  The first three columns show results for $p = 0$ in $T_n^p$ and $R_{n,\zeta}^p$ while the last column uses $p = 1$.  The regularization parameter $\alpha$ is constructed using \eqref{eq:lasso_param} with ${\bm x}_{est}$ calculated as a least squares solution.
(middle row) Corresponding $\log_{10}E^{abs}_j$, $j = 1,\dots,n$,  of the top row solutions given by \cref{eq:absolute_err}. (bottom row) The application of $R_{n,\zeta}^p$ to each corresponding solution in the top row (using the same parameters). 
}
\label{fig:denoising}
\end{figure}

\Cref{fig:denoising}(top) compares the solutions for \cref{ex:example_denoising} obtained by  \cref{eq:lasso_regression} using our new residual transform operator for regularization with the more standard first and third order differencing operators.  \Cref{fig:denoising}(top-left) shows that when SNR = 20 dB, the TV solution ($\LL = T_n^0$) barely differs from the given measurement. The ``staircasing'' artifact becomes more visibly apparent when  SNR = 10 or 5 dB (second and third columns), and it is increasingly difficult to recover the magnitude of $f_1(s)$. By contrast, the residual transform operator of the {\em same} order, that is $R_{n,\zeta}^0$, is able to capture the varying structural behavior.  

Since $f_1(s)$ is not piecewise constant, it is appropriate to use $p = 1$ (third order), and indeed \Cref{fig:denoising}(top-right) demonstrates a more accurate recovery for SNR = 20 dB. However, such  information regarding the variability in the smooth regions of the underlying signal is typically unknown.  We also note that oscillations start to increase with decreasing SNR (not shown here), since the solutions with third order operators are trying to generate piecewise quadratics.   \Cref{fig:denoising}(middle row) displays the spatial errors for the corresponding approximations in the top row.  Finally, \Cref{fig:denoising}(bottom) provides more insight into how each solution ${\bm x}_{loc}$ and ${\bm x}_R$ acts when operated on by $R_{n,\zeta}^p$.  In each case it is clear that $R_{n,\zeta}^p{\bm x}_R$ is more sparse than  $R_{n,\zeta}^p{\bm x}_{loc}$.

\begin{figure}[h!]
    \centering
\begin{subfigure}[b]{.4\textwidth}
\includegraphics[width=\textwidth]{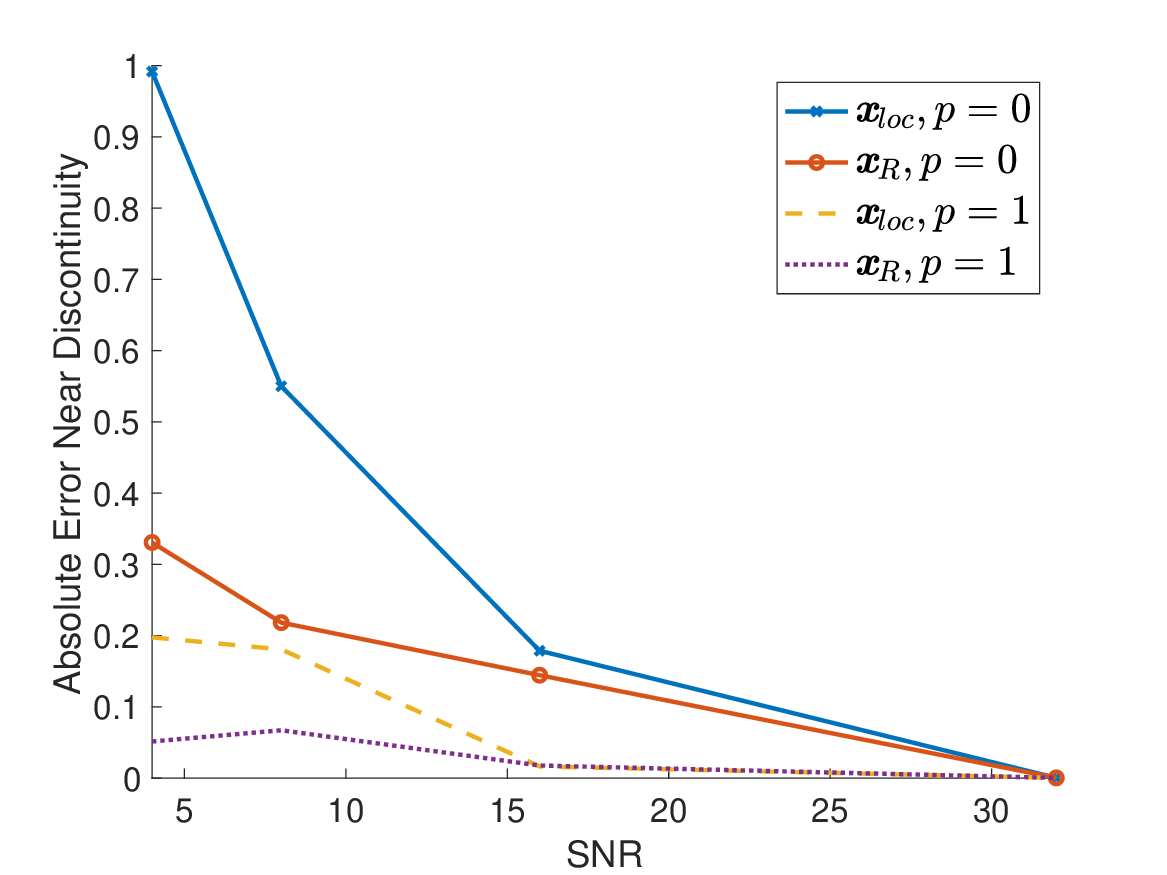}
\label{fig:snr_edge}
\end{subfigure}
\begin{subfigure}[b]{.4\textwidth}
\includegraphics[width=\textwidth]{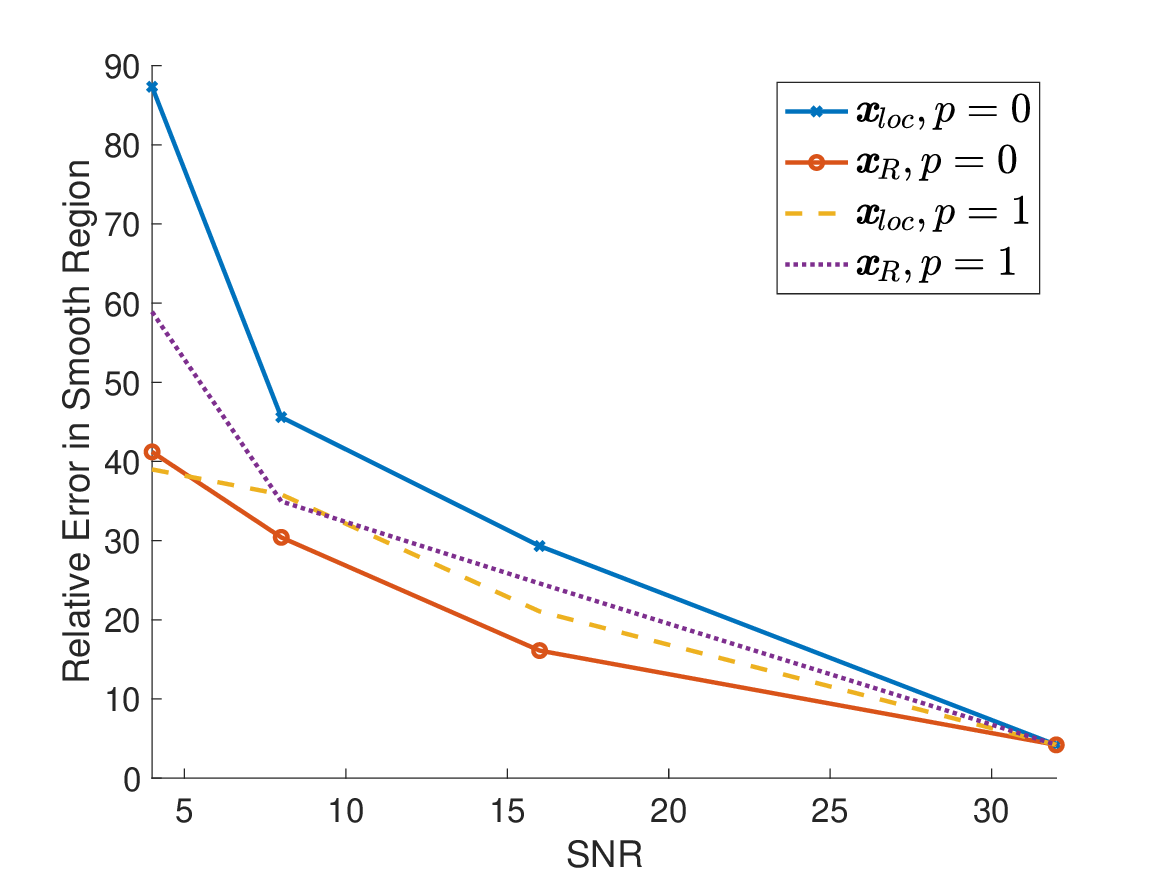}
\label{fig:snr_smooth}
\end{subfigure}
\caption{Error in the recovery of \cref{ex:example_denoising} as a function of $\text{SNR}=2^b$, $b=2,\dots,5$, in \eqref{eq:SNR}. (left) $E^{abs}(s_{n/4-2})$  given by  \eqref{eq:absolute_err} with  $s_{n/4-2} \approx-1.7$. (right) 
$E^{rel}$ given by \cref{eq:relative_err} with $j_{min} = 44$ and $j_{max} = 85$ so that  $s_j \in\left[-\frac{\pi}{3},\frac{\pi}{3}\right]$ (a smooth region of $f_1(s)$).}
\label{fig:snr_test}
\end{figure}

\Cref{fig:snr_test}(left) compares the absolute error using different choices of sparse transform operators near a point of discontinuity in \cref{ex:example_denoising} as a function of SNR $\in [2,32]$, while \Cref{fig:snr_test}(right) compares the relative error in each case within a smooth region. As already noted, the variability of the underlying signal suggests that  third order differencing ($\mathcal L = T_n^1$ in \cref{eq:lasso_regression}) would be appropriate to obtain sparsity, and indeed ${\bm x}_{loc}$ when $p = 1$ is more accurate than ${\bm x}_{loc}$ when $p=0$ for all tested levels of SNR in smooth regions as well as near a discontinuity.  Our new residual transform operator, $R_{n,\zeta}^p$, for $p = 0$ {\em  or} $p = 1$ clearly yields more accurate solutions than the standard TV operator does.  Interestingly, the new residual operator with $p = 0$ recovers the most accurate solution in smooth regions while the results using $p = 1$ are comparable for SNR $> 8$.  Most importantly we observe that the residual transform operator is robust with respect to noise and transform matrix order $2p+1$.

\subsection{Deblurring}
\label{sec:deblurring}
For the deblurring problem,  we apply $A=\tilde A$ in \eqref{eq: data_acquisition} with $\gamma = .01$ and $.05$.  Once again we compare the results of \cref{eq:lasso_regression} with $\mathcal{L} = T_{n}^0$ and $\mathcal{L} = R_{n,\zeta}^0$, $\zeta = \frac{1}{4}$, for $f_1(s)$ in \cref{ex:example_denoising}. Since a higher-order operator is often more appropriate for piecewise polynomial function $f_1(s)$ the corresponding results are displayed in the rightmost column for $p=1$ under blurring condition $\gamma=.01$.
As there is no additive Gaussian noise in the measurement, the regularization parameter $\alpha$ is heuristically set to 0.1 and 0.3 for $\mathcal{L} = T_{n}^{p}$ and $\mathcal{L} = R_{n,\zeta}^{p}$ respectively, $p=0,1$. 

\begin{remark}[Regularization parameter]
\label{rem:regparam} Many methods exist for choosing the regularization parameter $\alpha$ in \cref{eq:lasso_regression} for the deblurring problem.  As this is not the focus of our investigation, we simply note that  using our new residual transform operator in place of a more standard sparse transform operator intuitively suggests that the regularization term should be more heavily penalized in the solution.  This is confirmed in our numerical experiments.  Noise and ill-conditioning of the forward operator should also be considered in this parameter selection, and will be addressed in future work.
\end{remark}

\begin{figure}[h!]
    \centering
\begin{subfigure}[b]{.32\textwidth}
\includegraphics[width=\textwidth]{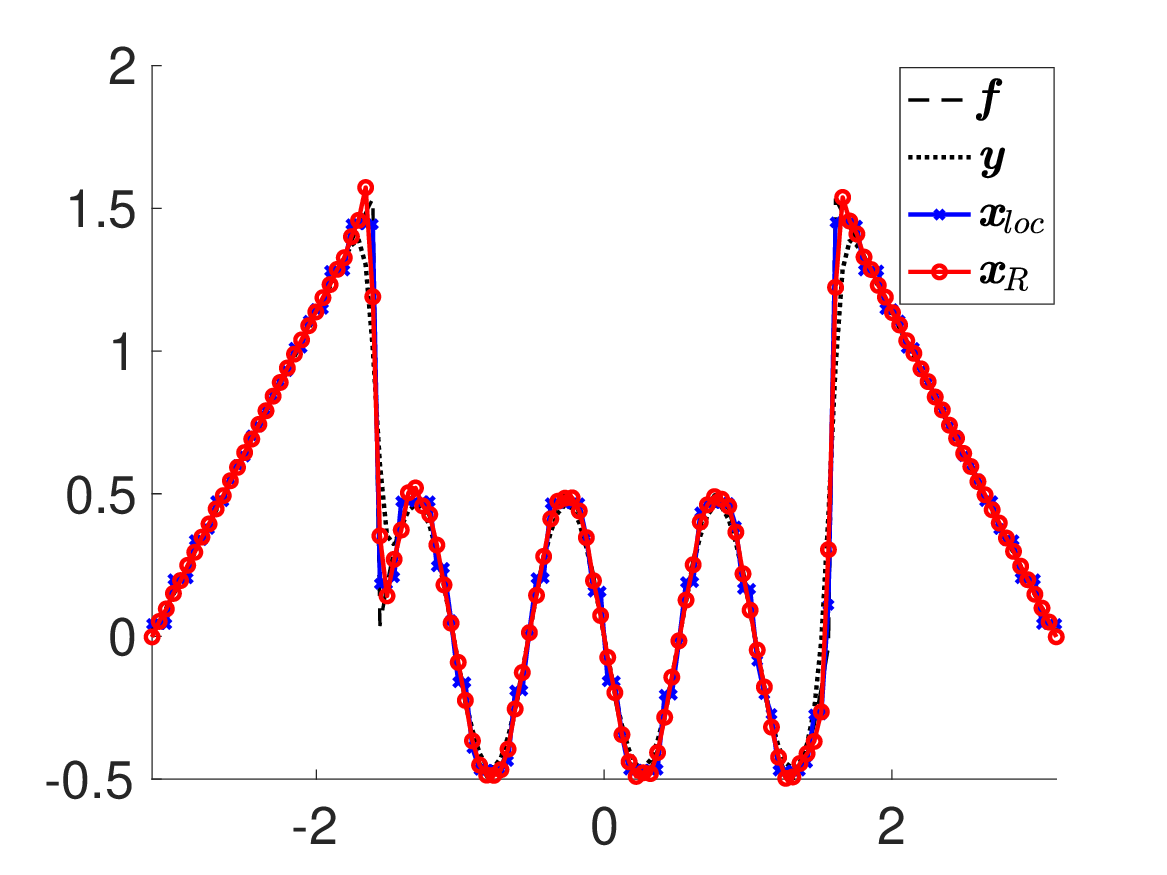}
\end{subfigure}
\begin{subfigure}[b]{.32\textwidth}
\includegraphics[width=\textwidth]{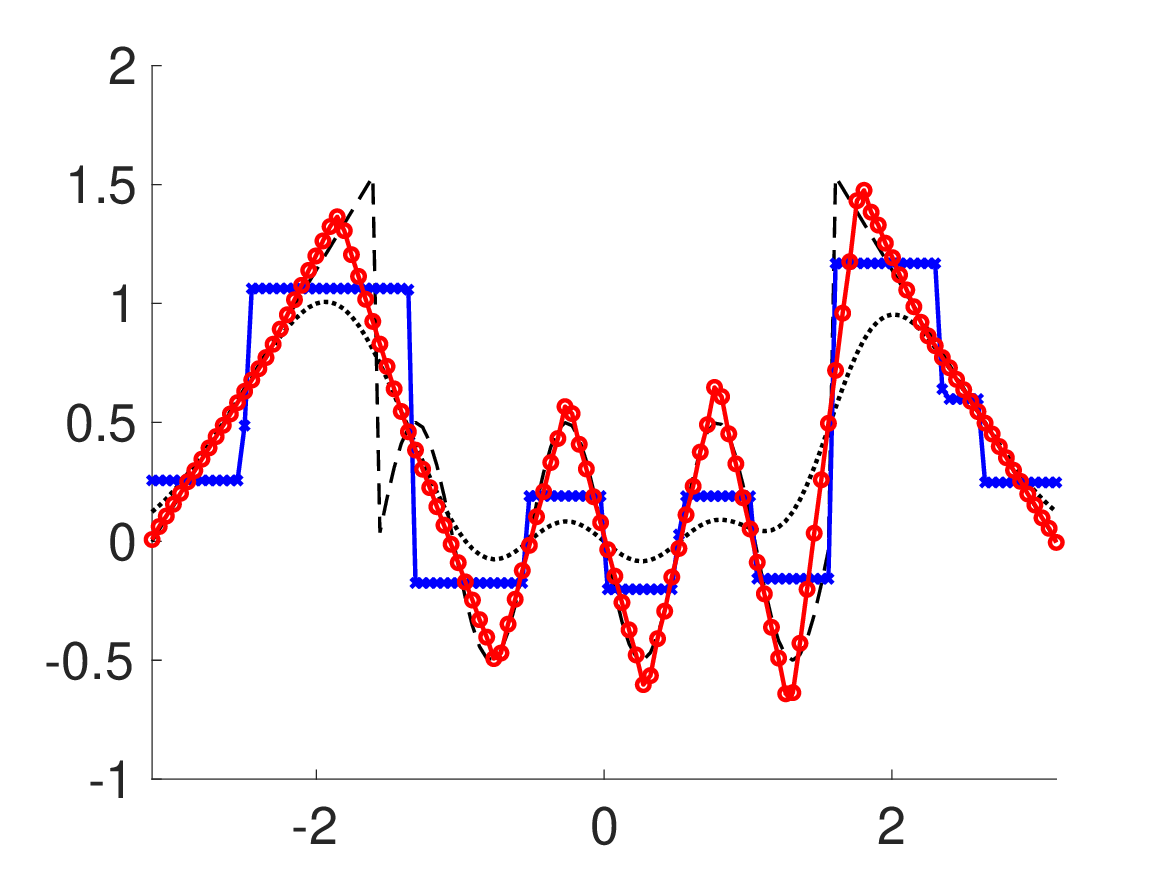}
\end{subfigure}
\begin{subfigure}[b]{.32\textwidth}
\includegraphics[width=\textwidth]{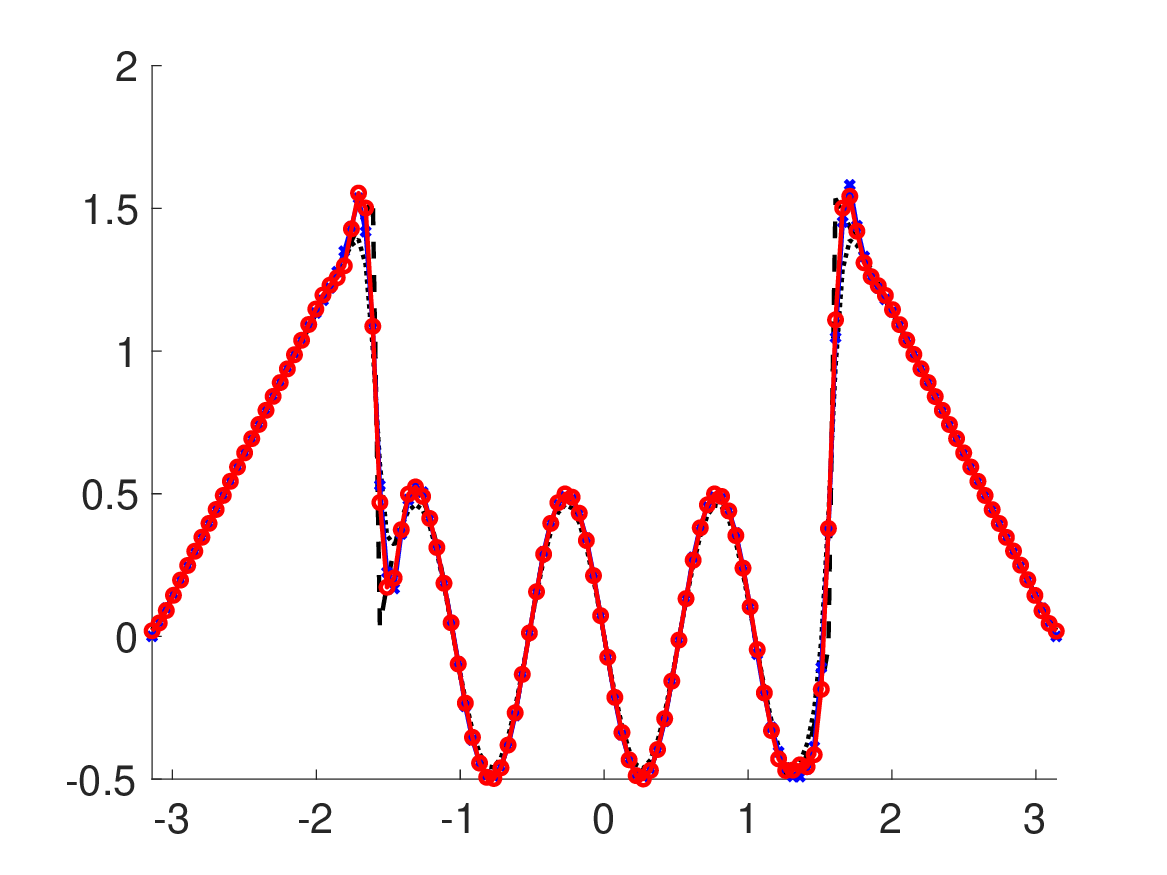}
\end{subfigure}\\
\begin{subfigure}[b]{.32\textwidth}
\includegraphics[width=\textwidth]{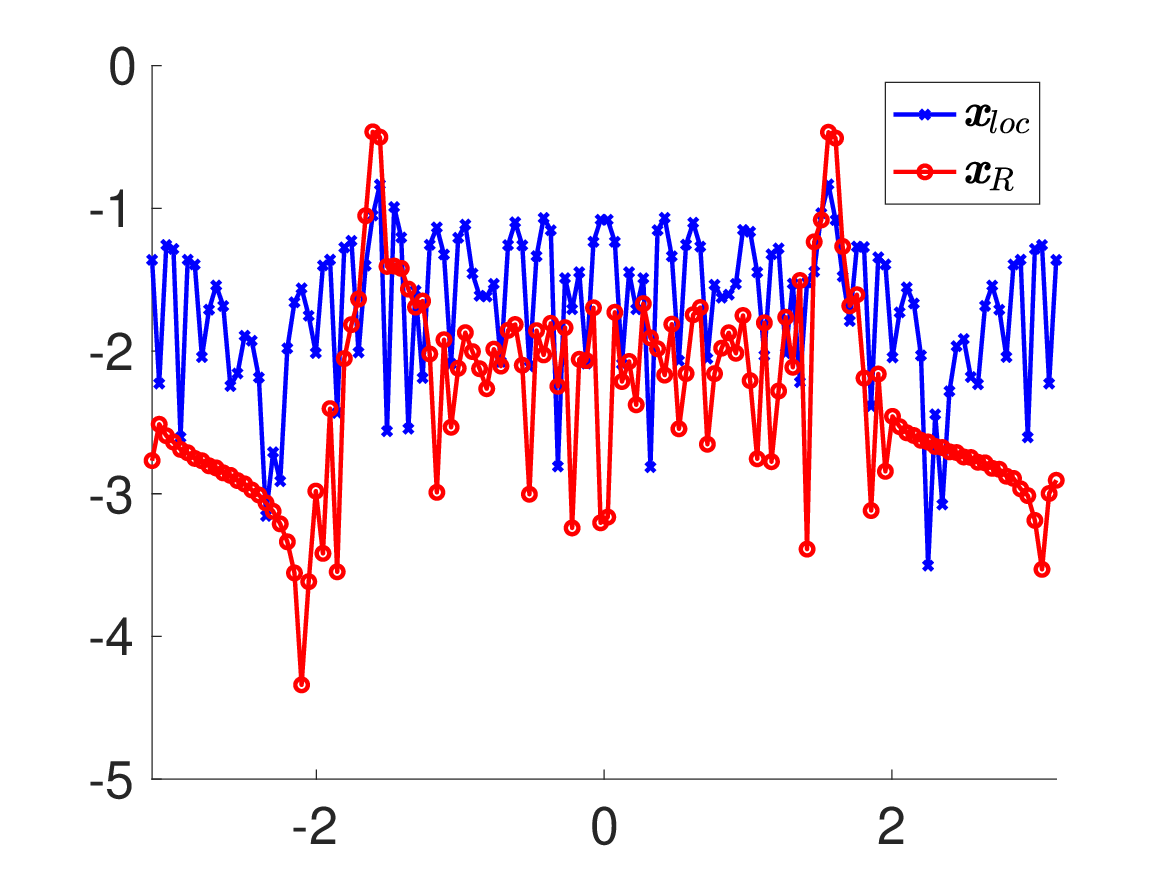}
\end{subfigure}
\begin{subfigure}[b]{.32\textwidth}
\includegraphics[width=\textwidth]{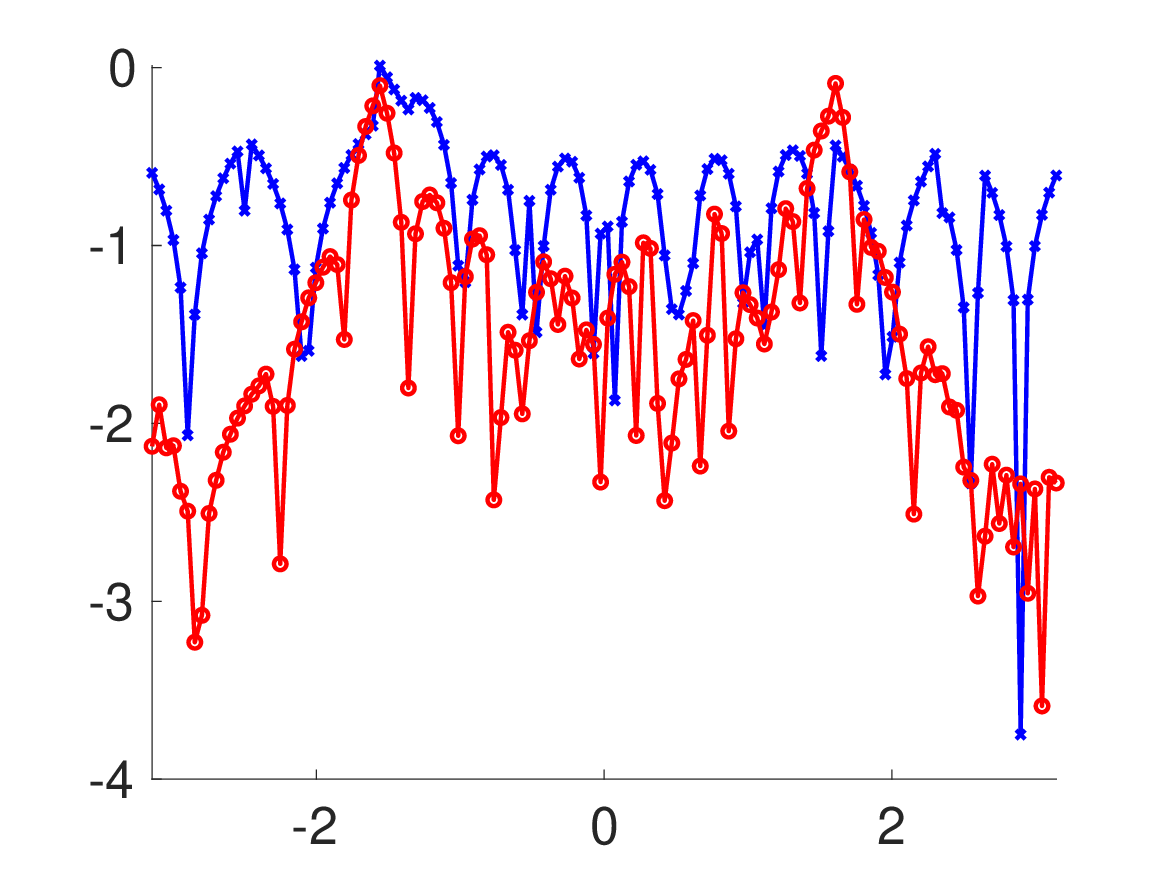}
\end{subfigure}
\begin{subfigure}[b]{.32\textwidth}
\includegraphics[width=\textwidth]{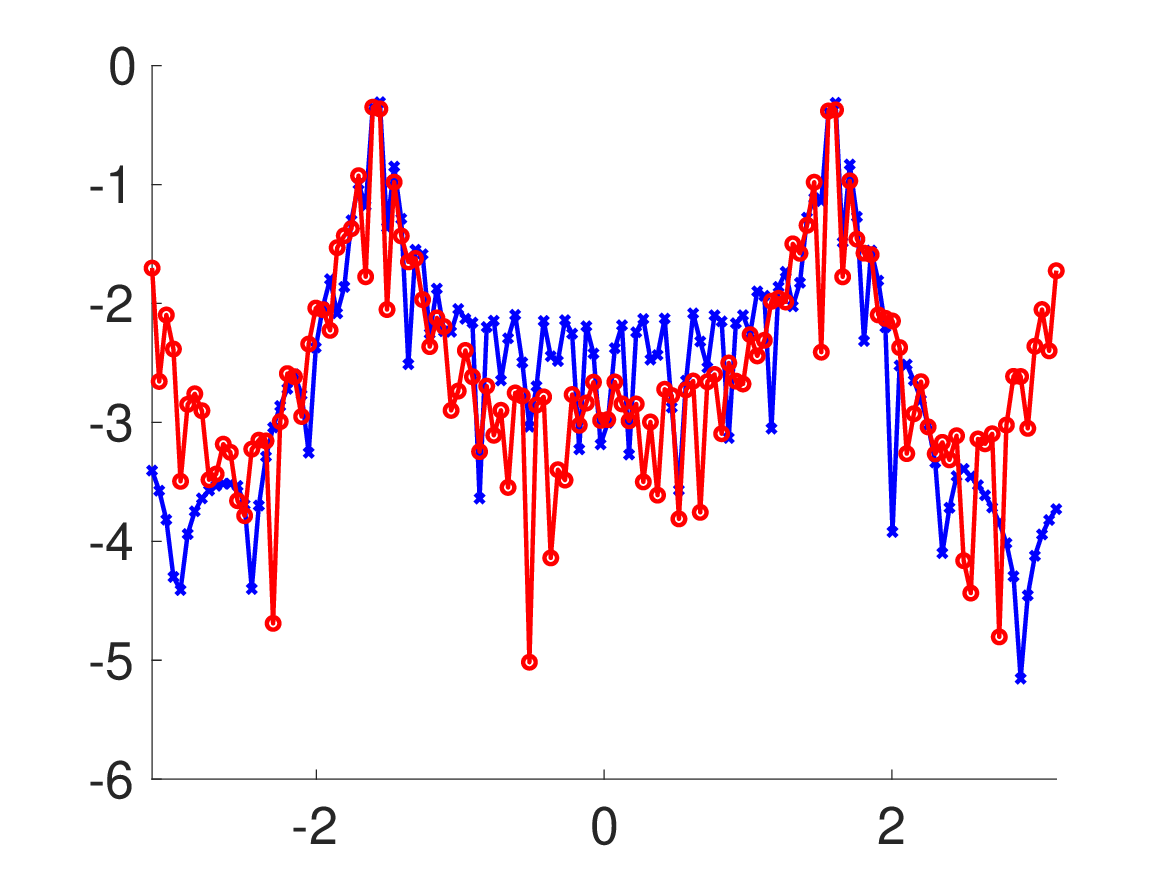}
\end{subfigure}
\caption{(top row) Recovery of \Cref{ex:example_denoising} from blurred data with  (left) $\gamma=0.01$ (middle) $\gamma=0.05$ and  (right) $\gamma= 0.01$ in \cref{eq:bluroperator} with $A=\tilde A$. The left and middle columns show  results for $T_n^p$ and $R_{n,\zeta}^p$  with $p = 0$ while the right column uses $p = 1$.  The regularization parameters are respectively $\alpha = .1$ and $\alpha = .3$ for ${\bm x}_{loc}$  and $\bm x_{R}$.  (bottom row) Corresponding $\log_{10}E^{abs}_j$, $j = 1,\dots,n$,  of the top row solutions given by \cref{eq:absolute_err}.}
\label{fig:deblurring}
\end{figure}

\begin{figure}[h!]
    \centering
\begin{subfigure}[b]{.4\textwidth}
\includegraphics[width=\textwidth]{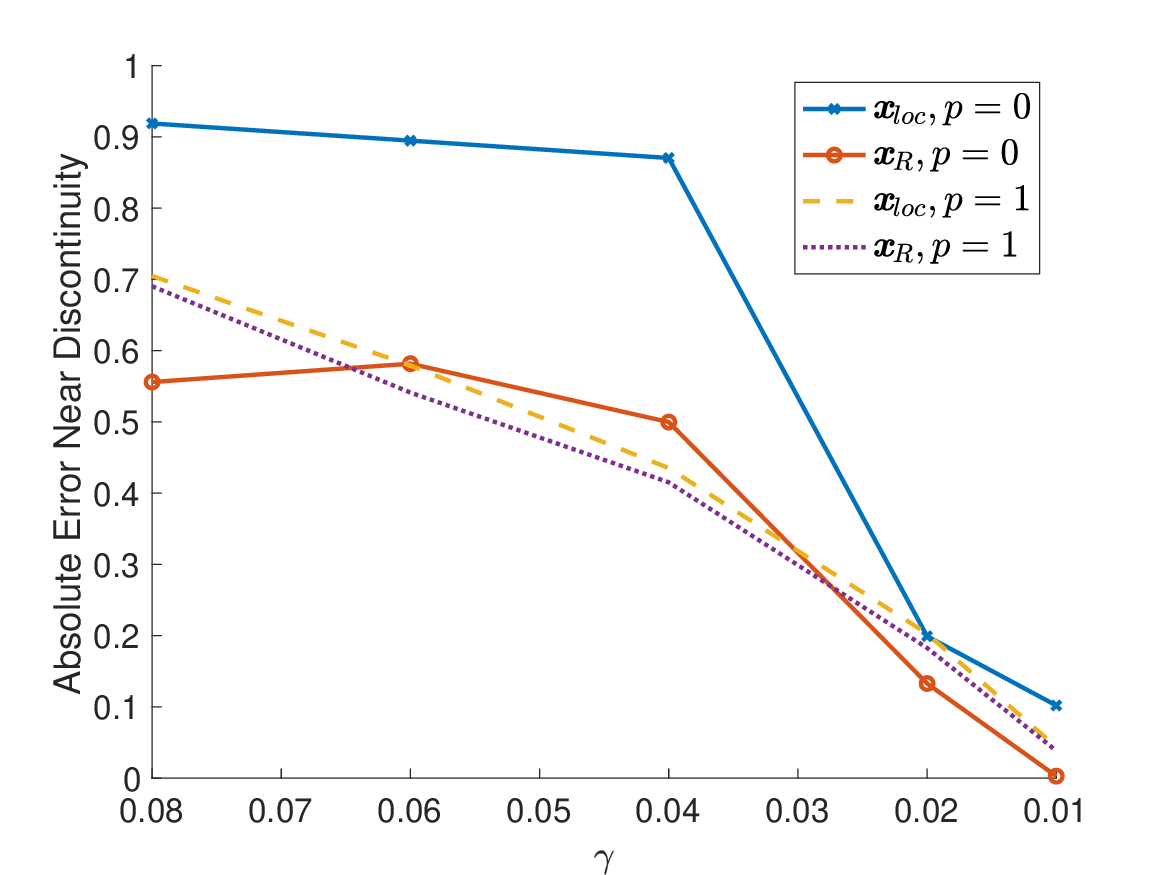}
\label{fig:blur_edge}
\end{subfigure}
\begin{subfigure}[b]{.4\textwidth}
\includegraphics[width=\textwidth]{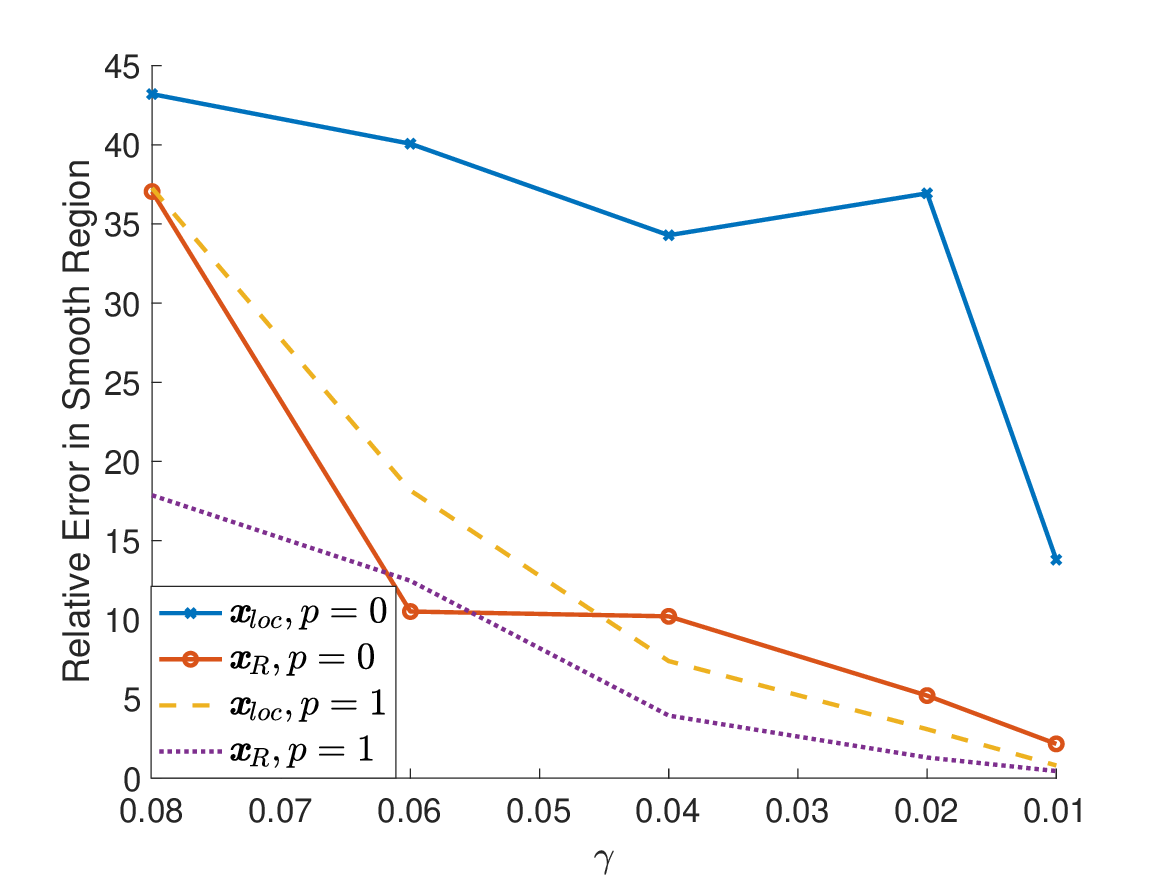}
\label{fig:blur_smooth}
\end{subfigure}
\caption{Error in the recovery of \cref{ex:example_denoising} for $\gamma \in [.01,.08]$ in \eqref{eq:bluroperator}. 
(left) $E^{abs}(s_{n/4+2})$ given by  \eqref{eq:absolute_err} with  $s_{n/4+2} \approx-1.5$. (right) 
$E^{rel}$ given by \cref{eq:relative_err} with $j_{min} = 44$ and $j_{max} = 85$ so that  $s_j \in\left[-\frac{\pi}{3},\frac{\pi}{3}\right]$ (a smooth region of $f_1(s)$).}
\label{fig:blur_test}
\end{figure}

\Cref{fig:deblurring}(top) compares the solutions for \cref{ex:example_denoising} obtained by  \cref{eq:lasso_regression}, contrasting our new residual transform operator $\LL=R^p_{n,\zeta}$ with the more standard first and third order differencing operators  $\LL=T^p_n$.  The third order differencing operator ($p=1$) is included for comparison in the right column for $\gamma=0.01$ as a better fit for a piecewise polynomial function $f_1(s)$. 
The proposed residual transform operator, $\LL=R^1_{n,\zeta}$, yields more accurate solutions than the standard third-order local differencing operator, $\LL=T^1_n$. While not displayed here, this superior performance is consistent in cases with significant blurring, as demonstrated in  \cref{fig:blur_test}.
Since the variability within the smooth regions of an underlying function is often unknown a priori, similar to the denoising problem, applying an unsuitable regularizer can lead to significant errors in the deblurring scenario. This is evident in the left and the middle columns for both cases of $\gamma=0.01$ and $\gamma=0.05$, where using the first order local differencing operator ($\LL = T_n^0$) results in considerable error and a pronounced ``staircasing'' artifact. By contrast, the residual transform operator {\em of the same order} successfully mitigates the ``staircasing'' effect by accurately capturing the local variability within smooth regions. The bottom row displays the spatial errors for the corresponding approximations in the top row. Notably, the error plot for $p=1$(right) confirms that the solution ${\bm x}_R$ continues to outperform ${\bm x}_{loc}$.

\Cref{fig:blur_test} then compares the errors resulting from use of  each regularization operator as a function of blurring parameter $\gamma\in[.01, .08]$ in \cref{eq:bluroperator}. As in \Cref{fig:snr_test}, we look at the absolute error near a discontinuity (left) and the relative error over a smooth region (right). As expected  for a signal with complex internal variability, the third order differencing ($\mathcal L = T_n^1$ in \cref{eq:lasso_regression}) produces a more accurate reconstruction, ${\bm x}_{loc}$, than the first order operator ($p=0$). This superior performance  is consistent across all tested levels of $\gamma$ and is observed in both smooth regions as well as near a discontinuity. Our new residual transform operator, $R_{n,\zeta}^p$, for $p = 0$ {\em  or} $p = 1$ consistently yields more accurate solutions than the standard TV operator ($\LL=T_n^0$) does. Specifically, the new residual operator ($\LL=R_{n,\zeta}^1$) recovers the most accurate solution in smooth regions when the blurring parameter is low ($\gamma<0.06$). Overall, the results in \cref{fig:blur_test}, similar to the previous SNR analysis, confirm that our residual transform operator is robust to varying levels of blur and effective across transform matrix order $2p+1$.

\subsection{Undersampling}
\label{sec:undersample}

\begin{figure}[h!]
    \centering
\begin{subfigure}[b]{.24\textwidth}
\includegraphics[width=\textwidth]{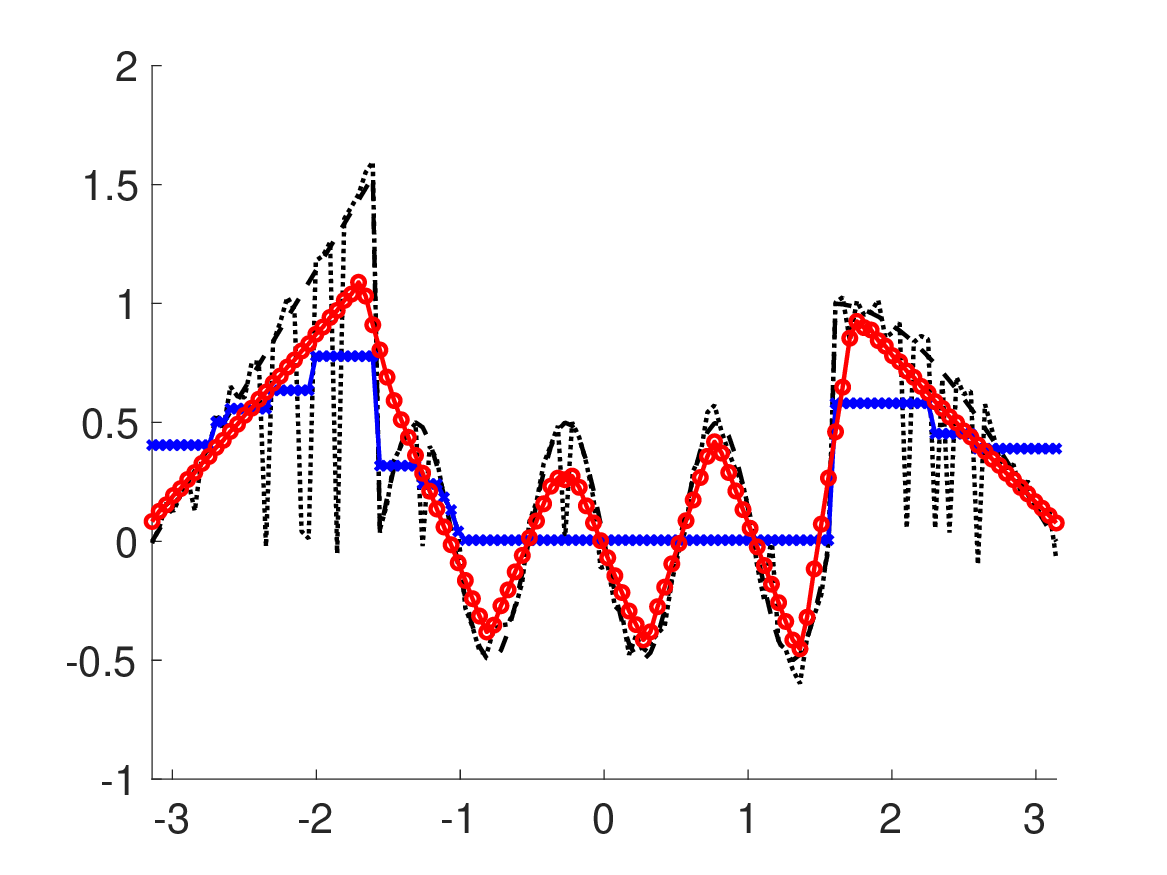}
\end{subfigure}
\begin{subfigure}[b]{.24\textwidth}
\includegraphics[width=\textwidth]{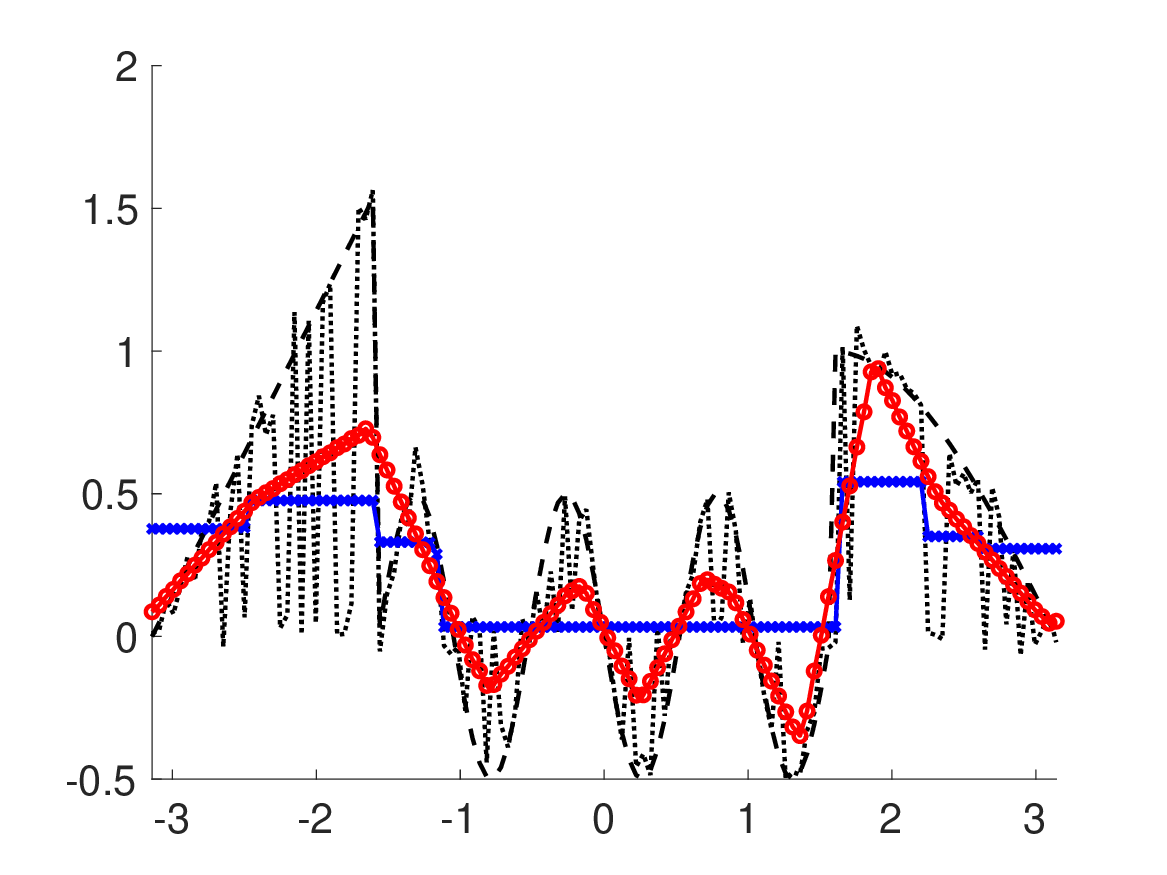}
\end{subfigure}
\begin{subfigure}[b]{.24\textwidth}
\includegraphics[width=\textwidth]{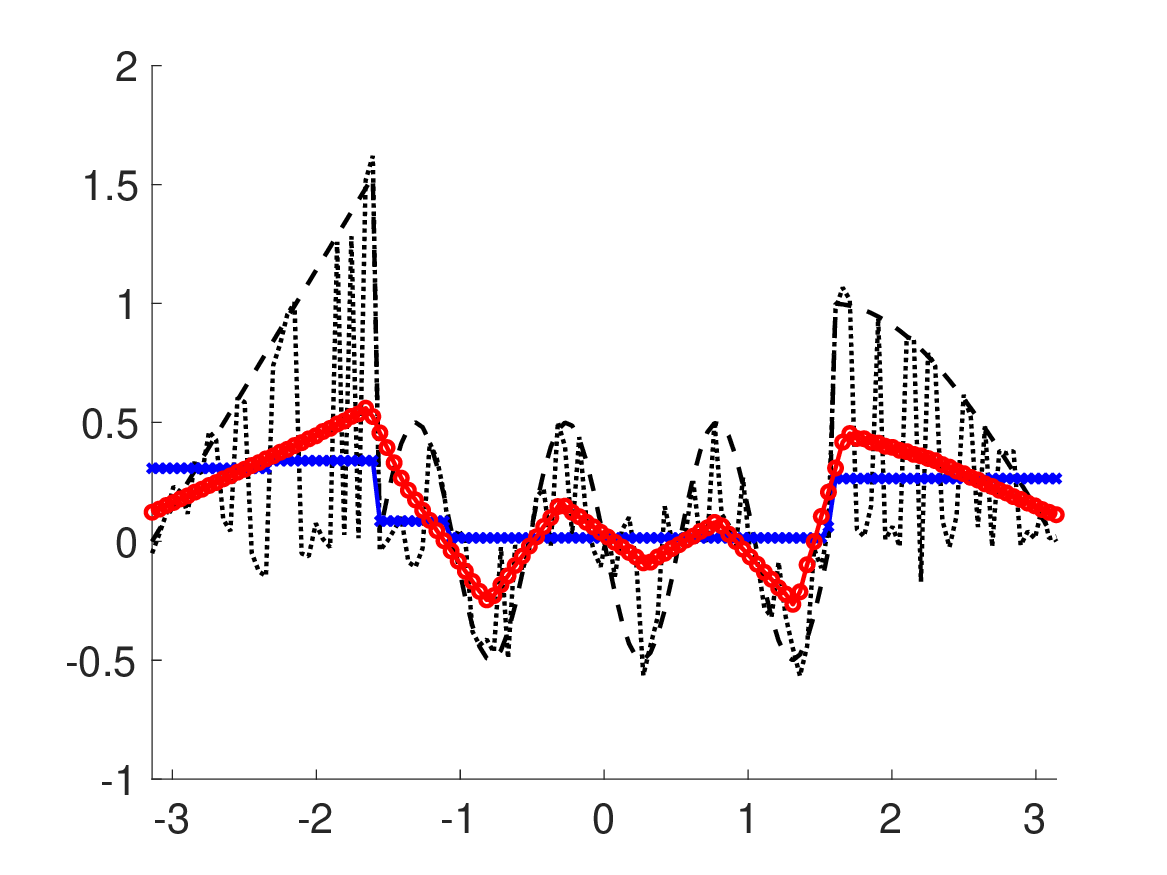}
\end{subfigure}
\begin{subfigure}[b]{.24\textwidth}
\includegraphics[width=\textwidth]{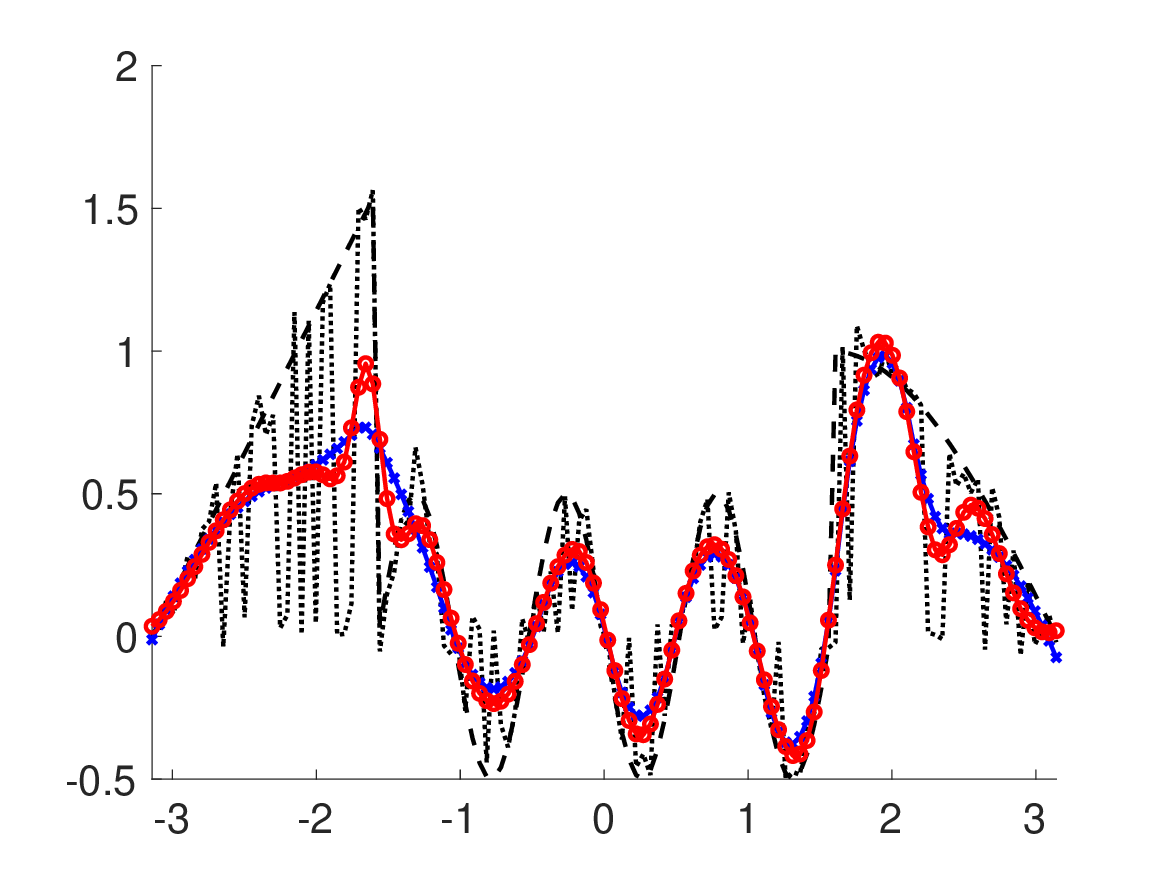}
\end{subfigure}\\
\begin{subfigure}[b]{.24\textwidth}
\includegraphics[width=\textwidth]{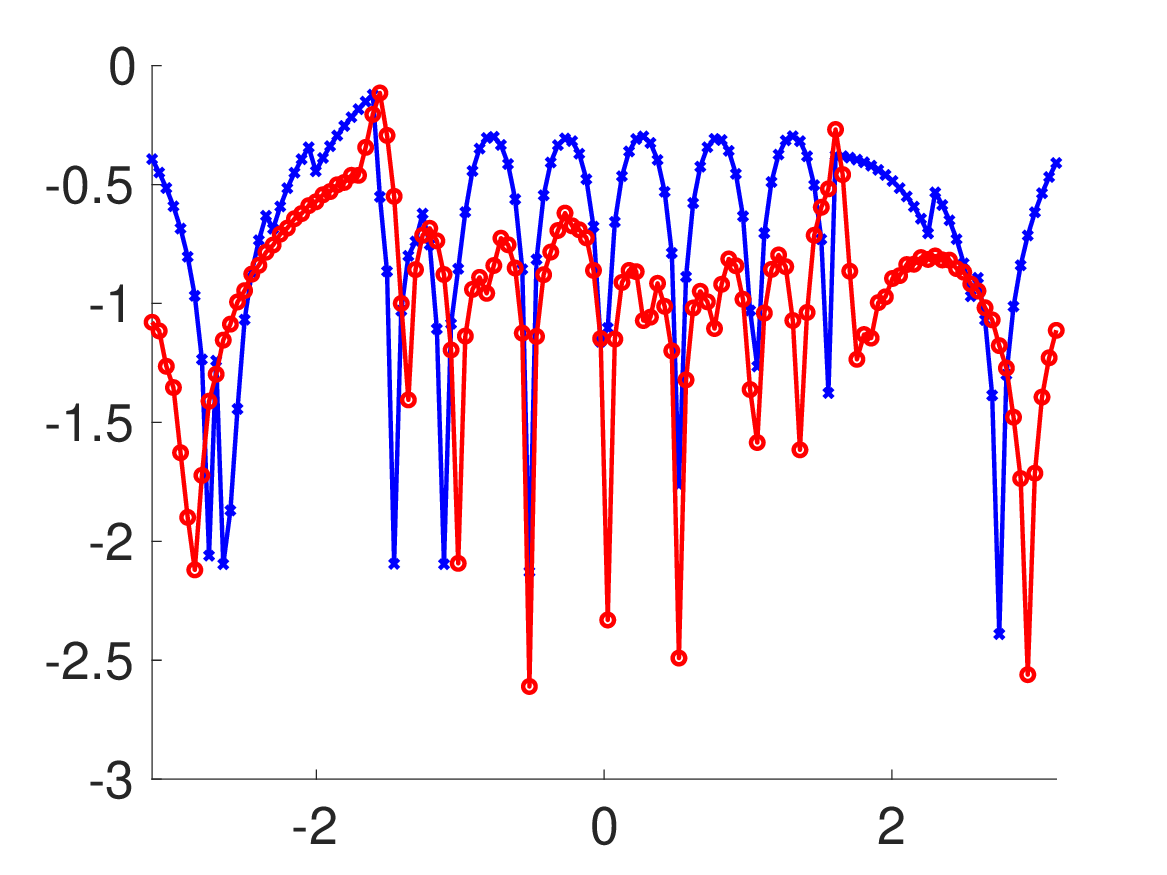}
\end{subfigure}
\begin{subfigure}[b]{.24\textwidth}
\includegraphics[width=\textwidth]{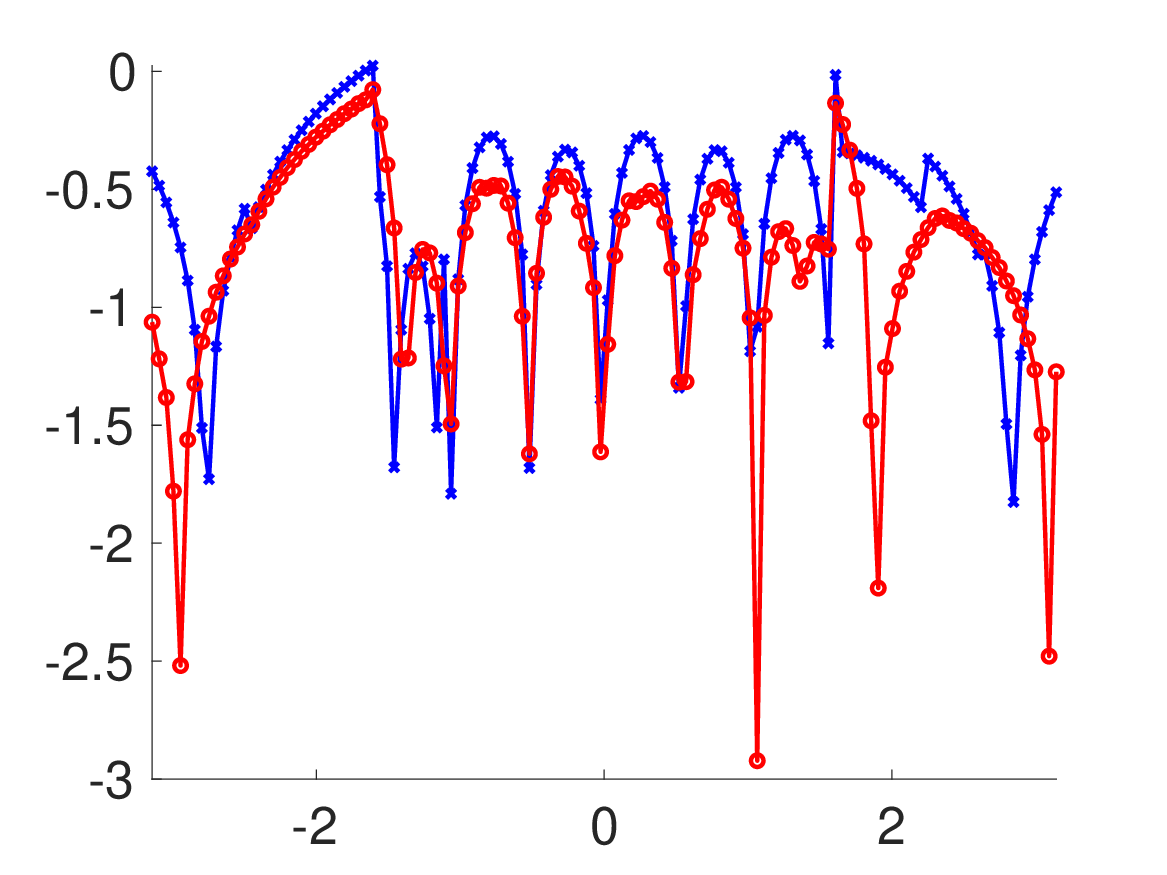}
\end{subfigure}
\begin{subfigure}[b]{.24\textwidth}
\includegraphics[width=\textwidth]{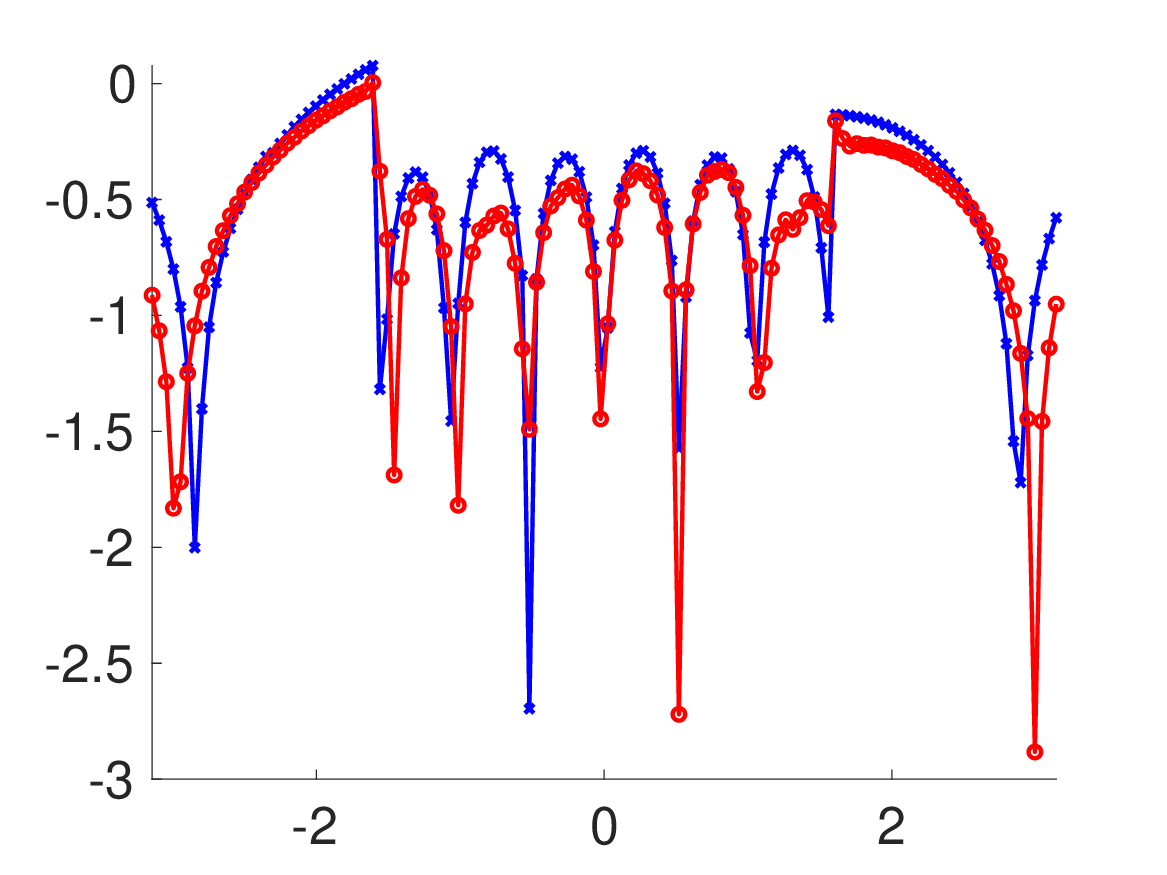}
\end{subfigure}
\begin{subfigure}[b]{.24\textwidth}
\includegraphics[width=\textwidth]{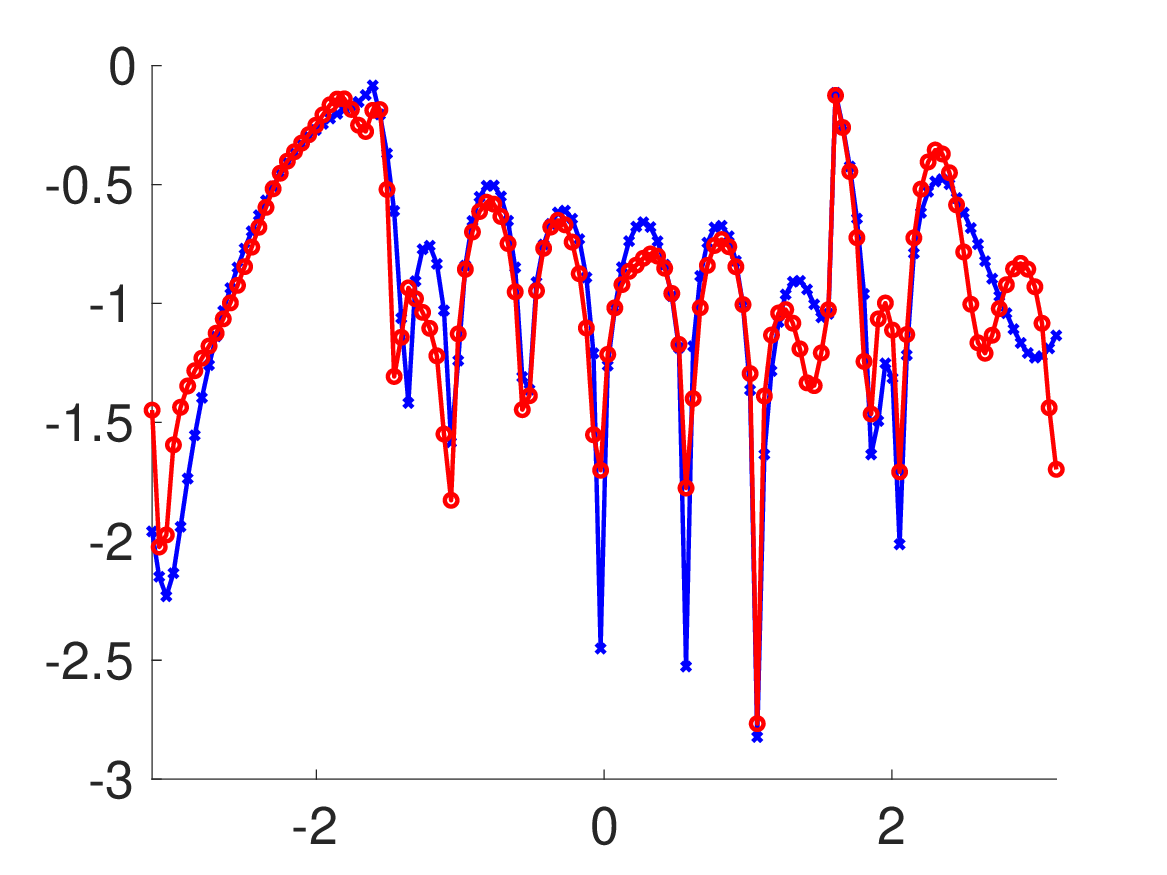}
\end{subfigure}
\caption{(top row) Recovery of \Cref{ex:example_denoising} by \cref{eq:lasso_regression} when $A=I_n$ from undersampled data with SNR = 20 dB in \cref{eq: data_acquisition}, where $A=HI_n$ with $h = rn$ randomly selected zeroes for (left) $r=0.1$; (middle left) $r=0.3$; (middle right); (right) $r=0.3$. The first three columns show  results for $p = 0$ in $T_n^p$ and $R_{n,\zeta}^p$ while the right column uses $p = 1$.  The regularization parameter $\alpha$ is determined from  \eqref{eq:lasso_param}  with ${\bm x}_{est}$ calculated as a least squares solution. (bottom row) Corresponding $\log_{10}E^{abs}_j$, $j = 1,\dots,n$,  of the top row solutions given by \cref{eq:absolute_err}.}
\label{fig:undersampling}
\end{figure}

Varying levels of undersampling for collected measurements is accomplished by applying an operator $H$ to $A = I_n$ in \eqref{eq: data_acquisition} that randomly selects $h\approx rn$, $r \in (0,1)$, entries from $\bm f\in\R^n$ to be of zero value. In our experiments we choose $r=0.1$, $0.3$, and $0.5$. While the measurement data, $\bm y$, is acquired according to \cref{eq: data_acquisition} using the undersampling operator $A=HI_n$, during the recovery step in  \cref{eq:lasso_regression}, we assume the forward operator is the identity $A=I_n$, effectively treating the reconstruction as a denoising problem where the undersampling artifacts are part of the noise term. That is, we do not assume that we know which entries from ${\bm f}$ are zero.
Moreover, each undersampled measurement is perturbed by 20 dB of additive Gaussian noise. The visual distortion exhibited in \Cref{fig:undersampling} for the measurement data $\bm y$ is mainly attributable to undersampling (as opposed to the additive noise). 

\Cref{fig:undersampling}(top) compares the solutions for \cref{ex:example_denoising} obtained by \cref{eq:lasso_regression}, contrasting our new residual transform operator for regularization with the more standard first and third order differencing operators. The primary artifact in solutions $\bm x_{loc}$ obtained with the local differencing operator, $\LL=T^0_n$, is a pronounced ``staircasing'' effect, which is present across all tested undersampling ratios. Notably, this artifact persists even for a small undersampling ratio of $r=0.1$ (115 sampled points), as shown in the leftmost column. Moreover, within the  interval $(-\tfrac\pi2,\tfrac\pi2)$, a region of high variability in the true signal, the solution $\bm x_{loc}$ exhibits a significant loss of both magnitudes and variability. By contrast, the solution $\bm x_R$  obtained using our proposed residual operator of the same order ($\LL=R_n^0$) successfully captures the variability in smooth regions although there is a clear loss in signal magnitude. As before, since a higher-order regularizer is better suited for the piecewise polynomial function $f_1(s)$, we also show the results using third order differencing in the rightmost column.   In this case we see that the solution $\bm x_{loc}$ becomes comparable to $\bm x_R$ both in the smooth regions and near a discontinuity. However, as already pointed out in our previous examples, such information regarding the signal's variability within its smooth regions is typically not known a priori.
\Cref{fig:undersampling}(bottom) displays the spatial errors for the corresponding approximations in the top row,  confirming some advantage in using our residual transform operator, especially in smooth regions. It is clear, however, that the main benefit of using the residual transform operator in this experiment is that the type of variability can be captured.

\subsection{Recovery of a piecewise constant function}

\begin{figure}[h!]
    \centering
\begin{subfigure}[b]{.32\textwidth}
\includegraphics[width=\textwidth]{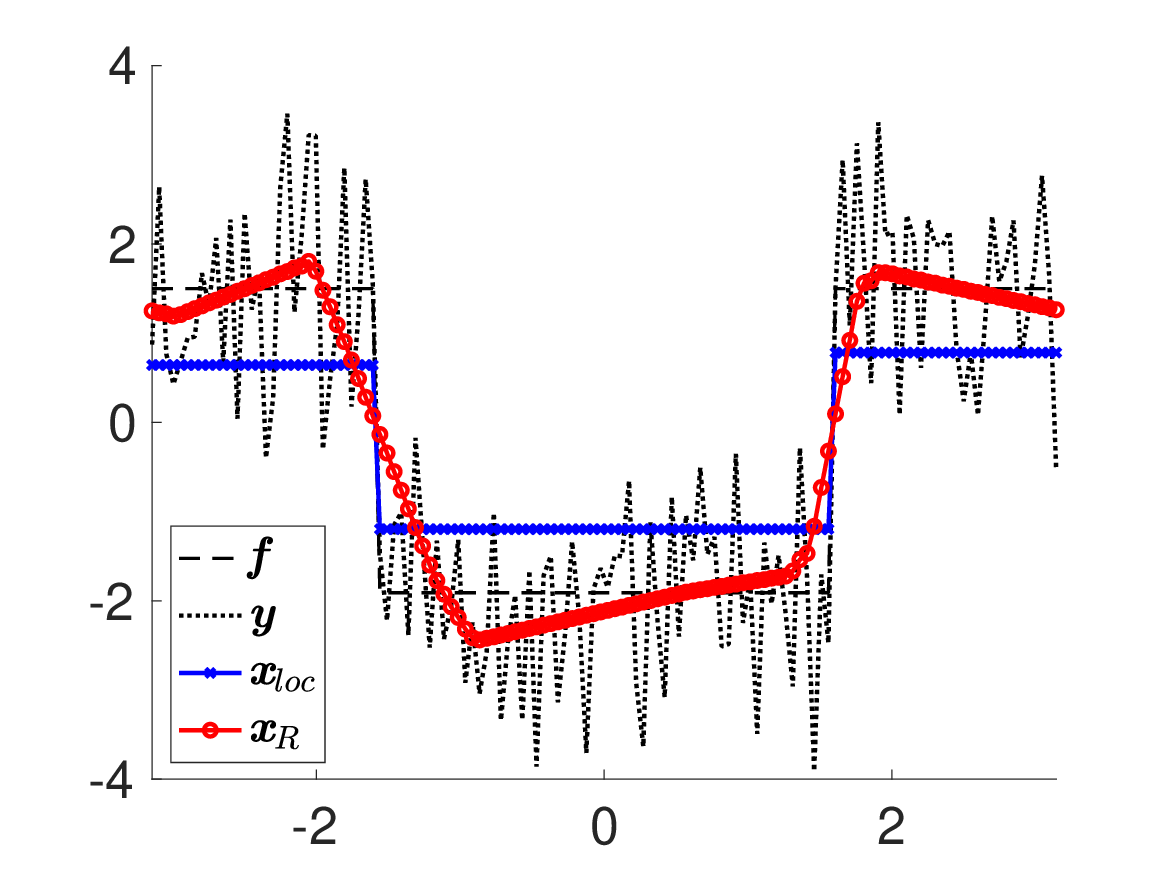}
\end{subfigure}
\begin{subfigure}[b]{.32\textwidth}
\includegraphics[width=\textwidth]{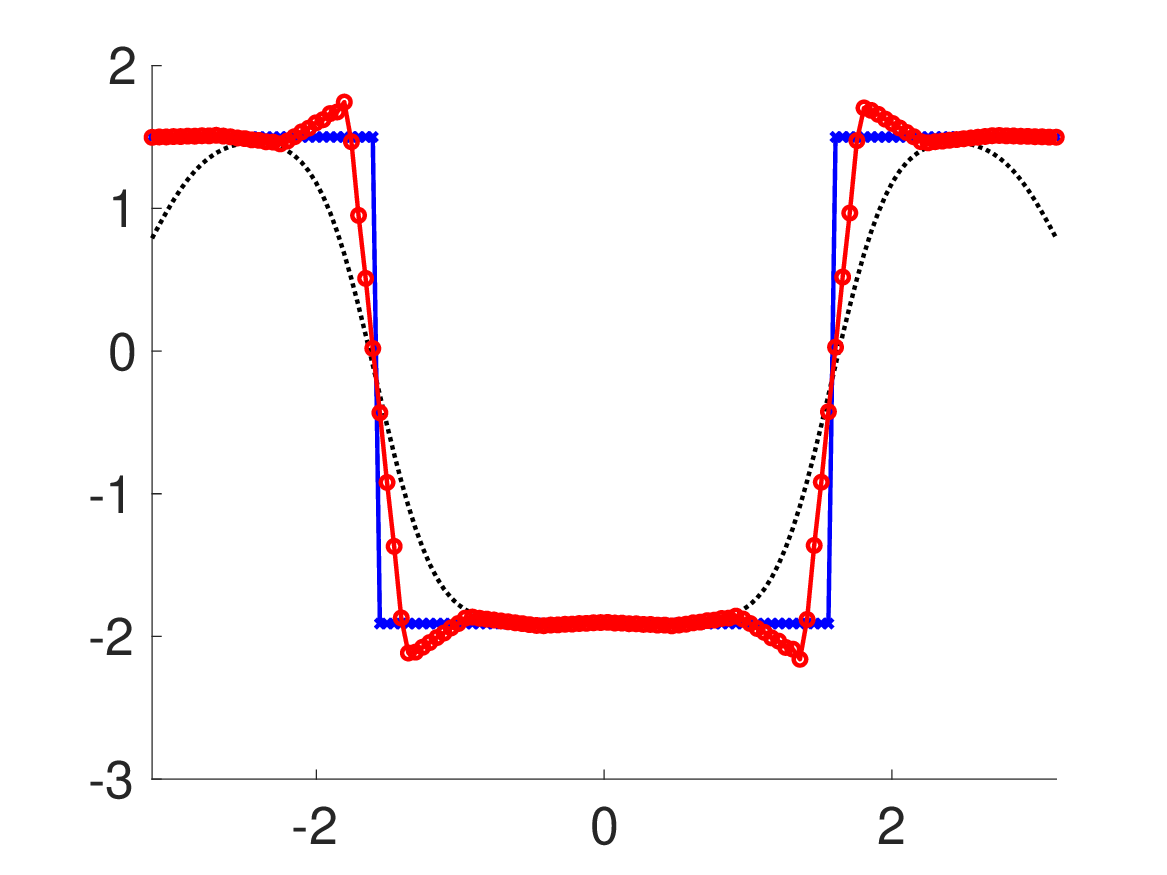}
\end{subfigure}
\begin{subfigure}[b]{.32\textwidth}
\includegraphics[width=\textwidth]{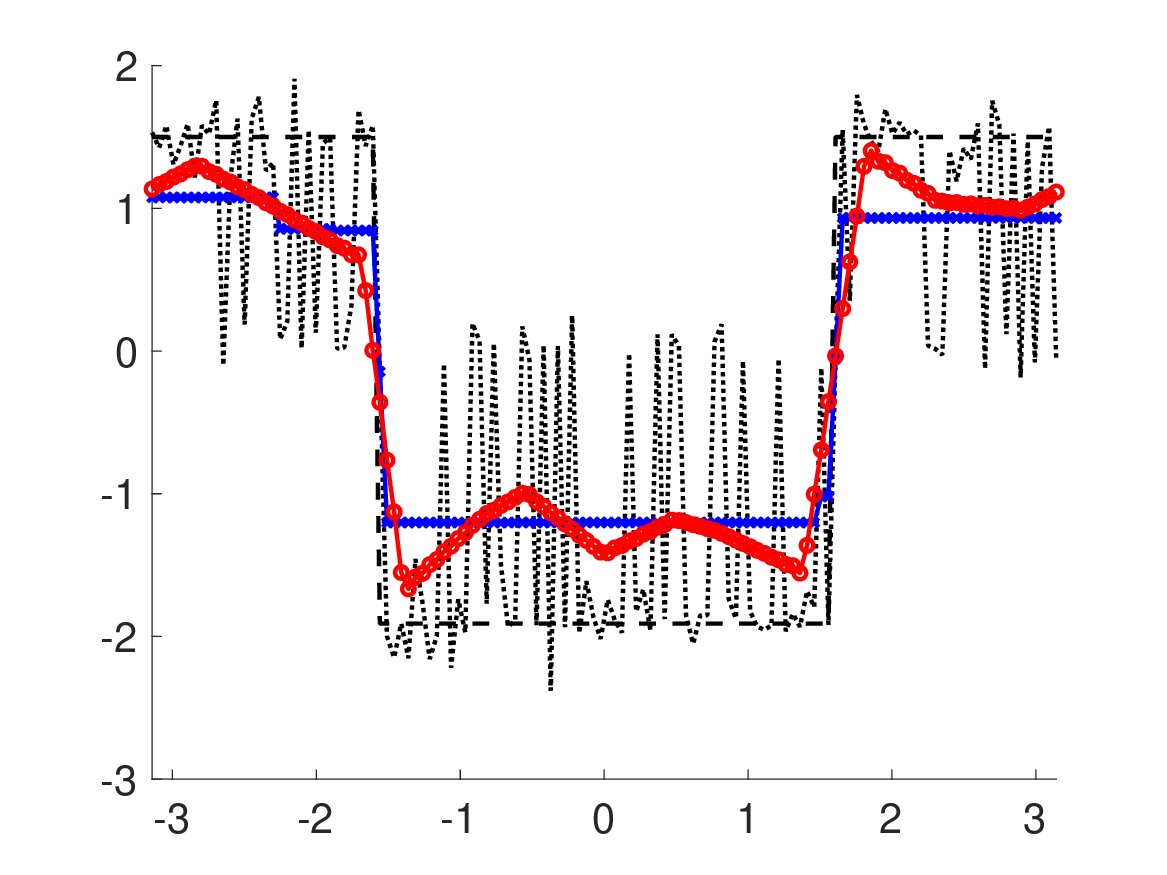}
\end{subfigure}\\
\begin{subfigure}[b]{.32\textwidth}
\includegraphics[width=\textwidth]{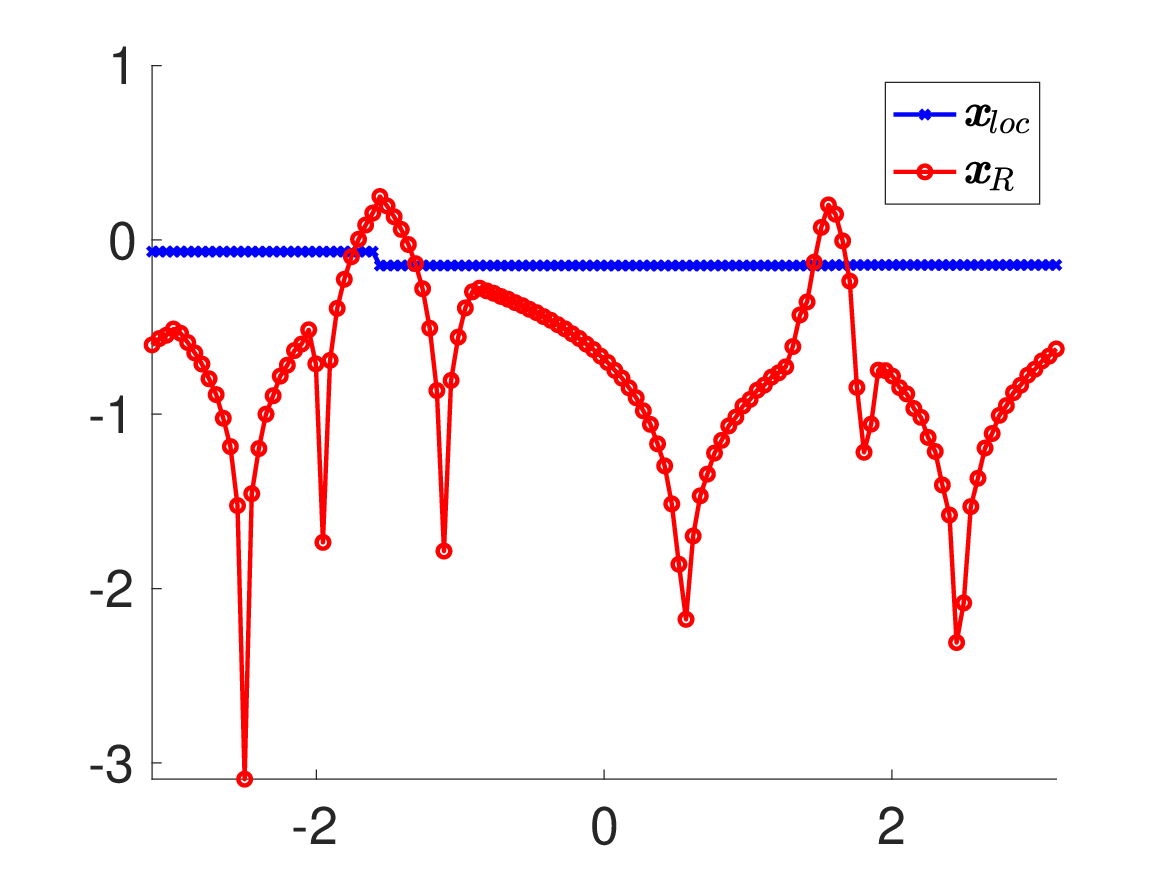}
\end{subfigure}
\begin{subfigure}[b]{.32\textwidth}
\includegraphics[width=\textwidth]{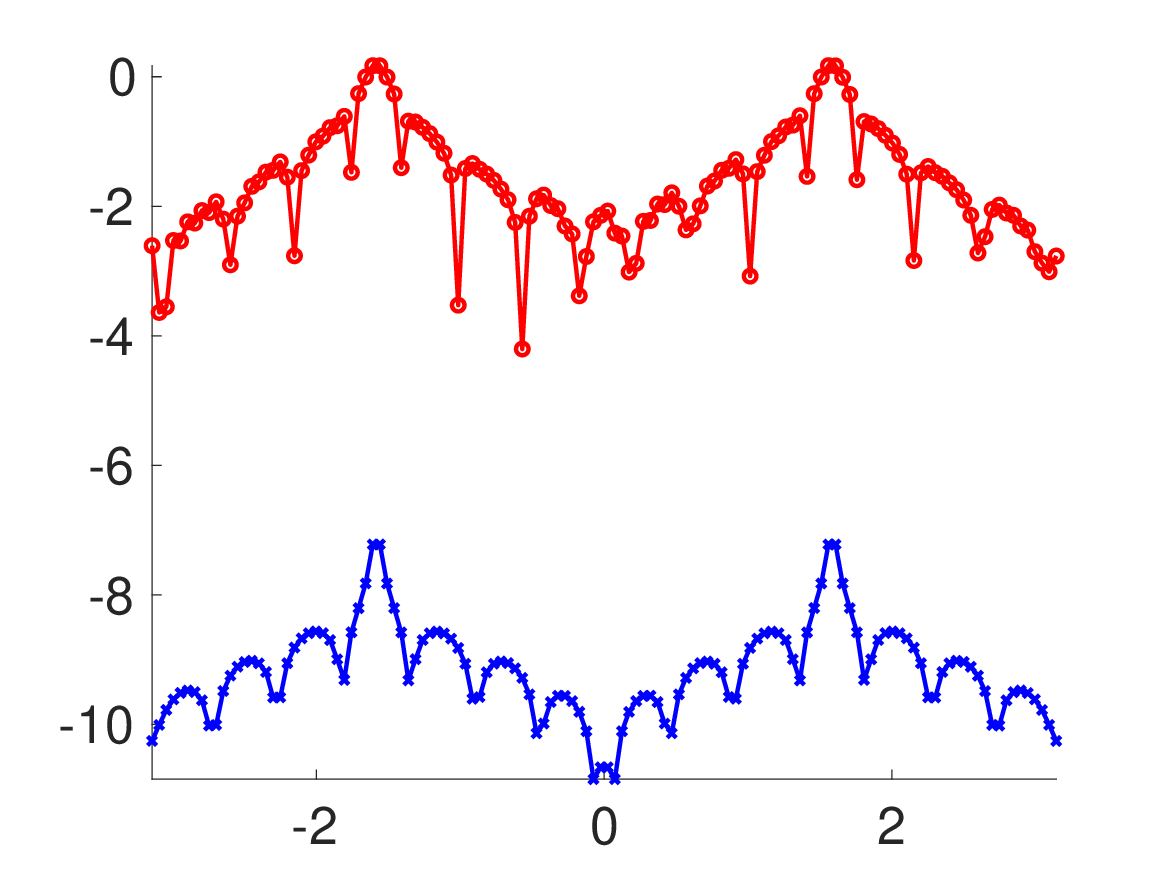}
\end{subfigure}
\begin{subfigure}[b]{.32\textwidth}
\includegraphics[width=\textwidth]{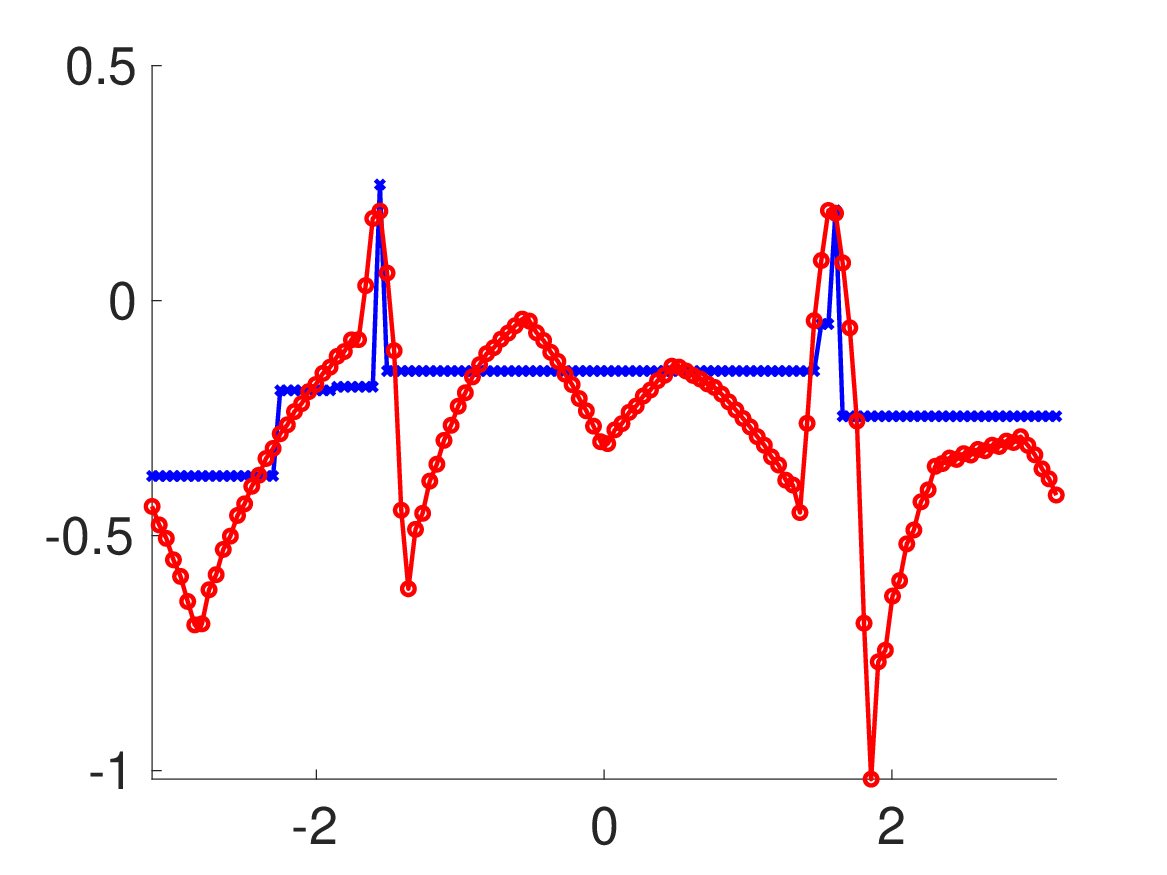}
\end{subfigure}
\caption{(Top row) Recovery of  \cref{ex:example2} for (left) denoising $\text{SNR}=5$ dB; (middle) deblurring for $\gamma=0.05$; (right) undersampling  with ratio $r=0.3$ and 20 dB additive Gaussian noise. All solutions use $p = 0$ and adopt regularization parameter $\alpha$ as in \cref{fig:denoising,fig:deblurring,fig:undersampling}. (Bottom row) Corresponding $\log_{10}E^{abs}_j$, $j = 1,\dots,n$,  of the top row solutions given by \cref{eq:absolute_err}.}
\label{fig:const_test}
\end{figure}

\begin{figure}[h!]
    \centering
\begin{subfigure}[b]{.32\textwidth}
\includegraphics[width=\textwidth]{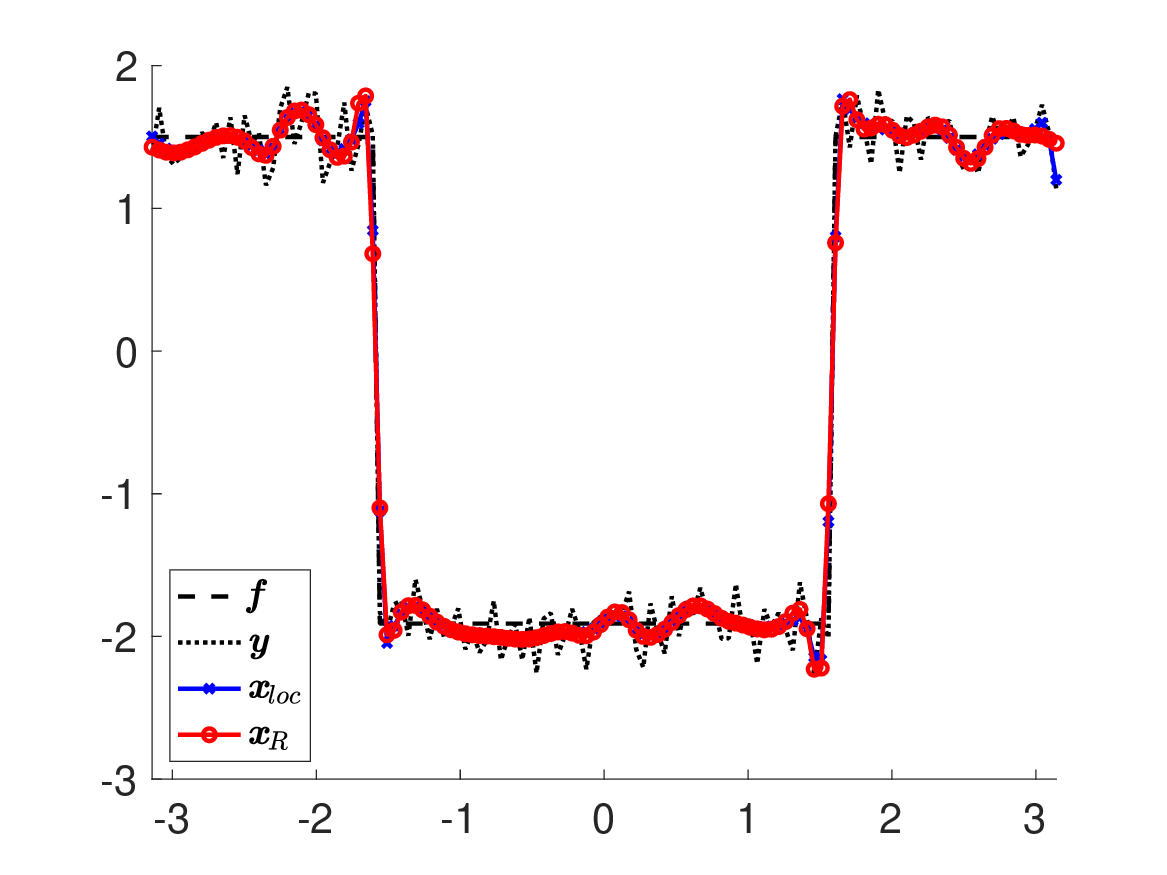}
\end{subfigure}
\begin{subfigure}[b]{.32\textwidth}
\includegraphics[width=\textwidth]{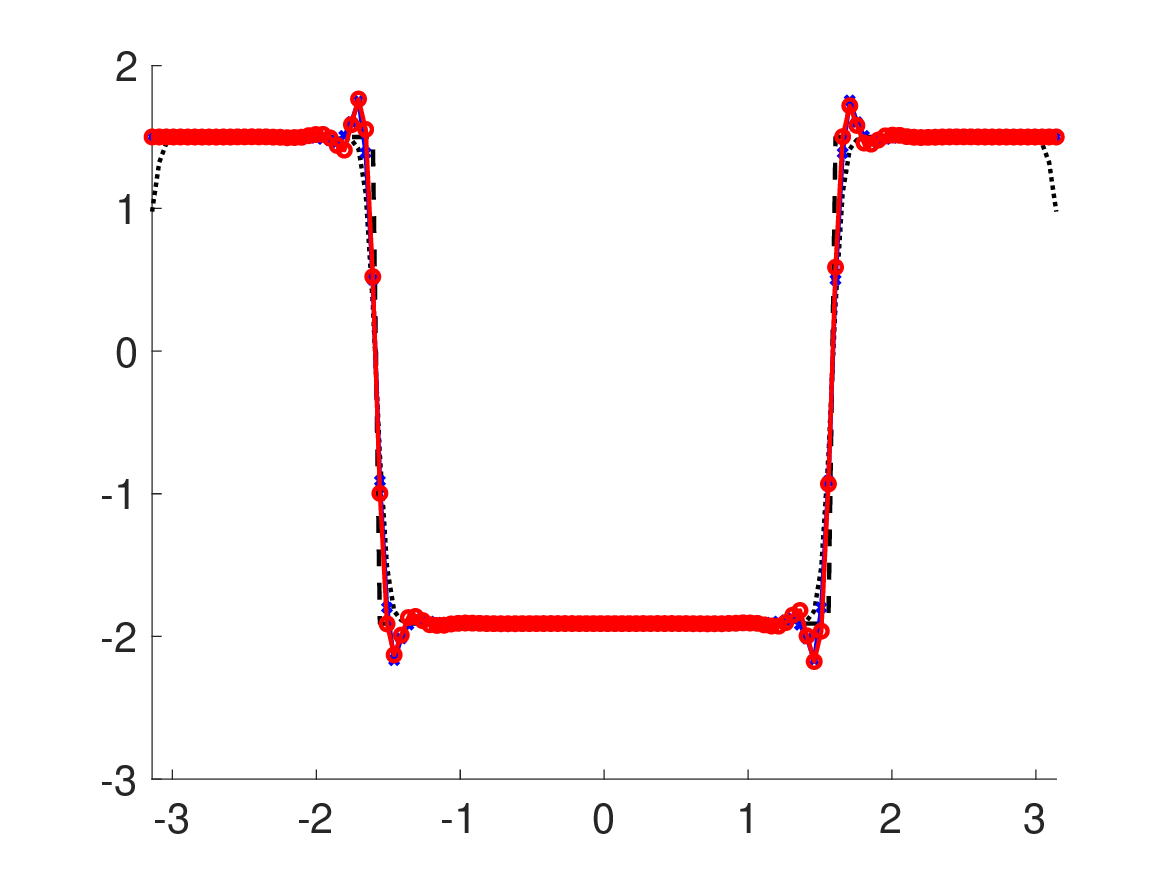}
\end{subfigure}
\begin{subfigure}[b]{.32\textwidth}
\includegraphics[width=\textwidth]{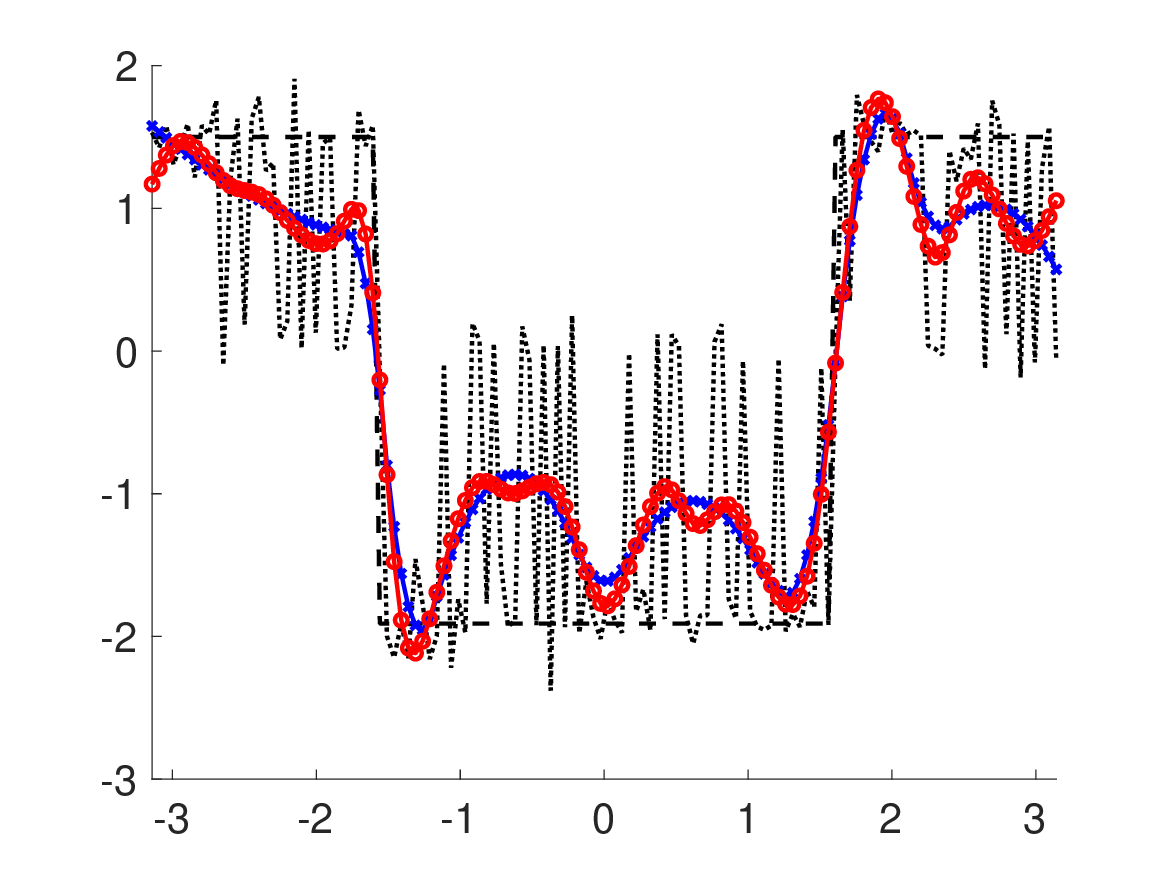}
\end{subfigure}\\
\begin{subfigure}[b]{.32\textwidth}
\includegraphics[width=\textwidth]{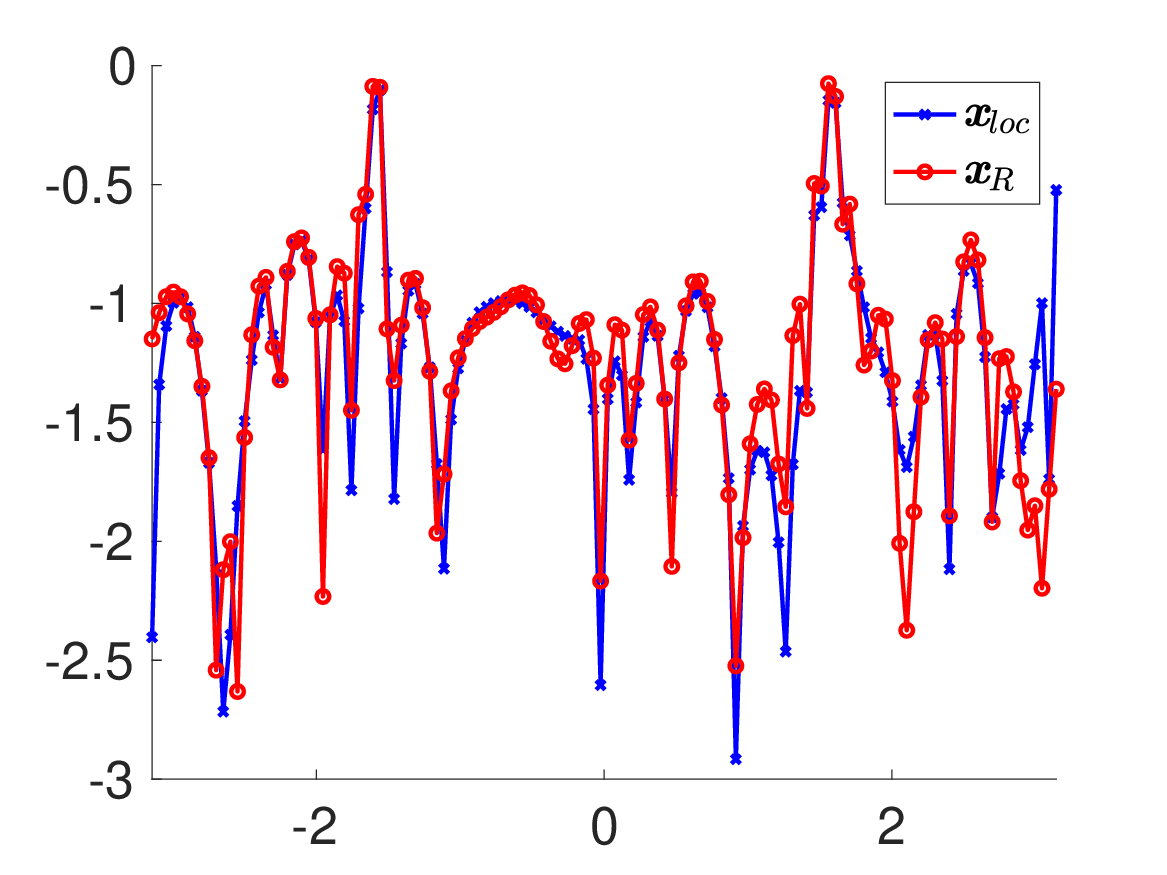}
\end{subfigure}
\begin{subfigure}[b]{.32\textwidth}
\includegraphics[width=\textwidth]{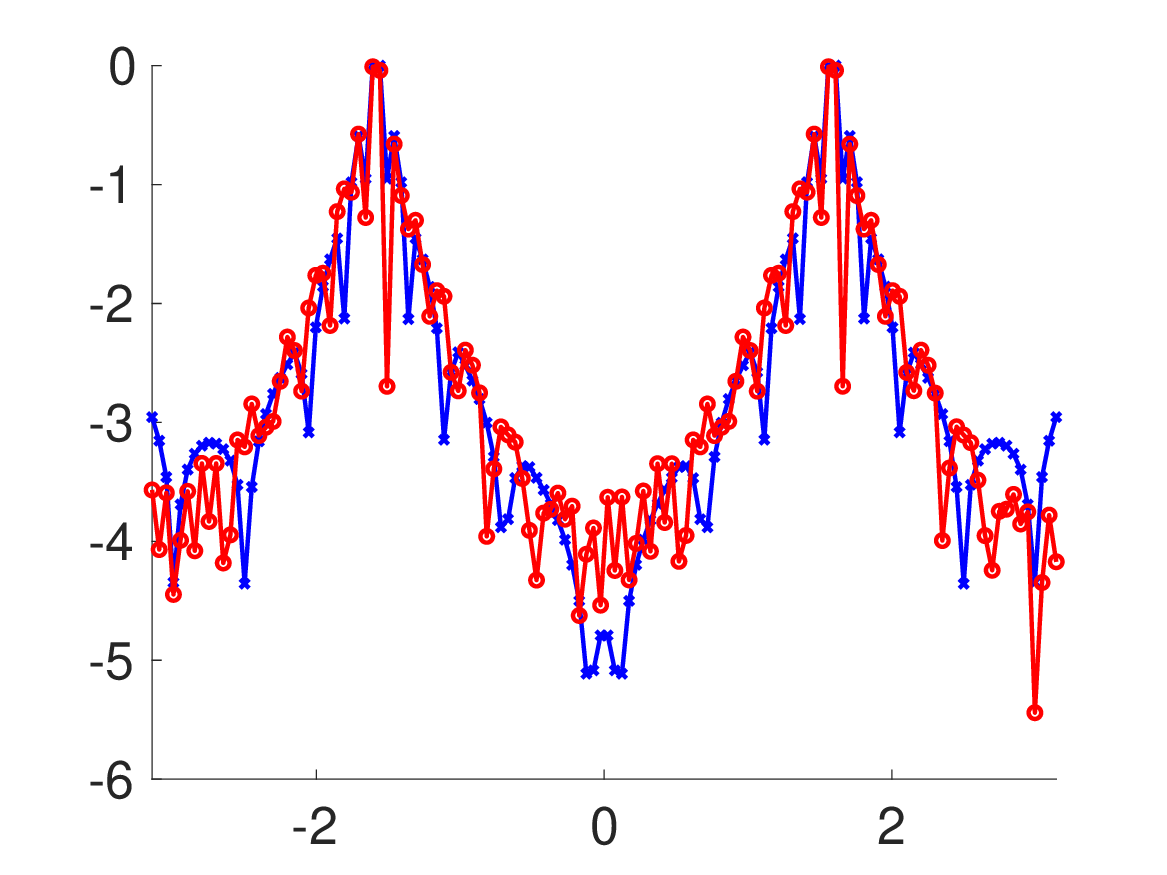}
\end{subfigure}
\begin{subfigure}[b]{.32\textwidth}
\includegraphics[width=\textwidth]{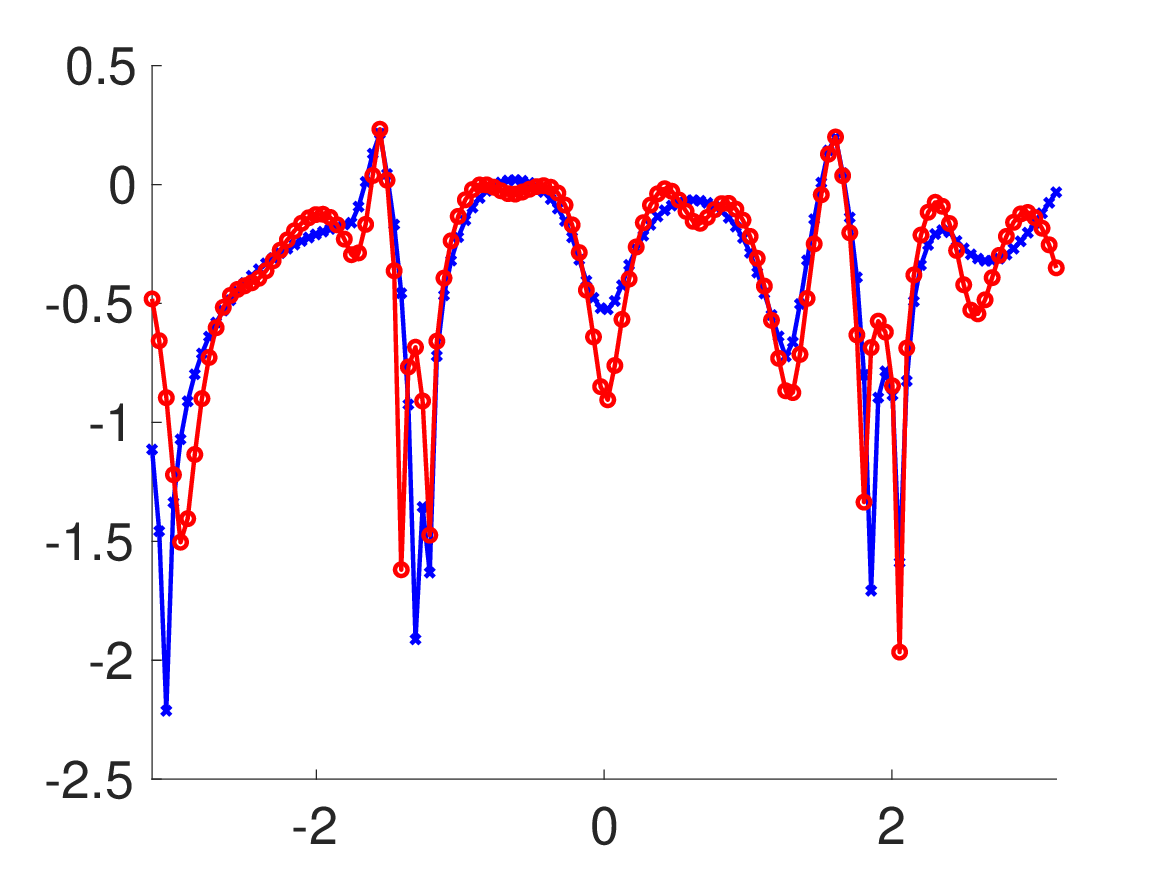}
\end{subfigure}
\caption{(Top row) Recovery of  \cref{ex:example2} for (left) denoising SNR=20 dB; (middle) deblurring for $\gamma=0.01$; (right) undersampling  with ratio $r=0.3$ and 20 dB additive Gaussian noise
. All solutions use $p = 1$ and adopt regularization parameter $\alpha$ as in \cref{fig:denoising,fig:deblurring,fig:undersampling}. (Bottom row) Corresponding $\log_{10}E^{abs}_j$, $j = 1,\dots,n$,  of the top row solutions given by \cref{eq:absolute_err}.}
\label{fig:const_test_p1}
\end{figure}
We include numerical experiments for \cref{ex:example2}  to demonstrate that while the TV operator is most suitable for the recovery of discretized piecewise constant functions, the proposed residual transform operator still yields reasonable results, a testament to the fact that when the type of variability is unknown, using $\mathcal L = R_{n,\zeta}^p$ is a viable and robust  option.  \Cref{fig:const_test} shows results corresponding to challenging cases for (left) additive noise  (middle) blurred data  and (right) undersampled data.   While using $\mathcal L = T_{n}^0$ in \cref{eq:lasso_regression} clearly yields more accurate results in the blurring case, there are still visible benefits to using $R_{n,\zeta}^0$, especially in obtaining the correct magnitude of the signal. 
\cref{fig:const_test_p1} displays results corresponding to the same suite of tests in \cref{fig:const_test} using $\LL=T^1_n$ and $\LL=R_{n,\zeta}^1$. Unlike what is observed in \cref{fig:const_test} for $p=0$, using $\LL=T^1_n$ in \eqref{eq:lasso_regression} no longer demonstrated an advantage in the blurring case.

\section{Concluding remarks}
\label{sec:summary}
This paper proposes a new sparsity-promoting operator, called the {\em residual transform operator}, which is constructed from the residuals of two distinct sparsity-promoting operators, $\LL_1$ and $\LL_2$, each of which promotes sparsity for the same type of variability (e.g. piecewise constant or piecewise linear behavior).  In our investigation we considered $\LL_1 = T_n^p$, a local differencing operator of order $2p+1$, and $\LL_2 = S_{n,\zeta}^p$, a global edge detector of order $2p+1$.  The resulting residual transform operator is then $R_{n, \zeta}^p = T_n^p - S_{n,\zeta}^p$.  The observation that motivates the construction of $R_{n, \zeta}^p$ is that $T_n^p{\bm f} = S_{n,\zeta}^p{\bm f}$ when  $\zeta = \frac{1}{2}$ for ${\bm f} = \{f(s_j)\}_{j = 1}^n$, with $s_j$ uniformly spaced.  For conditioning purposes we choose $\zeta = \frac{1}{4}$.  Our numerical experiments demonstrate that our new residual transform operator can be effectively used as the sparsity promoting operator in $\ell_1$ regularization.  Specifically, it is able to reduce errors that depend on signal variability—errors which often arise when using either $\LL_1$ or $\LL_2$ individually. Notably, this improvement holds true even when neither component operator provides an ideal sparse representation of the underlying signal. Furthermore, the construction of our residual operator does not require prior information about the signal's variability, making it broadly applicable.

Our numerical experiments further demonstrate that our new residual transform operator gains visible benefits compared to a standard local differencing operator, particularly for piecewise smooth signals with complex internal variations. This enhanced performance is evident across all test cases — denoising, deblurring, and undersampling — with notable improvements in both smooth regions and near discontinuities. Moreover, its performance on piecewise constant functions confirms that our operator is a robust and viable choice for Lasso regression, especially when prior knowledge of the signal's variability is unavailable.

While this paper demonstrates the benefits of employing the  residual transform operator in a Lasso regression, there are many worthwhile extensions to investigate. A natural next step is to incorporate the operator into a Bayesian learning framework, which would provide not only robust uncertainty quantification but also potentially more representative point estimates. Furthermore, the approach can be expanded to address a challenging problem of multiple measurement vectors that exhibit substantial variabilities with only shared locations of discontinuities. Finally, since we confirmed that an approximation to the sparse transform domain of the underlying signal is not necessarily germane to the construction of the residual transform operator, it is possible to identify myriad choices for $\LL_1$ and $\LL_2$ that are more suitable for multi-dimensional signals (although both operators in this case can be adapted to two-dimensions, \cite{archibald2005polynomial,xiao2022sequential}. Such extensions to our method will be explored  in future work.

\section*{Acknowledgments}
This work was partially supported by the DOD (ONR MURI) grant \#N00014-20-1-2595 and  the DOE ASCR grant \#DE-SC0025555.

\bibliographystyle{siamplain}
\bibliography{references}
\end{document}

%% file: main_shared.tex
\usepackage{lipsum}
\usepackage{amsfonts}
\usepackage{graphicx}
\usepackage{epstopdf}
\usepackage{algorithmic}
\ifpdf
  \DeclareGraphicsExtensions{.eps,.pdf,.png,.jpg}
\else
  \DeclareGraphicsExtensions{.eps}
\fi

\newsiamremark{remark}{Remark}
\newsiamremark{hypothesis}{Hypothesis}
\crefname{hypothesis}{Hypothesis}{Hypotheses}
\newsiamthm{claim}{Claim}
\newsiamremark{fact}{Fact}
\crefname{fact}{Fact}{Facts}

\headers{residual transform operator}{Y.\ Xiao, A.\ Gelb and A.\ Viswanathan}

\title{A new sparsity promoting residual transform operator for Lasso regression \thanks{\today}} 

\author{ 
Yao Xiao\thanks{
Department of Mathematics, Dartmouth College, USA 
(\email{yao.xiao@dartmouth.com}, \orcid{0000-0002-0850-9624})
}
\and
Anne Gelb\thanks{
Department of Mathematics, Dartmouth College, USA 
(\email{anne.E.gelb@dartmouth.edu}, \orcid{0000-0002-9219-4572})
}
\and
Aditya Viswanathan \thanks{
Mathematics and Statistics, University of Michigan-Dearborn, USA 
(\email{adityavv@umich.edu}, \orcid{ })
}
}

\usepackage{amsopn}
\usepackage[normalem]{ulem}

\DeclareMathOperator*{\argmin}{arg\,min} 

\newsiamremark{example}{Example}

\usepackage{amsmath}
\allowdisplaybreaks
\usepackage{amssymb}
\usepackage{commath}
\usepackage{mathtools}
\usepackage{bbm}
\usepackage{bm}

\usepackage{color}
\usepackage{graphicx}
\usepackage{subcaption}

\newcommand{\R}{\mathbb{R}} 
\newcommand{\LL}{\mathcal{L}}
